\documentclass[a4paper,10pt,]{article}

\usepackage[left=2cm, right=2cm, top=2.00cm, bottom=3cm]{geometry}

\usepackage{graphicx}
\usepackage{amssymb}
\usepackage{amsmath}
\usepackage{amsthm}
\usepackage{xcolor}
\usepackage{dsfont}
\usepackage{algpseudocode}
\usepackage{authblk}
\usepackage[ngerman, english]{babel}
\usepackage{makecell}

\usepackage{amsmath,bm}
\usepackage{enumerate}
\usepackage{algorithm} 
\usepackage{algpseudocode}
\usepackage{tabularx}
\usepackage{subcaption} 
\usepackage{booktabs} 
\usepackage{algpseudocode}
\usepackage{tikz}
\usetikzlibrary{decorations.pathreplacing}
\usepackage{caption}
\usepackage{subcaption}
\usepackage{bibentry}
\usepackage{colonequals}

\usepackage[colorlinks,
            pdffitwindow=false,
            plainpages=false,
            pdfpagelabels=true,
            pdfpagemode=UseOutlines,
            pdfpagelayout=SinglePage,
            bookmarks=false,
            colorlinks=true,
            hyperfootnotes=false,
            linkcolor=black,
            citecolor=black]
{hyperref}

\usepackage{mathtools}
\usepackage{enumerate}
\usepackage{thmtools}
\usepackage{cite}
\usepackage{float}

\allowdisplaybreaks

\usepackage{sectsty}
\sectionfont{\fontsize{10}{10}\selectfont}
\subsectionfont{\fontsize{10}{10}\selectfont}

\newcommand{\R}{\mathbb{R}}

\newcommand{\N}{\mathbb{N}}
\newcommand{\E}{\mathcal{E}}

\renewcommand{\H}{\mathcal{H}}
\newcommand{\G}{\mathcal{G}}

\renewcommand{\L}{\mathcal{L}}
\newcommand{\Pw}{\mathcal{P}_w}
\newcommand{\LPw}{\mathcal{L}\mathcal{P}_w}

\newcommand{\W}{\mathcal{W}}
\newcommand{\M}{\mathfrak{M}}

\DeclareMathOperator*{\argmin}{arg\,min}

\DeclareMathOperator{\conv}{conv}
\DeclareMathOperator{\proj}{proj}

\DeclareMathOperator{\dist}{dist}
\DeclareMathOperator{\closure}{cl}

\newcommand{\expnumber}[2]{{#1}\mathrm{e}{#2}}

\newtheorem{theorem}{Theorem}[section]
\newtheorem{myprop}[theorem]{Proposition}
\newtheorem{mydef}[theorem]{Definition}

\newtheorem{mylemma}[theorem]{Lemma}
\newtheorem{myremark}[theorem]{Remark}
\newtheorem{myexample}[theorem]{Example}

\numberwithin{equation}{section}

\makeatletter
\newcommand{\leqnomode}{\tagsleft@true\let\veqno\@@leqno}
\newcommand{\reqnomode}{\tagsleft@false\let\veqno\@@eqno}
\makeatother

\newcommand{\ignore}[1]{}
\newcommand{\nobibentry}[1]{{\let\nocite\ignore\bibentry{#1}}}

\let\OLDthebibliography\thebibliography
\renewcommand\thebibliography[1]{
	\OLDthebibliography{#1}
	\setlength{\parskip}{0pt}
	\setlength{\itemsep}{0pt plus 0.3ex}
}

\newenvironment{@abssec}[1]{
       \footnotesize
         {\upshape\bfseries #1. }\ignorespaces}
     {\if@twocolumn\else\par\vspace{.1in}\fi}
\newenvironment{keywords}{\begin{@abssec}{Key words}}{\end{@abssec}}
\newenvironment{AMS}{\begin{@abssec}{AMS subject classifications}}{\end{@abssec}}

\begin{document}

\title{Inertial dynamics with vanishing Tikhonov regularization for multiobjective optimization}
\author{Radu Ioan Bo\c{t}\thanks{Faculty of Mathematics, University of Vienna, Oskar-Morgenstern-Platz 1, 1090 Vienna, Austria, e-mail: \url{radu.bot@univie.ac.at}. Research partially supported by the Austrian Science Fund (FWF), projects W 1260 and P 34922-N.} \quad and \quad Konstantin Sonntag\thanks{Department of Mathematics, Paderborn University, Warburger Straße 100, 33098 Paderborn, Germany, e-mail: \url{konstantin.sonntag@upb.de}. Research partially support by the German Federal Ministry of Education and Research (BMBF) within the AI junior research group ``Multicriteria Machine Learning". Research initiated during a research visit at the University of Vienna in October 2023.}}

\date{\today}

\maketitle


\begin{abstract}
In this paper, we introduce, in a Hilbert space setting, a second order dynamical system with asymptotically vanishing damping and vanishing Tikhonov regularization that approaches a multiobjective optimization problem with convex and differentiable components of the objective function. Trajectory solutions are shown to exist in finite dimensions. We prove fast convergence of the function values, quantified in terms of a merit function. Based on the regime considered, we establish both weak and, in some cases, strong convergence of trajectory solutions towards a weak Pareto optimal point. To achieve this, we apply Tikhonov regularization individually to each component of the objective function. Furthermore, we conduct numerical experiments to validate the theoretical results and investigate the qualitative behavior of the dynamical system. This work extends results from convex single objective optimization into the multiobjective setting. The results presented in this paper lay the groundwork for the development of fast gradient and proximal point methods in multiobjective optimization, offering strong convergence guarantees.
\end{abstract}

{\begin{center}
\begin{minipage}[c]{0.898\textwidth}
\end{minipage}
\end{center}}



\section{Introduction}\label{sec:intro}
Let $\H$ be a real Hilbert space with inner product $\langle \cdot, \cdot \rangle$ and induced norm $\lVert \cdot \rVert$. Consider the problem 
\begin{align}
\label{eq:MOP}
\tag{MOP}
    \min_{x\in \H} F(x) \coloneqq \left[
    \begin{array}{c}
        f_1(x) \\
        \vdots \\
        f_m(x)
    \end{array}
    \right],
\end{align}
with $f_i:\H \to \R$, $i=1,\dots,m$, convex and continuously differentiable functions. In this paper we study the \emph{multiobjective Tikhonov regularized inertial gradient system} assigned to \eqref{eq:MOP} which is defined on $[t_0,+\infty)$ by
\begin{align}
    \label{eq:MTRIGS}
    \tag{MTRIGS}
    \frac{\alpha}{t^q} \dot{x}(t) + \proj_{C(x(t)) + \frac{\beta}{t^p} x(t) + \ddot{x}(t)}(0) = 0,
\end{align}
where $t_0 >0$, $\alpha,\, \beta > 0$ and $q \in (0,1],\, p \in (0,2]$ and $C(x) \coloneqq \conv \left( \{ \nabla f_i(x) \, : \, i = 1,\dots, m \} \right)$, with initial data $x(t_0) = x_0 \in \H$ and $\dot{x}(t_0) = v_0 \in \H$. Here, $\conv(\cdot)$ denotes the \emph{convex hull} of a set, and $\proj_K : \H \rightarrow \H,\, \proj_{K}(x) \coloneqq \argmin_{y \in K} \lVert y - x \rVert$, denotes the \emph{projection operator} onto a nonempty, convex and closed set $K \subseteq \H$. The development of the system \eqref{eq:MTRIGS} is motivated by the recent research on fast continuous gradient dynamics for single objective optimization problems with convex and differentiable objective functions.
In the latter case, namely, when $m = 1$ and $f:=f_1$ in \eqref{eq:MOP}, the system \eqref{eq:MTRIGS} reduces to the \emph{Tikhonov regularized inertial gradient system} 
\begin{align*}
    \label{eq:TRIGS_intro}
    \tag{TRIGS}
    \ddot{x}(t) + \frac{\alpha}{t^q} \dot{x}(t) + \nabla f(x(t)) + \frac{\beta}{t^p} x(t) = 0,
\end{align*}
which has recently been extensively studied in the literature (see \cite{Attouch2022, Attouch2021, Laszlo2023}). Assuming that $\argmin f$, the set of global minimizers of $f$, is not empty, if, for instance, $p \in (0,2)$, $q \in (0,1)$ and $p < q+1$, then for the trajectory solution $x(\cdot)$ of \eqref{eq:TRIGS_intro} it holds $f(x(t)) - \min f = \mathcal{O}\left(t^{-p}\right)$ as $t \to + \infty$, where $\min f$ denotes the minimal objective value of $f$. Thus, a convergence rate arbitrary close to $\mathcal{O}\left(t^{-2}\right)$ can be obtained. Additionally, the trajectory solution converges strongly to the element with the minimum norm in  $\argmin f$, that is, $x(t) \to \proj_{\argmin f} (0)$ as $t \to +\infty$. 

On the other hand, \eqref{eq:MTRIGS} is related to the \emph{multiobjective inertial gradient system with asymptotic vanishing damping}
\begin{align}
    \label{eq:MAVD_intro}
    \tag{MAVD}
    \frac{\alpha}{t} \dot{x}(t) + \proj_{C(x(t)) + \ddot{x}(t)}(0) = 0,
\end{align}
with $\alpha \ge 3$, which was introduced in \cite{Sonntag2024a} and further studied in \cite{Sonntag2024b}. The system \eqref{eq:MAVD_intro} builds on the \emph{inertial multiobjective gradient system}
\begin{align}
    \label{eq:IMOG'_intro}
    \tag{IMOG'}
    \gamma \dot{x}(t) + \proj_{C(x(t)) + \ddot{x}(t)}(0) = 0,
\end{align} 
with $\gamma > 0$, which has been examined in \cite{Sonntag2024a} and naturally extends the \emph{heavy ball with friction dynamical system} 
\begin{align}
    \label{eq:HBF_intro}
    \tag{HBF}
    \ddot{x}(t) + \gamma \dot{x}(t) + \nabla f(x(t)) = 0,
\end{align}
studied in \cite{Alvarez2000, Attouch2000, Polyak1964} in the context of single objective optimization. As shown in \cite{Sonntag2024a}, \eqref{eq:IMOG'_intro} has theoretical advantages over the dynamical system
\begin{align}
    \label{eq:IMOG_intro}
    \tag{IMOG}
    \ddot{x}(t) + \gamma \dot{x}(t) + \proj_{C(x(t))} = 0,
\end{align} 
which was introduced in \cite{Attouch2015} as the first multiobjective gradient-like dynamical system featuring an inertial term. As the asymptotic analysis of \eqref{eq:IMOG_intro} requires the condition $\gamma^2 \ge L$, where $L$ is a joint Lipschitz constant of the gradients of the components of the objective function, it is unclear whether \eqref{eq:IMOG_intro} can be adapted to systems with asymptotic vanishing damping, i.e., obtained by replacing $\gamma$ by $\frac{\alpha}{t}$. In \cite{Sonntag2024b}, it is shown that the \emph{merit function}
\begin{align}
    \label{eq:merit_function}
    \varphi: \H \to \R, \quad x \mapsto \varphi(x) \coloneqq \sup_{z \in \H} \min_{i=1,\dots, m}  f_i(x) - f_i(z),
\end{align}
exhibits fast convergence along the trajectory solutions of \eqref{eq:MAVD_intro}, namely,
 $\varphi(x(t)) = \mathcal{O}(t^{-2})$ as $t \to + \infty$, thus expressing fast convergence of the function values. In addition, for $\alpha > 3$, the trajectory solutions $x(\cdot)$ of \eqref{eq:MAVD_intro} weakly converge to a weak Pareto optimal points of \eqref{eq:MOP}. In the single objective case, when $m = 1$ and $f:=f_1$, the system \eqref{eq:MAVD_intro} reduces to the celebrated \emph{inertial gradient system with asymptotic vanishing damping}
\begin{align*}
    \label{eq:AVD_intro}
    \tag{AVD}
    \ddot{x}(t) + \frac{\alpha}{t} \dot{x}(t) + \nabla f(x(t)) = 0,
\end{align*}
which was introduced in \cite{Su2016} as the continuous counterpart of Nesterov's accelerated gradient method \cite{Nesterov1983}. The system  \eqref{eq:AVD_intro} has further been studied in several papers, including \cite{Cabot2009, Cabot2009_long, May2017, Attouch2018}. It holds that $f(x(t)) - \min f = \mathcal{O}(t^{-2})$ as $t \to +\infty$ and, for $\alpha > 3$, the trajectory solutions weakly converge to a global minimizer of $f$, provided that $\argmin f$ is not empty. Due to its convergence properties, \eqref{eq:MAVD_intro} is the natural counterpart of \eqref{eq:AVD_intro} when considering multiobjective optimization problems.

The dynamical system \eqref{eq:TRIGS_intro} enhances the asymptotic properties of \eqref{eq:AVD_intro} by ensuring, depending on the chosen regime, weak and even strong convergence of the trajectory to the minimum norm solution, while retaining the rapid convergence of function values. The dynamical system \eqref{eq:MTRIGS} we introduce in this paper aims to provide a similar improvement over \eqref{eq:MAVD_intro} in the context of multiobjective optimization. The main results regarding the asymptotic behavior \eqref{eq:MTRIGS} obtained in this paper are summarized in Table \ref{tab:Convergence_Results}. In principal, we obtain convergence rates for the function values which can be arbitrarily close to $\mathcal{O}(t^{-2})$ as $t \to + \infty$. Furthermore, for $p \in (0,2),\,q \in (0,1)$ and $p < q +1$ the trajectory solution $x(\cdot)$ converges strongly to a weak Pareto optimal solution which has the minimal norm in the set $\bigcap_{i=1}^m \L\left( f_i, f_i^{\infty} \right) \subseteq \Pw$, with $f_i^{\infty} \coloneqq \lim_{t \to + \infty} f_i(x(t))$, $\L\left( f_i, f_i^{\infty} \right)$ the lower level set of $f_i$ with respect to $f_i^{\infty}$ for $i = 1,\dots,m$, and $\Pw$ the set of weak Pareto optimal solutions of \eqref{eq:MOP}. For $p \in (0,2),\,q \in (0,1)$ and $p > q +1$, we show that the trajectory converges weakly to a weak Pareto optimal solution. The case $p = q + 1$ is critical, as it seems that convergence results for the trajectories cannot be obtained. In addition, we treat some boundary cases for the parameters $p$ and $q$, which require additional conditions on the parameters $\alpha$ and $\beta$.
\renewcommand{\arraystretch}{1.925}
\begin{table}[!h]
\begin{center}
\begin{tabular}{ |c|c|c|c|c|l| } 
    \hline
    \makecell{Conditions on \\ $p, q, \alpha, \beta$}
    & $\varphi(x(t))$ & $\lVert \dot{x}(t) \rVert$ & $\lVert x(t) - z (t) \rVert$ & $x(t)$  & \makecell[c]{Theorem}\\ 
    \hline
    $p \in (0,2], \, 2q < p$ & $\mathcal{O}\left( t^{-2q} \right)$ & $\mathcal{O}\left( t^{-q} \right)$  & $\mathcal{O}\left( 1 \right)$  & - & Thm. \ref{thm:main_convergence_special_r=q}\\ 
    \hline
    \makecell{$q \in(0,1),$ $p < q + 1$} & $\mathcal{O}\left( t^{-p} \right)$ & $\mathcal{O}\left( t^{\frac{\max(q, p - q) - (p + 1)}{2}} \right)$ & $\mathcal{O}\left( t^{\frac{\max(q, p - q) - 1}{2}} \right)$ & \makecell{strong \\ convergence} & \makecell[l]{Thm. \ref{thm:main_convergence_special_p<q+1},\\Thm. \ref{thm:strong_convergence}}\\ 
    \hline
    \makecell{$q = 1$, $\alpha \ge 3$} & $\mathcal{O}\left( t^{-p} \right)$ & $\mathcal{O}\left( t^{-\frac{p}{2}} \right)$  & $\mathcal{O}\left( 1 \right)$  & - & Thm. \ref{thm:main_convergence_special_p<q+1_r=1}\\ 
    \hline
    \makecell{$p \in (0,2),\,q + 1 < p$} & $\mathcal{O}\left( t^{-2q} \right)$ & \makecell{$\mathcal{O}\left( t^{-q} \right)$,\\$\int_{t_0}^{+\infty}s \lVert \dot{x}(s) \rVert^2 < + \infty$} & $\mathcal{O}\left( 1 \right)$ & \makecell{weak \\ convergence} & \makecell[l]{Thm. \ref{thm:main_convergence_special_r=q},\\Thm. \ref{thm:intergral_estimate},\\ Thm. \ref{thm:weak_convergence}}\\ 
    \hline
    \makecell{$q \in(0,1),$ $p = 2$, \\ $\beta \ge q(1 - q)$} & $\mathcal{O}\left( t^{-2q} \right)$ & \makecell{$\mathcal{O}\left( t^{-q} \right)$,\\$\int_{t_0}^{+\infty}s \lVert \dot{x}(s) \rVert^2 < + \infty$} & $\mathcal{O}\left( 1 \right)$ & \makecell{weak \\ convergence} & \makecell[l]{Thm. \ref{thm:main_convergence_special_r=q},\\Thm. \ref{thm:intergral_estimate2},\\ Thm. \ref{thm:weak_convergence}}\\
    \hline
\end{tabular}
\end{center}
\caption{Summary of main asymptotic results for \eqref{eq:MTRIGS}. The function $z(\cdot)$ is the generalized regularization path, that will be introduced in Section \ref{sec:tikh_moo}. The merit function $\varphi(\cdot)$ measures the decay of the function values and gets introduced in Subsection \ref{subsec:pareto optimality and merit functions}. All results have to be understood asymptotically, i.e., as $t \to + \infty$.}
\label{tab:Convergence_Results}
\end{table}
\renewcommand{\arraystretch}{1}

To this end, we extend the concept of Tikhonov regularization, initially developed in order to handle ill-posed integral equations in \cite{Tikhonov1963a, Tikhonov1963b}, to multiobjective optimization. The Tikhonov regularization of a convex optimization problem $$\min_{x \in \H} f(x)$$ reads
\begin{align*}
    \min_{x \in \H} f(x) + \frac{\varepsilon}{2}\lVert x \Vert^2,
\end{align*}
where $\varepsilon > 0$ is a positive constant.  Denoting for all $\varepsilon > 0$ its unique minimizer by
\begin{align*}
    x_{\varepsilon} \coloneqq \argmin_{x \in \H} \left \{f(x) + \frac{\varepsilon}{2}\lVert x \rVert^2 \right \}, 
\end{align*}
it holds that $x_{\varepsilon}$ converges strongly to  $\proj_{\argmin f} (0)$ as $\varepsilon \to 0$, given $\argmin f \neq \emptyset$. The set $\{ x_{\varepsilon} \,:\, \varepsilon > 0 \}$ forms a smooth curve called \emph{regularization path}. This is one of  the key ingredients used to prove the strong convergence of the trajectory solution of \eqref{eq:TRIGS_intro} to the element of minimum norm in $\argmin f$. To extend this approach to the multiobjective optimization setting, we need to define an appropriate generalization of the regularization path. Although there are a few studies addressing Tikhonov regularization in multiobjective optimization (see \cite{Chen2009, Chen2011, Chuong2010, Chuong2010b}), these works are limited to the finite dimensional case and impose stringent assumptions, such as the compactness of the set of weak Pareto optima. Furthermore, these studies do not address whether a Pareto optimum with the minimum norm is achieved and are thus not suitable for our convergence analysis.

Therefore, given a regularization function $\varepsilon(\cdot)$ and a solution $x(\cdot)$ to \eqref{eq:MTRIGS}, we define the \emph{generalized regularization path} for our problem as
\begin{align}
\label{eq:moo_reg_path_intro}
    z(t) \coloneqq \argmin_{z \in \H} \max_{i=1,\dots,m} f_i(z) - f_i(x(t)) + \frac{\varepsilon(t)}{2}\lVert z \Vert^2.
\end{align}
The optimization problem in \eqref{eq:moo_reg_path_intro} can be seen as a regularization of an adaptive Pascoletti-Serafini scalarization of \eqref{eq:MOP} (see \cite{Eichfelder2008}). It will turn out that $z(\cdot)$ strongly converges to the weak Pareto optimal point of \eqref{eq:MOP} with minimal norm in a particular lower level set of the objective function. This result will allow us to conclude that the trajectory solutions $x(\cdot)$ of \eqref{eq:MTRIGS} strongly converges to the same weak Pareto optimal point of \eqref{eq:MOP}. These investigations lay the groundwork for developing fast gradient and proximal point methods in multiobjective optimization with strong convergence guarantees for the iterates. This parallels recent advances in single objective optimization \cite{Laszlo2023, Laszlo2024, Laszlo2025a, Bot2024, Karapetyants2024, Laszlo2025}.

The paper is organized as follows. In the remainder of this section, we summarize the basic definitions of multiobjective optimization and introduce the standing assumptions for this study.  Section \ref{sec:tikh_moo} is dedicated to Tikhonov regularization. We discuss the single objective case, provide a brief overview of existing work for the multiobjective setting, and prove the strong convergence of the generalized regularization path to the weak Pareto optimal point of \eqref{eq:MOP} with minimal norm in a particular lower level set of the objective function. Section \ref{sec:MTRIGS}  formally introduces the system \eqref{eq:MTRIGS}, where we prove the existence of solutions in finite dimensions, discuss uniqueness, and gather preliminary results on the trajectories. Section \ref{sec:convergence_anaylsis} contains the asymptotic analysis of solutions of \eqref{eq:MTRIGS}. The main results of this section concern the fast convergence rate of the function values in terms of the merit function and the strong convergence of the trajectory solutions.  We conclude our work in Section \ref{sec:conclusion} and propose possible directions for future research.

\subsection{Pareto optimality and merit function}
\label{subsec:pareto optimality and merit functions}

The notions of optimality under consideration for the multiobjective optimization problem \eqref{eq:MOP} are introduced below.

\begin{mydef}
\label{def:Pareto_opt}
\begin{enumerate}[i)]
    \item An element $x^* \in \H$ is called \emph{Pareto optimal} for \eqref{eq:MOP} if there does not exist  $x \in \H$ such that $f_i(x) \le f_i(x^*)$ for all $i = 1,\dots,m$ and $f_j(x) < f_j(x^*)$ for at least one $j = 1,\dots,m$. The set of Pareto optimal points is called the \emph{Pareto set}, and will be denoted by $\mathcal{P}$. 
    \item An element $x^* \in \H$ is called \emph{weak Pareto optimal} if there does not exist $x \in \H$ such that $f_i(x) < f_i(x^*)$ for all $i = 1,\dots, m$. The set of all weak Pareto optimal points is called the \emph{weak Pareto set}, and will be denoted by $\Pw$.
\end{enumerate}
\end{mydef}
Obviously, every Pareto optimal element is weak Pareto optimal. The following definition extends the concept of a level set to vector valued functions.
\begin{mydef}
\label{def:level_set}
Let $F: \H \to \R^m,\, F(x) = (f_1(x), \dots, f_m(x))^{\top}$ be a vector valued function, and $a \in \R^m$.
\begin{enumerate}[i)]
    \item We define
    \begin{align*}
        \L(F, a) \coloneqq \left\lbrace x \in \H \, : \, F(x) \leqq a \right\rbrace = \bigcap_{i = 1}^m \left\lbrace x \in \H \, : \, f_i(x) \le a_i \right\rbrace,
    \end{align*}
where ``$\leqq$'' denotes the partial order on $\R^m$ induced by $\R^m_+$. For $a,b \in \R^m$ it holds $a \leqq b$ if and only if $a_i \le b_i$ for all $i = 1,\dots,m$.
    \item We denote
    \begin{align*}
        \LPw(F, a) \coloneqq \L(F,a) \cap \Pw.
    \end{align*}
\end{enumerate}
\end{mydef}

In addition to proving strong convergence for the trajectory solutions of \eqref{eq:MTRIGS}, we are interested in quantifying the speed of convergence in terms of the objective function values. In multiobjective optimization, a useful and meaningful notion used  for this purpose (see \cite{Sonntag2024a, Sonntag2024b, Tanabe2022, Tanabe2022_2, Tanabe2022_3, Tanabe2022_4, Yang2002, Liu2009}) is the merit function $\varphi: \H \to \R,\, x\mapsto \varphi(x) \coloneqq \sup_{z \in \H} \min_{i=1,\dots, m}  f_i(x) - f_i(z)$, see \eqref{eq:merit_function}.
The following result, given in \cite[Theorem 3.1]{Tanabe2022_3}, gives a complete description of the set of weak Pareto optimal points of \eqref{eq:MOP}.
\begin{theorem}
Let $\varphi(\cdot)$ be defined by \eqref{eq:merit_function}. For all $x \in \H$ it holds that $\varphi(x)\ge 0$. Furthermore, $x \in \H$ is a weak Pareto optimal element for \eqref{eq:MOP} if and only if $\varphi(x) = 0$.
\label{thm:varphi_properties}
\end{theorem}

Since $f_i$ is weakly lower semicontinuous for $i = 1,\dots,m$, the function $x \mapsto \min_{i=1,\dots, m}  f_i(x) - f_i(z)$ is weakly lower semicontinuous for every $z \in \H$ and therefore $\varphi( \cdot)$ is also weakly lower semicontinuous. This means that every weak accumulation point of a trajectory $x(\cdot)$ that satisfies $\lim_{t \to + \infty} \varphi(x(t)) = 0$ is weakly Pareto optimal. In the single objective case, i.e., for $m=1$ and $f_1:=f$, it holds $\varphi(x) = f(x) - \inf_{z \in \H} f(z)$ for all $x \in \H$. This provides another justification for using $\varphi(\cdot)$ as a measure of the convergence speed in multiobjective optimization. One should also note that, even if all objective functions are smooth, the function $\varphi(\cdot)$ is not smooth in general. The following lemma provides a useful characterization of $\varphi(\cdot)$.
\begin{mylemma}
    For $x_0 \in \H$ and $a \in \R_{+}^{m}$, assume that $\LPw(F, F(x)) \neq \emptyset$ holds for all $x \in \L(F, F(x_0) + a)$. Then, 
    \begin{align*}
        \varphi(x) = \sup_{z \in \LPw(F,F(x_0) + a)} \min_{i=1,\dots,m} f_i(x) - f_i(z) \quad \forall x \in \L(F, F(x_0) + a).
    \end{align*}
    \label{lem:sup_inf_u_0}
\end{mylemma}
\begin{proof}
    Let $x \in \L(F, F(x_0) + a)$ be fixed. Obviously,
    \begin{align}
    \label{eq:sup_phi_ineq1}
        \begin{split}
            \sup_{z \in \LPw(F,F(x_0) + a)} \min_{i=1,\dots,m} f_i(x) - f_i(z)
            \le \sup_{z \in \H} \min_{i=1,\dots,m} f_i(x) - f_i(z) = \varphi(x).
        \end{split}
    \end{align}
Next, we show that $\min_{i=1,\dots,m} f_i(x) - f_i(z) \le \sup_{z' \in \L(F, F(x))} \min_{i=1,\dots,m} f_i(x) - f_i(z')$ holds for all $z \in \H$.  We assume that there exists $z \not\in \L(F, F(x))$ with $\min_{i=1,\dots,m} f_i(x) - f_i(z) > \sup_{z' \in \L(F, F(x))} \min_{i=1,\dots,m} f_i(x) - f_i(z')$. Since $z \not\in \L(F, F(x))$, there exists $j \in \{1,\dots,m\}$ with $f_j(z) > f_j(x)$. Therefore
    \begin{align*}
        0 > \min_{i = 1,\dots,m} f_i(x) - f_i(z) \ge \sup_{z' \in \L(F, F(x))} \min_{i=1,\dots,m} f_i(x) - f_i(z') \ge 0,
    \end{align*}
    which leads to a contradiction. Hence, 
    \begin{align}
    \label{eq:sup_phi_ineq5}
        \sup_{z \in \H} \min_{i=1,\dots,m} f_i(x) - f_i(z) \le \sup_{z \in \L(F, F(x))} \min_{i=1,\dots,m} f_i(x) - f_i(z).
    \end{align}
Next, we show that $\sup_{z \in \L(F, F(x))} \min_{i=1,\dots,m} f_i(x) - f_i(z) \le \sup_{z \in \LPw(F, F(x))} \min_{i=1,\dots,m} f_i(x) - f_i(z)$. By assumption, for all $z \in \L(F, F(x))$ there exists $z' \in \LPw(F, F(z)) \subseteq \LPw(F, F(x))$. Since $z' \in \LPw(F, F(z)))$, it holds $f_i(z') \le f_i(z)$ for all $i = 1,\dots,m$, hence
    \begin{align}
    \label{eq:sup_phi_ineq6}
        \min_{i=1,\dots,m} f_i(x) - f_i(z) \le \min_{i=1,\dots,m} f_i(x) - f_i(z'). 
    \end{align}
    From \eqref{eq:sup_phi_ineq6}, we conclude
    \begin{align}
    \label{eq:sup_phi_ineq7}
        \sup_{z \in \L(F, F(x))} \min_{i=1,\dots,m} f_i(x) - f_i(z) \le \sup_{z \in \LPw(F, F(x))} \min_{i=1,\dots,m} f_i(x) - f_i(z).
    \end{align}
    Since $x \in \L(F, F(x_0) + a)$, we have $ \LPw(F,F(x)) \subseteq \LPw(F,F(x_0) + a)$, hence
    \begin{align}
    \label{eq:sup_phi_ineq3}
        \sup_{z \in \LPw(F,F(x))} \min_{i=1,\dots,m} f_i(x) - f_i(z) \le \sup_{z \in \LPw(F,F(x_0) + a)} \min_{i=1,\dots,m} f_i(x) - f_i(z).
    \end{align}
    Combining \eqref{eq:sup_phi_ineq5}, \eqref{eq:sup_phi_ineq7} and \eqref{eq:sup_phi_ineq3}, it yields 
    \begin{align}
    \label{eq:sup_phi_ineq8}
        \varphi(x) \le \sup_{z \in \LPw(F,F(x_0) + a)} \min_{i=1,\dots,m} f_i(x) - f_i(z),
    \end{align}
which proves the statement.
\end{proof}

\subsection{Assumptions}

The research presented in this paper is conducted within the context of the following standing assumptions, which apply throughout the paper.

\begin{enumerate}
    \renewcommand{\labelenumi}{\textbf{\theenumi}}
    \renewcommand{\theenumi}{$(\mathcal{A}_{\arabic{enumi}})$}
    \item \label{ass:A1}The component functions $f_i:\H \to \R,\, i=1, \dots, m,$ are convex and continuously differentiable with Lipschitz continuous gradients.
    \item \label{ass:A3} Given the initial data $t_0 > 0$ and $x_0, v_0 \in \H$, define $a \in \R^m$ with $a_i \coloneqq \frac{\beta}{2t_0^p}\lVert x_0 \rVert^2 + \frac{1}{2}\lVert v_0 \rVert^2$ for $i = 1, \dots, m$. For all $x \in \L(F, F(x_0) + a)$ it holds that $\LPw(F, F(x)) \neq \emptyset$ and further
    \begin{align}
        R \coloneqq  \sup_{F^* \in F(\LPw(F,F(x_0) + a))} \inf_{z \in F^{-1}(\{F^*\})} \lVert z \rVert < + \infty.
    \label{eq:assumption_par_front}
    \end{align}
    \item \label{ass:A4}The set $S(q) \coloneqq \argmin_{z \in \H} \max_{i=1,\dots,m}  f_i(z) - q_i \neq \emptyset$ is nonempty for all $q \in \R^m$ and the mapping $z_0:\R^m \to \H,\, q \mapsto \proj_{S(q)}(0)$, is continuous.
\end{enumerate}

\subsubsection{Discussion of assumption \ref{ass:A3}}

The assumption \ref{ass:A3} is in the spirit of a hypothesis used in the literature (see \cite{Sonntag2024a, Sonntag2024b, Tanabe2022, Tanabe2022_2, Tanabe2022_3, Tanabe2022_4}) in the asymptotic analysis of continuous and discrete time gradient methods for multiobjective optimization. There, the assumption is formulated only for $a = 0$, which is recovered in our setting if we restrict the initial conditions to $x_0 = v_0 = 0$. For arbitrary initial conditions, our analysis requires the assumption to hold for $a \in \R_{+}^{m}$ by $a_i := \frac{\beta}{2t_0^p}\lVert x(t_0) \rVert + \frac{1}{2}\lVert \dot{x}(t_0) \rVert^2 \ge 0$ for $i = 1,\dots,m$,  as for this choice of $a$, the solutions of \eqref{eq:MTRIGS} can be shown to remain in $\L(F, F(x(t_0)) + a)$. This expansion of the level set is necessary because of the additional Tikhonov regularization which can produce trajectories that leave the initial level set $\L(F, F(x(t_0)))$. We visualize \ref{ass:A3} in Figure \ref{fig:ass:A3}, which shows the schematic image space for an \eqref{eq:MOP} with two objective functions. Given an initial point $x_0 \in \H$ and $a \in \R^m$ from \ref{ass:A3}, the set $F(\LPw(F(x_0) + a))$ is shown in blue. For all function values $F^* \in F(\LPw(F(x_0) + a))$ the constant $R$ gives a uniform bound on the minimum norm element in the preimage $F^{-1}(\{F^*\})$. For the single objective case ($m = 1$) this assumption is naturally satisfied if a solution to the optimization problem exists.

\begin{figure}
    \centering


    \begin{tikzpicture}[scale=1.2]
        
        \draw[black, thick, -stealth] (-.5,-.25) -- (4,-.25);
        \draw[black, thick, -stealth] (0,-.75) -- (0,4);

        \node at (3.66, -.58) {$f_1$};
        \node at (-.33, 3.66) {$f_2$};

        \draw [white, fill=gray!40] plot [smooth, tension=0] coordinates { (1,1)(1, 3.75) (3.75, 3.75) (3.75, 1)};
        \draw [black, thick, fill=gray!40] plot [smooth, tension=2] coordinates { (1,3.75) (1,1) (3.75,1)};

        \node at (3.35,3.5)  {$F(\H)$};

        \draw [black, line width=0.25mm] plot [smooth, tension=0.95] coordinates { (2.5,1.5) (2.7,1.75) (2.5,1.85) (2, 1.55) (1.9, 1.7) (1.6, 1.5) (1.45, 1.6) (1.3, 1.2) (1.15, 1.3) (1, 1)};
        \node at (1.15, 1.775)  {$F(x(t))$};  

        \fill (2.5,1.5) circle (0.05);
        \node at (2.4,1.25)  {$F(x_0)$};        
        \fill (3,2) circle (0.05);
        \node at (2.9,2.25)  {$F(x_0) + a$};

        \draw[black, densely dotted] (3,2) -- (0,2);
        \draw[black, thick] (.05,2) -- (-.05,2);
        \draw[black, densely dotted] (3,2) -- (3,-.25);
        \draw[black, thick] (3,-.2) -- (3,-.3);

        \begin{scope}
        \clip (0,0) rectangle (2.4 ,2);
        \draw [blue, line width=0.35mm] plot [smooth, tension=2] coordinates { (1,3.75) (1,1) (3.75,1)};
        \end{scope}

        \node[blue] at (1.5,0.05)  {$F(\LPw(F, F(x_0) + a))$};
        
    \end{tikzpicture}

    \caption{Visualization of \ref{ass:A3} with a trajectory $x(t) \in \LPw(F, F(x_0) + a)$.}
    \label{fig:ass:A3}
\end{figure}

\subsubsection{Discussion of assumption \ref{ass:A4}}

We need assumption \ref{ass:A4} to show the strong convergence of the generalized regularization path for multiobjective optimization problems.  We illustrate the necessity of this assumption with an example in Section \ref{sec:tikh_moo}. In the following we show that the continuity of the projection $q \mapsto z_0(q) \coloneqq \proj_{S(q)}(0)$ is closely connected with the continuity of the set-valued map (see \cite{Beer1988, Beer1992, Terazono2015, Mosco1969, Borwein1989, Aubin2009} for related discussions)
\begin{align*}
    S:\R^m \rightrightarrows \H, \quad q \mapsto S(q) \coloneqq \argmin_{z \in \H} \max_{i=1,\dots,m} f_i(z) - q_i.
\end{align*}
To this end, we recall the notion of Mosco convergence (see \cite{Beer1988}).

\begin{mydef}
\label{def:mosco_convergence}
    Let $\{C^k\}_{k \geq 0}, C^* \subseteq \H$ be nonempty, convex and closed sets. We say that the sequence $\{C^k\}_{k \geq 0}$ is Mosco convergent to $C^*$ if
    \begin{enumerate}[i)]
        \item for any $x^* \in C^*$ there exists $\lbrace x^k\rbrace_{k \ge 0}$ with $x^k \to x^*$ such that $x^k \in C^k$ for all $k \ge 0$;
        \item for any sequence $\lbrace k_l\rbrace_{l \ge 0} \subseteq \N$ with $x^{k_l} \in C^{k_l}$ for all $l \ge 0$ such that $x^{k_l} \rightharpoonup x^{*}$ as $l \to +\infty$, it holds $x^{*} \in C^*$.
    \end{enumerate}
\end{mydef}

Here we use $\to$ to denote strong convergence and $\rightharpoonup$ to denote weak convergence. The following theorem can be used to derive the continuity of $z_0(\cdot)$ from the \emph{Mosco continuity} of $S(\cdot)$. We recall that a set-valued map $S(\cdot)$ is said to be Mosco continuous if for all $q^* \in \R^m$ and any sequence $\lbrace q^k \rbrace_{k \ge 0} \subseteq \R^m$ with $q^k \to q^*$ the sequence $\{S(q^k)\}_{k \geq 0}$ is Mosco convergent to $S(q^*)$.

\begin{theorem}(\cite[Sonntag-Attouch Theorem]{Beer1988})
    Let $\{C^k\}_{k \geq 0}, C^* \subseteq \H$ be nonempty, convex and closed sets. The following statements are equivalent:
    \begin{enumerate}[i)]
        \item $\{C^k\}_{k \geq 0}$ is Mosco convergent to $C^*$;
        \item $\{C^k\}_{k \geq 0}$ is Wijsman convergent to $C^*$, i.e., for all $x \in \H$, it holds $\lim_{k \to + \infty} \dist(x, C^k) = \dist(x, C^*) $;
        \item for all $x \in \H$, it holds $\lim_{k \to + \infty} \proj_{C^k}(x) = \proj_{C^*}(x)$.
    \end{enumerate}
    \label{thm:sonntag_attouch}
\end{theorem}

The following proposition shows that for all $q^* \in \R^m$ and for any sequence $\lbrace q^k \rbrace_{k \ge 0} \subseteq \R^m$ with $q^k \to q^*$, condition \emph{ii)} in the definition of the Mosco convergence of $\{S(q^k)\}_{k \geq 0}$ to $S(q^*)$ is always fulfilled.

\begin{myprop}\label{prop17}
 Let $q^* \in \R^m$ and $\{ q^k \}_{k \ge 0} \subseteq \R^m$ be a sequence with $q^k \to q^*$ as $k \to + \infty$. Let $\{x^k\}_{k \ge 0} \subseteq \H$ be a sequence with $x^k \in S(q^k)$ for all $k \ge 0$ such that $x^k \rightharpoonup x^{*} \in \H$ as $k \to + \infty$. Then, $x^{*} \in S(q^*)$.
\end{myprop}

\begin{proof}
We show that
    \begin{align*}
        \max_{i=1,\dots,m}f_i(x^{*}) - q_i^* \le \max_{i=1,\dots,m}f_i(z) - q_i^* \quad \forall z \in \H.
    \end{align*}
    Let $z \in \H$ be arbitrary. We use the weak lower semicontinuity of $\max_{i=1,\dots,m} f_i(\cdot) - q_i^*$ to conclude
    \begin{align*}
        \max_{i = 1,\dots, m} f_i(x^{*}) - q_i^* \le & \liminf_{k \to +\infty} \max_{i = 1,\dots, m} f_i(x^k) - q_i^* \le \liminf_{k \to +\infty} \left(\max_{i = 1,\dots, m} f_i(x^k) - q_i^k + \max_{i = 1,\dots,m} q_i^k - q_i^* \right)\\
        = & \liminf_{k \to +\infty} \max_{i = 1,\dots, m} f_i(x^k) - q_i^k \le \liminf_{k \to +\infty} \max_{i = 1,\dots, m} f_i(z) - q_i^k \\
        \le & \liminf_{k \to +\infty} \left(\max_{i = 1,\dots, m} f_i(z) - q_i^* + \max_{i = 1,\dots,m} q_i^* - q_i^k \right) = \max_{i = 1,\dots, m} f_i(z) - q_i^*.
    \end{align*}
    Hence $x^{*} \in S(q^*)$, which completes the proof.
\end{proof}

The condition \emph{i)} in the definition of the Mosco convergence of $\{S(q^k)\}_{k \geq 0}$ to $S(q^*)$ when $q^k \to q^*$ as $k \to + \infty$ does not hold in general, but can be show to be satisfied under various circumstances. One of these is when the function $x \mapsto \max_{i = 1,\dots,m} f_i(x) - q_i$ exhibits a growth property uniformly for $q \in \R^m$ along approximating sequences.

\begin{mydef}(growth property uniformly along approximating sequences) Assume $S(q) \neq \emptyset$ for all $q \in \R^m$. We say that the function $x \mapsto \max_{i = 1,\dots,m} f_i(x) - q_i$ satisfies the growth property uniformly along approximating sequences if for all $q^* \in \R^m$ there exists a strictly increasing function $\psi:[0,+\infty) \to [0,+\infty)$ with $\psi(0) = 0$ such that for all sequences $\{ q^k \}_{k \ge 0} \subseteq \R^m$ with $q^k \to q^*$ as $k \to + \infty$ it holds
    \begin{align*}
        \max_{i=1,\dots,m} f_i(x^*) - q^k_i - \inf_{z \in \H} \max_{i=1,\dots,m} f_i(z) - q^k_i \ge \psi\left( \dist(x^*, S(q^k)) \right) \quad \forall x^* \in S(q^*) \ \forall k \geq 0.
    \end{align*}
\end{mydef}

The following lemma states the Lipschitz continuity of the optimal value function arising in the definition of the set-valued map $S(\cdot)$.

\begin{mylemma}
\label{lem:cont_opt_val_fun}
Assume $S(q) \neq \emptyset$ for all $q \in \R^m$. Then, the optimal value function
    \begin{align*}
        v: \R^m \to \R, \quad q \mapsto v(q) \coloneqq \inf_{z \in \H} \max_{i=1,\dots,m} f_i(z) - q_i,
    \end{align*}
is Lipschitz continuous. 
\end{mylemma}

\begin{proof}
    Let $q^1, q^2 \in \R^m$ and choose $x^1 \in S(q^1)$ and $x^2 \in S(q^2)$. It holds
    \begin{align*}
    \begin{split}
        v(q^1) = & \max_{i=1,\dots,m} f_i(x^1) - q_i^1 \le \max_{i=1,\dots,m} f_i(x^2) - q_i^1 \\
        \le & \max_{i=1,\dots,m} f_i(x^2) - q_i^2 + \max_{i = 1,\dots, m} q_i^2 - q_i^1 \le v(q^2) + \lVert q^1 - q^2 \Vert_{\infty}.
    \end{split}
    \end{align*}
Analogously,
    \begin{align*}
        v(q^2) \le v(q^1) + \lVert q^1 - q^2 \Vert_{\infty},
    \end{align*}
 thus,
    \begin{align*}
        \lvert v(q^1) - v(q^2) \rvert \le \lVert q^1 - q^2 \Vert_{\infty}.
    \end{align*}
\end{proof}

The next theorem shows that the uniform growth property indeed guarantees that for all $q^* \in \R^m$ and for any sequence $\lbrace q^k \rbrace_{k \ge 0} \subseteq \R^m$ with $q^k \to q^*$, the sequence $\{S(q^k)\}_{k \geq 0}$  is Mosco convergent to $S(q^*)$. Therefore, in the light of Theorem \ref{thm:sonntag_attouch}, assumption \ref{ass:A4} is fulfilled.

\begin{theorem}
\label{thm:uniform_growth_to_mosco}
Assume $S(q) \neq \emptyset$ for all $q \in \R^m$ and that $x \mapsto \max_{i = 1,\dots,m} f_i(x) - q_i$ satisfies the growth property uniformly along approximating sequences. Let $q^* \in \R^m$ and $\{q^k\}_{k \ge 0} \subseteq \R^m$ be a sequence with $q^k \to q^*$ as $k \to + \infty$. Then, $\{S(q^k)\}_{k \geq 0}$  is Mosco convergent to $S(q^*)$.
\end{theorem}
\begin{proof}
Condition \emph{ii)} in Definition \ref{def:mosco_convergence} is satisfied according to Proposition \ref{prop17}. We prove by contradiction that condition \emph{i)} is also satisfied. Let $x^* \in S(q^*)$ be such that for any sequence  $\{ x^k \}_{k \ge 0}$ with $x^k \in S(q^k)$ for all $k \geq 0$, it holds $x^k \not\to x^*$ as $k \to + \infty$. Hence, there exist $\delta >0$ and a subsequence $\{k_l\}_{l \ge 0} \subseteq \N$  such that $\dist(x^*, S(q^{k_l})) > \delta$ for all $l \ge 0$. We use the growth property to conclude
    \begin{align*}
        \max_{i=1,\dots,m} f_i(x^*) - q^{k_l}_i - \inf_{z \in \H} \max_{i=1,\dots,m} f_i(z) - q^{k_l}_i \ge \psi\left( \dist(x^*, S(q^{k_l})) \right) \ge \psi(\delta) > 0 \quad \forall l \geq 0,
    \end{align*}
which yields
    \begin{align*}
       \max_{i=1,\dots,m} q^*_i - q^{k_l}_i + v(q^*) -v(q^{k_l}) \ge  \psi(\delta) > 0 \quad \forall l \geq 0.
    \end{align*}
We let $l \to +\infty$ and use $q^{k_l} \to q^*$ and the continuity of the optimal value function to derive a contradiction.
\end{proof}


\section{Tikhonov regularization for multiobjective optimization}
\label{sec:tikh_moo}

In this section we extend the concept of Tikhonov regularization from single objective to multiobjective optimization and study the properties of the associated regularization path. The obtained results will play a crucial role in the asymptotic analysis we perform in the following sections for \eqref{eq:MTRIGS}.

A fundamental concept in the study of Tikhonov regularization when minimizing a convex and differentiable function $f:\H \to \R$, is the regularization path. This path, defined as $\{x_\varepsilon: \varepsilon > 0\}$, is a smooth and bounded curve where each $x_\varepsilon$ is the unique minimizer of $f +  \frac{\varepsilon}{2}\lVert \cdot \rVert^2$. As $\varepsilon \to 0$, it holds $x_\varepsilon \to \proj_{\argmin f}(0)$  (see, for instance, \cite[Theorem 27.23]{Bauschke2017}). The regularization path is crucial in the asymptotic analysis conducted in  \cite{Attouch2022} for \eqref{eq:TRIGS_intro}, where the convergence of the trajectory solution $x(\cdot)$ to the minimum norm solution was demonstrated by showing that  $\lim_{t \to + \infty} \lVert x(t) - x_{\varepsilon(t)} \rVert = 0$. We aim to extend this idea to the multiobjective setting when studying \eqref{eq:MOP} and the dynamical system \eqref{eq:MTRIGS}.

Although the analysis presented in this section holds in a more general form for any continuously differentiable function $\varepsilon:[t_0, +\infty) \to (0,+\infty)$ that is nonincreasing and satisfies $\lim_{t \to +\infty} \varepsilon(t) = 0$, we restrict the analysis in this paper to the case $\varepsilon(t) = \frac{\beta}{t^p}$ in order to be
consistent with the formulation of the system \eqref{eq:MTRIGS}. Define for all $t \ge t_0$
\begin{align}
\label{eq:MOP_Tikhonov}
\tag{MOP$_{\frac{\beta}{t^p}}$}
    \min_{x\in \H} \left[
    \begin{array}{c}
        f_{t, 1}(x)\\
        \vdots \\
        f_{t, m}(x)
    \end{array}
    \right]
    \coloneqq 
    \left[
    \begin{array}{c}
        f_1(x) + \frac{\beta}{2 t^p}\lVert x \rVert^2 \\
        \vdots \\
        f_m(x) + \frac{\beta}{2 t^p}\lVert x \rVert^2
    \end{array}
    \right],
\end{align}
where 
\begin{align*}
    f_{t, i}:\H \to \R,\quad x \mapsto f_i(x) + \frac{\beta}{2 t^p}\lVert x \rVert^2, \ \mbox{for} \ i = 1,\dots,m.
\end{align*}

Although the functions $f_{t,i}$ are strongly convex, one cannot expect \eqref{eq:MOP_Tikhonov} to have a unique Pareto optimal solution. This necessitates a suitable concept of a regularization path. To address this, we utilize the merit function defined in \eqref{eq:merit_function} for the regularized problem \eqref{eq:MOP_Tikhonov}, that we define for all $t \ge t_0$ as
\begin{align}
\label{eq:varphi_t}
    \varphi_t:\H \to \R, \quad x \mapsto  \sup_{z \in \H} \min_{i = 1, \dots, m} f_{t,i}(x) - f_{t,i}(z) = \sup_{z \in \H} \min_{i = 1, \dots, m} f_{i}(x) - f_{i}(z) + \frac{\beta}{2 t^p}\lVert x \rVert^2 - \frac{\beta}{2 t^p}\lVert z \rVert^2.
\end{align}
The merit function can be interpreted as the Pascoletti-Serafini scalarization of the problem \eqref{eq:MOP_Tikhonov} (see, for instance, \cite[Section 2.1]{Eichfelder2008}). Inspired by the formulation of the merit function and by the Tikhonov regularization in the single objective case, we consider for all $t \geq t_0$ the unique minimizer of the problem
\begin{align}
\label{eq:varphi_t_opt_prob_3}
    \min_{z \in \H} \max_{i=1,\dots,m} f_{i}(z) - f_{i}(x(t)) + \frac{\beta}{2 t^p}\lVert z \rVert^2
\end{align}
as an element of the regularization path, where $x:[t_0, +\infty) \to \H$ is a trajectory which will be specified later. Note that for the single objective case, namely when $m=1$, we recover the classical regularization path independent of the trajectory $x(\cdot)$. Since the function $z \mapsto  \max_{i=1,\dots,m}  f_{i}(z) - f_{i}(x(t))$ depends on $t$, we cannot make use of the properties of the regularization path in the single objective case to characterize the asymptotic behavior of this new path. This will be done in the following result.

\begin{theorem}
\label{thm:moo_tikh_conv}
 Let $q:[t_0, + \infty) \to \R^m$ be a continuous function with $q(t) \to q^* 
 \in \R^m$ as $t \to + \infty$, and
    \begin{align}
    \label{eq:z_q_z0}
        \begin{split}
            z(t) & \coloneqq \argmin_{z \in \H} \max_{i=1,\dots,m} f_i(z) - q_i(t) + \frac{\beta}{2 t^p}\lVert z \rVert^2 \text{ for all } t \ge t_0,\\
            S(q) & \coloneqq \argmin_{z \in \H} \max_{i=1,\dots,m} f_i(z) - q_i \text{ for all } q \in \R^m,\\
            z_0(q) & \coloneqq \proj_{S(q)}(0) \text{ for all } q \in \R^m.
        \end{split}
    \end{align}
    Then, $z(t) \to z_0(q^*)$ strongly converges as $t \to + \infty$.
\end{theorem}

\begin{proof}
    Let $(t_k)_{k \ge 0} \subset [t_0, + \infty)$ be an arbitrary sequence with $t_k \to + \infty$ as $k \to +\infty$. For all $k \geq 0$, we denote  $\varepsilon_k \coloneqq \frac{\beta}{(t_k)^p}$, $q^k \coloneqq q(t_k)$, $z^k \coloneqq z(t_k)$, and $z_0^{k} \coloneqq z_0(q^k)$. For all $k \geq 0$ it holds
    \begin{align}
        \begin{split}
            \label{eq:tikhonov0}
            \max_{i=1,\dots,m} f_i(z^k) - q_i^k + \frac{\varepsilon_k}{2}\lVert z^k \rVert^2 \le \max_{i=1,\dots,m} f_i(z_0^{k}) - q_i^k + \frac{\varepsilon_k}{2}\lVert z_0^{k} \rVert^2 \le \max_{i=1,\dots,m} f_i(z^k) - q_i^k + \frac{\varepsilon_k}{2}\lVert z_0^{k} \rVert^2,
        \end{split}
    \end{align}
  hence,
    \begin{align}
    \label{eq:reg_sol_bound}
        \lVert z^k \rVert \le \lVert z_0^{k} \rVert.
    \end{align}
    According to assumption \ref{ass:A4}, $z_0(\cdot)$ is continuous, consequently, $\{z_0^k\}_{k \ge 0}$ is bounded. This implies that $\{z^k\}_{k \ge 0}$ is also bounded and hence possesses a weak sequential cluster point. We show that this point is unique, which will imply that $\{z^k\}_{k \ge 0}$ is weakly convergent. 
    
    Let $z^{\infty}$ be an arbitrary weak sequential cluster point of $\{z^k\}_{k \ge 0}$, and a subsequence $z^{k_l} \rightharpoonup z^{\infty}$ as $l \to + \infty$. For all $z \in \H$ it holds
    \begin{align}
    \label{eq:tikhonov1}
        \begin{split}
            & \max_{i=1,\dots,m} \left( f_i(z^{\infty}) - q_i^* \right) \le \liminf_{l \to +\infty} \max_{i=1,\dots,m} \left( f_i(z^{k_l}) - q_i^* \right) + \frac{\varepsilon_{k_l}}{2}\lVert z^{k_l} \rVert^2\\
            \le & \liminf_{l \to +\infty} \left(\max_{i=1,\dots,m} \left( f_i(z^{k_l}) - q_i^{k_l} \right) + \frac{\varepsilon_{k_l}}{2}\lVert z^{k_l} \rVert^2 + \max_{i=1,\dots,m} \left( q_i^{k_l} - q_i^* \right) \right) \\
            \le & \liminf_{l \to +\infty} \left(\max_{i=1,\dots,m} \left( f_i(z) - q_i^{k_l} \right) + \frac{\varepsilon_{k_l}}{2}\lVert z \rVert^2 \right)\\
            \le & \liminf_{l \to +\infty} \left(\max_{i=1,\dots,m} \left( f_i(z) - q_i^{*} \right) + \frac{\varepsilon_{k_l}}{2}\lVert z \rVert^2 + \max_{i=1,\dots,m} \left( q_i^* - q_i^{k_l} \right) \right)\\
            = & \max_{i=1,\dots,m} \left( f_i(z) - q_i^{*} \right).
        \end{split}
    \end{align}
    From here, $z^{\infty} \in S(q^*)$ follows. Next, we show that $z^{\infty} = z_0(q^*)$. From the continuity of $z_0(\cdot)$ we have
        \begin{align}
        \label{eq:projection_se_continuous}
            z_0^{k_l} = z_0(q^{k_l}) \to z_0(q^*) \, \text{ as }\, l \to + \infty,
        \end{align}
    and the weak lower semicontinuity of the norm gives
        \begin{align}
        \label{eq:norm_ineq_weak_limits}
            \lVert z^{\infty} \rVert \le \liminf_{l \to + \infty} \lVert z^{k_l} \rVert \le \limsup_{l \to + \infty} \lVert z^{k_l} \rVert \le \limsup_{ l \to + \infty} \lVert z_0^{k_l} \rVert = \lVert z_0(q^*) \rVert.
        \end{align}
    Since $z^{\infty} \in S(q^*)$ and $z_0(q^*) = \proj_{S(q^*)}(0)$, we get $z^{\infty} = z_0(q^*)$. This proves that $\lbrace z^k \rbrace_{k \ge 0}$  weakly converges to $z_0(q^*)$. Using again \eqref{eq:norm_ineq_weak_limits}, we get
    \begin{align*}
        \lim_{k \to + \infty} \lVert z^k \rVert = \lVert z_0(q^*) \rVert, 
    \end{align*}
    from which we conclude that $z^k \to z_0(q^*)$ strongly converges as $k \to +\infty$. 
\end{proof}

\begin{myremark}
    The continuity of $z_0( \cdot )$ formulated in assumption \ref{ass:A4} can be seen as a regularity condition on the objective functions $f_i$ for $i = 1, \dots, m$. It is satisfied for convex single objective optimization problems as long as the set of minimizers is not empty. In this setting the mapping $q \to z_0(q)$ is constant. The following example shows that the assumption \ref{ass:A4} is crucial for obtaining convergence of $z(t)$ as $t \to +\infty$.
\end{myremark}

\begin{myexample}
\label{ex:counter_example_cont_proj}
    Define the functions
    \begin{align}
    \label{eq:example_functions}
    \begin{split}
        \phi:\R \to \R, \quad & y \mapsto \frac{1}{2}\max\left(y - 3, 0\right)^2 + \frac{1}{2}\max\left(2 - y, 0\right)^2, \\\\
        g:\R^2 \to \R, \quad & x \mapsto \left\lbrace\begin{array}{lll}
        \frac{1}{2} x_1^2 + \frac{1}{2} x_2^2, & \text{ if }\, \lvert x_1 \rvert \le 1,& x_2 + 1 \le \sqrt{1 - x_1^2},\\[4pt] 
       \lvert x_1 \rvert + \frac{1}{2}x_2^2 - \frac{1}{2}, & \text{ if }\, \lvert x_1 \rvert  > 1,& x_2 + 1 \le 0,\\[2pt]
       \sqrt{x_1^2 + (x_2 + 1)^2} - (x_2 + 1),& \text{ else, }& \\
    \end{array} \right.\\\\
        f_1:\R^2 \to \R , \quad & x \mapsto \frac{1}{2}(x_1 - 1)^2 +\phi(x_2) + g(x) , \\\\
        f_2:\R^2 \to \R, \quad & x \mapsto \frac{1}{2}(x_1 + 1)^2 + \phi(x_2) + g(x),
    \end{split}
    \end{align}
which are all convex and differentiable with Lipschitz continuous gradients (see \ref{appendix_c}). We consider the multiobjective optimization problem
    \begin{align}
    \label{eq:example_MOP}
    \begin{split}
        \min_{x \in \H} \left[ \begin{array}{c}
            f_1(x)  \\
            f_2(x)
        \end{array} \right],
    \end{split}
    \end{align}
    and the Tikhonov regularized problem
    \begin{align}
    \label{eq:example_MOP_Tikh}
    \begin{split}
        \min_{x \in \H} \left[ \begin{array}{c}
            f_1(x) + \frac{\varepsilon}{2} \lVert x \rVert^2 \\
            f_2(x) + \frac{\varepsilon}{2} \lVert x \rVert^2
        \end{array} \right].
    \end{split}    
    \end{align}
    Figure \ref{subfig:example_pareto_set_z(t)_a} illustrates the weak Pareto set $\Pw$ of the problem \eqref{eq:example_MOP} alongside the Pareto set of the regularized problem \eqref{eq:example_MOP_Tikh} for various values of $\varepsilon > 0$ denoted by $\mathcal{P}_{w, \varepsilon}$. As $\varepsilon$ decreases, the weak Pareto set of \eqref{eq:example_MOP_Tikh} ``converges'' to the weak Pareto set of \eqref{eq:example_MOP}. Due to the T-shape of the weak Pareto set, the edges of the regularized weak Pareto sets become sharper as $\varepsilon$ diminishes. For this problem the map
    \begin{align*}
        z_0: \R^2 \to \R^2,\quad q \mapsto z_0(q) = \proj_{S(q)}(0),
    \end{align*}
    with $S(q) = \argmin\limits_{z \in \R^2} \max\left( f_1(z) - q_1 , f_2(z) - q_2 \right)$ is not continuous everywhere. Indeed,
    \begin{align*}
        z_0(q_1, 0) \to \left(0, 3\right) \neq (0, 2) = \proj_{\{ 0\} \times [2,3]}(0) = z_0((0,0)) \text{ as } q_1 \to 0.
    \end{align*}
 We define, for $t_0 := \left( 192 \beta \right)^{\frac{1}{p}}$,
    \begin{align*}
        q:[t_0, +\infty) \to \R, \quad t \mapsto \left[ \begin{array}{c}
            q_1(t)  \\
            q_2(t)
        \end{array} \right] \coloneqq \left[\begin{array}{c}
       2 (\omega(t) + 1)\sqrt{\left(\frac{t^p}{ t^p - \beta \omega(t)}\right)^2 - 1}\\
        0   
        \end{array}\right],
    \end{align*}
    with $\omega(t) \coloneqq \frac{10 + \sin(\eta t)}{4}$, where $\eta > 0$ is a positive scaling parameter. It holds $q(t) \to q^* = (0,0)^{\top}$ as $t \to +\infty$. For this example the regularization path is given for all $t \geq t_0$ by
    \begin{align}
    \label{eq:reg_path_z(t)}
        z(t) = \left[\begin{array}{c}
        - (\omega(t) + 1)\sqrt{\left(\frac{t^p}{ t^p - \beta \omega(t)}\right)^2 - 1}\\
        \omega(t)   
        \end{array}\right] \in \argmin_{z \in \R^2} \max \left( f_1(z) - q_1(t), f_2(z) - q_2(t) \right) + \frac{\beta}{2 t^p}\lVert z \rVert^2.
    \end{align}
    In Figure \ref{fig:example_pareto_set_z(t)} (b), the regularization path $z(\cdot)$ given by \eqref{eq:reg_path_z(t)} is depicted. One can observe that it oscillates in the $x_2$-coordinate between the values $2.25$ and $2.75$ as $t \to +\infty$. The function $z(t)$ does not converge as $t \to +\infty$, although all accumulation points are weak Pareto optimal and global minimizers of $\max \left( f_1(z) - q_1^*, f_2(z) - q_2^* \right)$. The minimal norm solution $z_0(q^*) = \left(0, 2\right)$ is not an accumulation point of $z(\cdot)$. This example clearly shows that the continuity of $z_0(\cdot)$ is essential to derive Theorem \ref{thm:moo_tikh_conv}.
    \begin{figure}
    \begin{center}
        \begin{subfigure}[t]{.41\textwidth}
            \centering
            \includegraphics[width=\linewidth]{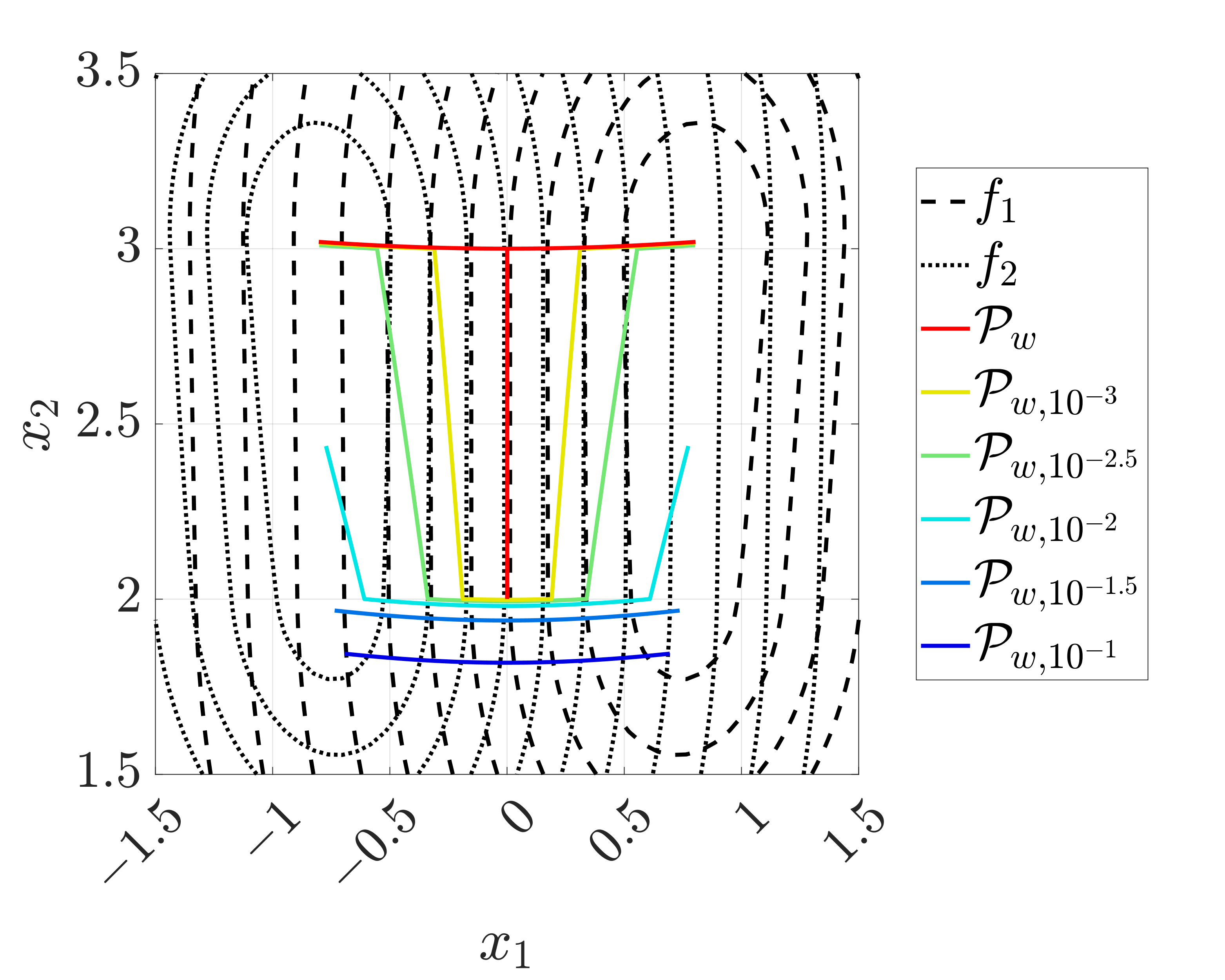}
            \caption{}
            \label{subfig:example_pareto_set_z(t)_a}
        \end{subfigure}
        \vspace{5mm}
        \begin{subfigure}[t]{.4\textwidth}
            \centering
            \includegraphics[width=\linewidth]{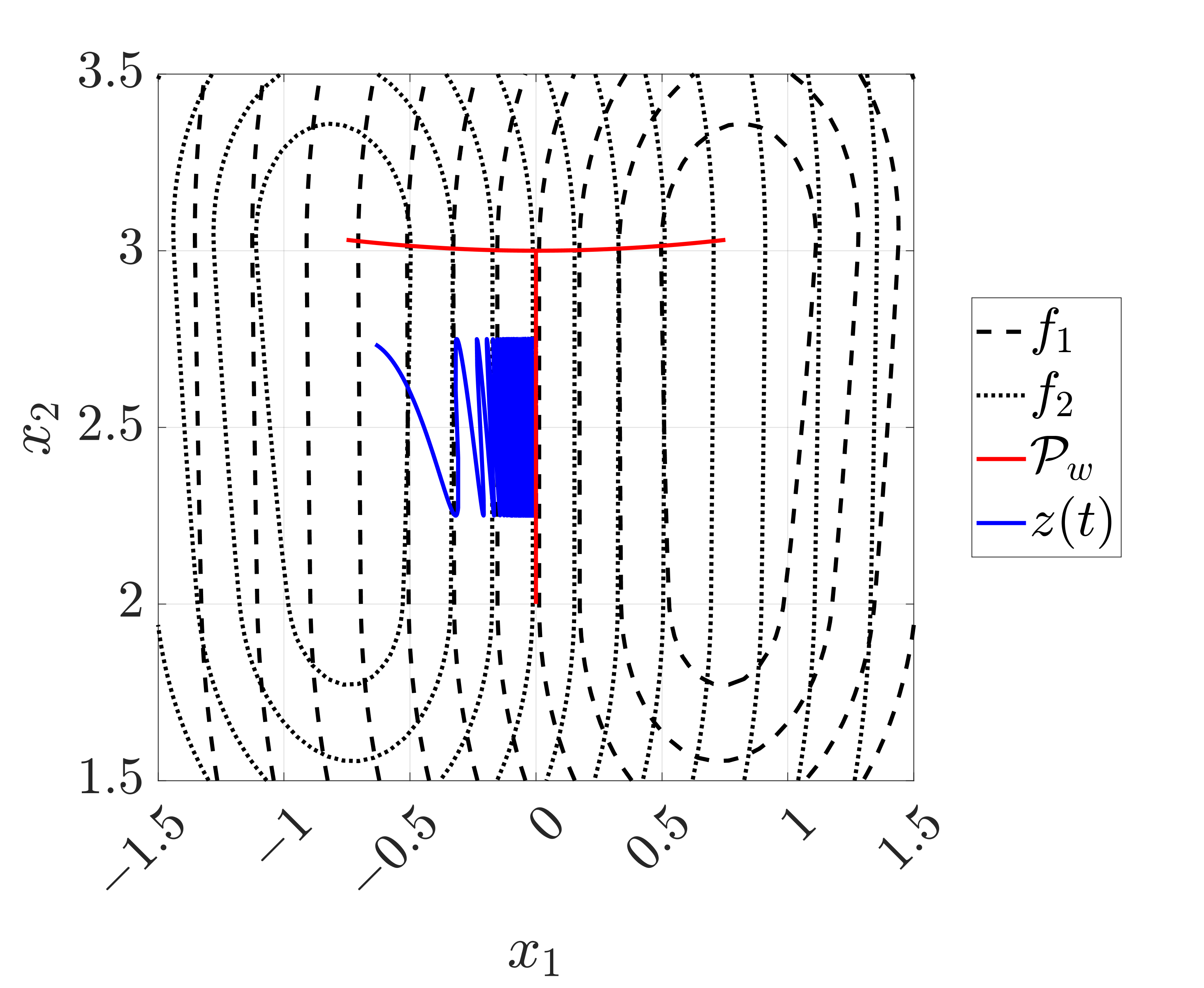}
            \caption{}
            \label{subfig:example_pareto_set_z(t)_b}
        \end{subfigure}    
    \end{center}
    \caption{Contour plots of the functions $f_1$ and $f_2$ defined in \eqref{eq:example_functions}: (a) The weak Pareto sets of \eqref{eq:example_MOP} and \eqref{eq:example_MOP_Tikh} for $\varepsilon \in \{ 10^{-1}, 10^{-1.5},10^{-2}, 10^{-2.5},10^{-3} \}$. (b) The weak Pareto set of \eqref{eq:example_MOP} and the regularization path $z(\cdot)$ defined in \eqref{eq:reg_path_z(t)} with parameters $p = 1$, $\beta = \frac{1}{2}$, $\eta = \frac{1}{50}$.}
    \label{fig:example_pareto_set_z(t)}
\end{figure}
\end{myexample}

We conclude this section by introducing three propositions that summarize the main properties of $z(\cdot)$. 

\begin{myprop}
\label{prop:z_bounded}
Let $a \in \R_{+}^{m}$ and assume that the trajectory solution $x:[t_0, +\infty) \to \H$ fulfills $x(t) \in \L(F, F(x(t_0)) + a)$ for all $t \ge t_0$. Then, the regularization path,
     \begin{align*}
         z(t) \coloneqq \argmin_{z \in \H} \max_{i=1,\dots,m} f_i(z) - f_i(x(t)) + \frac{\beta}{2 t^p}\lVert z \rVert^2,\quad\text{for all }\,t \geq t_0,
     \end{align*}
is bounded. Specifically, $z(t) \in B_{R}(0)$ for all $t \ge t_0$, where $R$ is defined in \ref{ass:A3}.
\end{myprop}

\begin{proof}
    By \ref{ass:A4}, it holds $S(F(x(t))) \coloneqq \argmin\limits_{z \in \H} \max_{i=1,\dots,m} \left(f_i(z) - f_i(x(t))\right) \neq \emptyset$ for all $t \ge t_0$. Fix some $t \ge t_0$. 
    
    From the properties of Tikhonov regularization in single objective optimization (cf. \cite[Theorem 27.23]{Bauschke2017}), we know
    \begin{align}
    \label{eq:ineq_z_t_<_z}
        \lVert z(t) \rVert \le \lVert z \rVert \quad \forall z \in S(F(x(t))).
    \end{align}

    Next, we show that 

    \begin{align}
    \label{eq:F_-1_F*_subset_S_x_t}
        F^{-1}(\{F^*\}) \subseteq S(x(t)) \quad \forall F^* \in F(S(F(x(t))) .
    \end{align}

    Let $F^* \in F(S(F(x(t)))$. Then, there exists $z \in S(F(x(t))$ with $F(z) = F^*$. Let $w \in F^{-1}(\{F^*\})$ then $F(w) = F(z)$ and hence
    \begin{align*}
        \max_{i=1,\dots,m} f_i(w) - f_i(x(t)) = \max_{i=1,\dots,m} f_i(z) - f_i(x(t)) = \inf_{z \in \H} \max_{i = 1,\dots,m} f_i(z) - f_i(x(t)).
    \end{align*}
    This shows $w \in S(F(x(t))$ and hence \eqref{eq:F_-1_F*_subset_S_x_t} holds. From \eqref{eq:ineq_z_t_<_z} and \eqref{eq:F_-1_F*_subset_S_x_t} we conclude that for all $F^* \in F(S(F(x(t))))$ we get
    \begin{align*}
        \lVert z(t) \rVert \le \lVert z \rVert \quad \forall z \in F^{-1}(\{F^*\}),
    \end{align*}
    and hence
    \begin{align*}
        \lVert z(t) \rVert \le \inf_{z \in F^{-1}(\{F^*\})} \lVert z \rVert \quad \forall F^* \in F(S(F(x(t)))).
    \end{align*}
    Since this bound holds for all $F^* \in F(S(F(x(t))))$, we get
    
   \begin{align}
   \label{eq:sup_inf_ineq_0}
        \lVert z(t) \rVert \le \inf_{z \in F^{-1}(F(S(F(x(t)))))} \lVert z \rVert = \inf_{\{z \in \H : F(z) \in F(S(F(x(t))))\}} \lVert z \rVert \le \sup_{F^* \in F(S(F(x(t)))} \inf_{z \in F^{-1}(\{F^*\})} \lVert z \rVert.
    \end{align}

    Next, we prove that 
    \begin{align}
    \label{inclusion}
    S(F(x(t))) \subseteq \LPw(F, F(x(t_0)) + a). 
    \end{align}
    Let $z \in S(F(x(t)))$. Then,
        \begin{align*}
            \max_{i=1,\dots,m} f_i(z) - f_i(x(t)) \le \max_{i=1,\dots,m} f_i(x(t)) - f_i(x(t)) = 0,
        \end{align*}
    hence
        \begin{align*}
            f_i(z) \le f_i(x(t)) \le f_i(x(t_0)) + a_i \quad \forall i=1, \ldots, m,
        \end{align*}
        and therefore $z \in \L(F, F(x(t_0)) + a)$. Assuming that $z \not\in \LPw(F, F(x(t_0)) + a)$, it follows that $z \not\in \Pw$ and hence there exists some $y \in \H$ with
        \begin{align*}
            f_i(y) < f_i(z) \text{ for all } i = 1,\dots,m.
        \end{align*}
    Therefore,
        \begin{align*}
            \max_{i=1,\dots,m} f_i(x) - f_i(x(t)) < \max_{i=1,\dots,m} f_i(z) - f_i(x(t)),
        \end{align*}
    which is a contradiction to $z \in S(F(x(t)))$. This proves inclusion \eqref{inclusion}. Consequently, according to \eqref{eq:sup_inf_ineq_0} and \eqref{inclusion},
        \begin{align*}
            \lVert z(t) \rVert \le \sup_{F^* \in F(\LPw(F, F(x(t_0)) + a)} \inf_{z \in F^{-1}(\{F^*\})} \lVert z \rVert = R < + \infty,
        \end{align*}
        where the upper bound $R$ is given by \ref{ass:A3}.
\end{proof}

\begin{myprop}
Let $q:[t_0, + \infty) \to \R^m$ be a continuous function and
    \begin{align*}
         z(t) \coloneqq \argmin_{z \in \H} \max_{i=1,\dots,m} f_i(z) - q_i(t) + \frac{\beta}{2 t^p}\lVert z \rVert^2 \ \mbox{for all} \ t \geq t_0.
     \end{align*}
     Then, $z(\cdot)$ is a continuous mapping.
\end{myprop}

\begin{proof}
We fix an arbitrary $\overline{t} \ge t_0$ and show that $z(\cdot)$ is continuous  (continuous from the right if $\overline{t}=t_0$) in $\overline{t}$. Let $t \in \left[\overline{t} - \kappa, \overline{t}+ \kappa\right] \cap \left[t_0, +\infty\right)$ for some $\kappa > 0$. Then, by strong convexity and the minimizing properties of $z(t)$ and $z(\overline{t})$, we get
    \begin{align}
    \label{eq:z_t_continuous_1}
        \begin{split}
            & \max_{i=1,\dots,m} \left( f_i(z(\overline{t})) - q_i(t) \right) + \frac{\beta}{2 t^p}\lVert z(\overline{t}) \rVert^2\\ - & \max_{i=1,\dots,m} \left( f_i(z(t)) - q_i(t) \right) - \frac{\beta}{2 t^p}\lVert z(t) \rVert^2 \ge \frac{\beta}{2 t^p} \lVert z(\overline{t}) - z(t) \rVert^2,
        \end{split}
    \end{align}
    and
    \begin{align}
    \label{eq:z_t_continuous_2}
        \begin{split}
            & \max_{i=1,\dots,m} \left( f_i(z(t)) - q_i(\overline{t}) \right) + \frac{\beta}{2 \overline{t}^p}\lVert z(t) \rVert^2 \\ - & \max_{i=1,\dots,m} \left( f_i(z(\overline{t}) - q_i(\overline{t}) \right) - \frac{\beta}{2 \overline{t}^p}\lVert z(\overline{t}) \rVert^2 \ge \frac{\beta}{2 \overline{t}^p} \lVert z(t) - z(\overline{t}) \rVert^2,
        \end{split}
    \end{align}
    respectively. Using the monotonicity of $t \mapsto \frac{\beta}{2 t^p}$, \eqref{eq:z_t_continuous_1} and \eqref{eq:z_t_continuous_2} lead to
    \begin{align}
    \label{eq:z_t_continuous_3}
        \begin{split}
            &\max_{i=1,\dots,m} \left( f_i(z(\overline{t}) - q_i(\overline{t}) \right) + \max_{i=1,\dots,m} \left(  q_i(\overline{t}) - q_i(t) \right) + \frac{\beta}{2 t^p}\lVert z(\overline{t}) \rVert^2 \\ - & \max_{i=1,\dots,m} \left( f_i(z(t) - q_i(t) \right) - \frac{\beta}{2 t^p}\lVert z(t) \rVert^2 \ge \frac{\beta}{2 (\overline{t} + \kappa)^p} \lVert z(\overline{t}) - z(t) \rVert^2,
        \end{split}
    \end{align}
respectively,
    \begin{align}
    \label{eq:z_t_continuous_4}
        \begin{split}
            &\max_{i=1,\dots,m} \left( f_i(z(t)) - q_i(t) \right) + \max_{i=1,\dots,m} \left( q_i(t) - q_i(\overline{t}) \right)+ \frac{\beta}{2 \overline{t}^p}\lVert z(t) \rVert^2\\ - & \max_{i=1,\dots,m} \left( f_i(z(\overline{t}) - q_i(\overline{t}) \right) - \frac{\beta}{2 \overline{t}^p}\lVert z(\overline{t}) \rVert^2 \ge \frac{\beta}{2 (\overline{t} + \kappa)^p} \lVert z(t) - z(\overline{t}) \rVert^2.
        \end{split}
    \end{align}
Adding \eqref{eq:z_t_continuous_3} and \eqref{eq:z_t_continuous_4} yields
    \begin{align}
    \label{eq:z_t_continuous_5}
        2\lVert q(t) - q(\overline{t}) \rVert_{\infty} + \frac{1}{2}\left( \frac{\beta}{\overline{t}^p} - \frac{\beta}{t^p} \right) \left( \Vert z(t) \rVert^2 -  \lVert z(\overline{t}) \rVert^2 \right) \ge \frac{\beta}{(\overline{t} + \kappa)^p} \lVert z(t) - z(\overline{t}) \rVert^2.
    \end{align}
    By Proposition \ref{prop:z_bounded}, the function $z(\cdot)$ is bounded, so by the continuity of $q(\cdot)$ the left-hand-side of \eqref{eq:z_t_continuous_5} vanishes as $t \to \overline{t}$. This demonstrates the continuity of $z(\cdot)$ in $\overline{t}$.
\end{proof}

In the next proposition, we describe the connection between the original merit function $\varphi(\cdot)$ and the merit function $\varphi_t(\cdot)$ of the regularized problem. This will allow us to derive asymptotic convergence results on $\varphi(x(t))$ for $t \to +\infty$.

\begin{myprop}
\label{prop:basic_properties_tikh}
Let $a \in \R_{+}^{m}$ be the vector introduced in assumption \ref{ass:A3} and assume that $x:[t_0, +\infty) \to \H$ fulfills $x(t) \in \L(F, F(x(t_0)) + a)$ for all $t \ge t_0$. We define 
     \begin{align*}
         z(t) \coloneqq \argmin_{z \in \H} \max_{i=1,\dots,m} f_i(z) - f_i(x(t)) + \frac{\beta}{2 t^p}\lVert z \rVert^2 \ \mbox{for all} \ t \geq t_0.
     \end{align*}
     Then, the following statements hold:
     \begin{enumerate}
         \item[i\emph{)}] For all $t \ge t_0$ and all $y \in \H$
             \begin{align*}
                \min_{i=1,\dots,m} f_i(x(t)) - f_i(y) \le \min_{i=1,\dots,n} f_{t,i}(x(t)) - f_{t,i}(z(t)) + \frac{\beta}{2 t^p}\lVert y \rVert^2,
             \end{align*}
         hence
             \begin{align*}
                 \varphi(x(t)) \le \varphi_t(x(t)) + \frac{\beta R^2}{2 t^p} ,
             \end{align*}
             where $R$ is defined in \ref{ass:A3}.
        \item[ii\emph{)}] For all $t \ge t_0$
     \begin{align*}
         \lVert x(t) - z(t) \rVert^2 \le \frac{t^p \varphi_t(x(t))}{\beta}.
     \end{align*}
     \end{enumerate}
\end{myprop}
\begin{proof} 

$i)$ Fix $t \ge t_0$ and $y \in \H$. From the definition of $z(t)$, we have
    \begin{align*}
        \max_{i=1,\dots,m} f_{t,i}(y) - f_{t,i}(x(t)) \ge \max_{i=1,\dots,m} f_{t,i}(z(t)) - f_{t,i}(x(t)),
    \end{align*}
 hence
    \begin{align*}
        \min_{i=1,\dots,m} f_{i}(x(t)) - f_{i}(y) + \frac{\beta}{2 t^p} \lVert x(t) \rVert^2 - \frac{\beta}{2 t^p} \lVert y \rVert^2 \le \min_{i=1,\dots,m} f_{t,i}(x(t)) - f_{t,i}(z(t)).
    \end{align*}
    Using the definition of $\varphi_t(\cdot)$, we get
    \begin{align}
        \label{eq:prop 2.7_i_1}
        \min_{i=1,\dots,m} f_{i}(x(t)) - f_{i}(y) \le \varphi_t(x(t)) + \frac{\beta}{2 t^p} \lVert y \rVert^2.
    \end{align}
    By \ref{ass:A3}, it holds $\LPw(F,F(x(t_0))+a) \neq \emptyset$, therefore,
    \begin{align}
        \label{eq:prop 2.7_i_3}
        & \sup_{F^* \in  F(\LPw(F, F(x(t_0)) + a))}\inf_{y \in F^{-1}(\{F^*\})} \min_{i=1,\dots,m} f_{i}(x(t)) - f_{i}(y) \nonumber \\
        \le & \ \varphi_t(x(t)) + \frac{\beta}{2 t^p} \sup_{F^* \in  F(\LPw(F, F(x(t_0)) + a))} \inf_{y \in F^{-1}(\{F^*\})} \lVert y \rVert^2.
    \end{align}
    Additionally, we have
    \begin{align}
        \label{eq:prop 2.7_i_4}
        \sup_{y \in \LPw(F, F(x(t_0)) + a)} \min_{i=1,\dots,m} f_{i}(x(t)) - f_{i}(y) = \sup_{F^* \in F(\LPw(F, F(x(t_0)) + a))}\inf_{y \in F^{-1}(\{F^*\})} \min_{i=1,\dots,m} f_{i}(x(t)) - f_{i}(y).
    \end{align}
    Note that \eqref{eq:prop 2.7_i_4} holds since for all $y \in \LPw(F,F(x(t_0)) + a)$ there exists $F^* = F(y) \in F(\LPw(F,F(x(t_0)) + a))$ with $\min_{i=1,\dots,m} f_i(x(t)) - f_i(y) = \min_{i=1,\dots,m} f_i(x(t)) - f_i(z)$ for all $z \in F^{-1}(\{F^*\})$. On the other hand, for all $F^* \in F(\LPw(F,F(x(t_0)) + a))$ any $y \in  \LPw(F, F(x(t_0)) + a)$ with $F(y) = F^*$ satisfies  $\min_{i=1,\dots,m} f_i(x(t)) - f_i(y) = \inf_{z \in F^{-1}(\{F^*\})} \min_{i=1,\dots,m} f_i(x(t)) - f_i(z)$. Combining \eqref{eq:prop 2.7_i_3} and \eqref{eq:prop 2.7_i_4}, and using Lemma \ref{lem:sup_inf_u_0} and \ref{ass:A3}, it yields
    \begin{align*}
        \varphi(x(t)) \le \varphi_t(x(t)) + \frac{\beta R^2}{2 t^p}.
    \end{align*}
$ii)$ From the strong convexity of $f_{t,i}$ with modulus $\frac{\beta}{t^p}$, we conclude the strong convexity of $z \mapsto \max_{i = 1,\dots,m} f_{t,i}(z) - f_{t,i}(x(t))$ with modulus $\frac{\beta}{t^p}$. This gives for all $t \ge t_0$
    \begin{align*}
        \varphi_t(x(t)) = & \ \min_{i=1,\dots,m} f_{t,i}(x(t)) - f_{t,i}(z(t)) \\ 
        = & \ \max_{i=1,\dots,m} f_{t,i}(x(t)) - f_{t,i}(x(t)) - \max_{i=1,\dots,m} f_{t,i}(z(t)) - f_{t,i}(x(t)) \\
        \ge & \ \frac{\beta}{t^p} \lVert x(t) - z(t) \rVert^2,
    \end{align*}
    and the desired inequality follows.
\end{proof}


\section{Existence of solutions and some preparatory results for the asymptotic analysis}
\label{sec:MTRIGS}

In this section, we discuss the existence of solution trajectories of the dynamical system \eqref{eq:MTRIGS} and derive their properties which will be used in the asymptotic analysis.


\subsection{Existence of trajectory solutions}
\label{subsec:existence_uniqueness}

The existence of solutions of \eqref{eq:MTRIGS} follows analogously to that shown for the system \eqref{eq:MAVD_intro} in \cite{Sonntag2024b} and requires the Hilbert space $\H$ to be finite dimensional. We only give the definition of solutions and the main existence theorem in this subsection and move the proof to Appendix \ref{sec:ap_exist_sol_MTRIGS}. 

Due to the implicit structure of the differential equation \eqref{eq:MTRIGS}, we do not expect the trajectory solutions $x(\cdot)$ to be twice continuously differentiable in general. However, we show that there are continuously differentiable solutions with an absolutely continuous first derivative. The following definition describes what we understand by a solution of \eqref{eq:MTRIGS}.

\begin{mydef}
\label{def:sol_CP}
    We call a function $x:[t_0, +\infty) \to \H,\, t\mapsto x(t)$ a solution to \eqref{eq:MTRIGS} if it satisfies the following conditions:
    \begin{enumerate}[(i)]
        \item $x(\cdot) \in C^1([t_0, +\infty))$, i.e., $x(\cdot) $ is continuously differentiable on $[t_0, +\infty)$;
        \item $\dot{x}(\cdot)$ is absolutely continuous on $[t_0, T]$ for all $T \ge t_0$;
        \item There exists a (Bochner) measurable function $\ddot{x}:[t_0, +\infty) \to \H$ with $\dot{x}(t) = \dot{x}(t_0) + \int_{t_0}^t \ddot{x}(s) ds$ for all $t \ge t_0$;
        \item $\dot{x}(\cdot)$ is differentiable almost everywhere and $\frac{d}{dt} \dot{x} (t) = \ddot{x}(t)$ holds for almost all $t \in [t_0, +\infty)$;
        \item $\frac{\alpha}{t^q}\dot{x}(t) + \proj_{C(x(t)) + \frac{\beta}{t^p} x(t)+ \ddot{x}(t)}{(0)} = 0$ holds for almost all $t \in [t_0, +\infty)$;
        \item $x(t_0) = x_0$ and $\dot{x}(t_0) = v_0$.
    \end{enumerate}
\end{mydef}

Next, we give the main existence theorem for solution to \eqref{eq:MTRIGS}.

\begin{theorem}
   Assume $\H$ is finite dimensional. Then, for all initial values $(x_0, v_0) \in \H \times \H$ there exists a function $x(\cdot)$ which is a solution of \eqref{eq:MTRIGS} in the sense of Definition \ref{def:sol_CP}.
\end{theorem}

\begin{proof}
    See the proof of Theorem \ref{thm:main_existence_theorem_appendix} in Appendix \ref{sec:ap_exist_sol_MTRIGS}.
\end{proof}

\begin{myremark}
    The uniqueness of the trajectory solutions of \eqref{eq:MTRIGS} remains an open problem. There are two major difficulties in deriving uniqueness, as for the dynamical system \eqref{eq:MAVD_intro}. First, the multiobjective steepest descent direction is not Lipschitz continuous, but only Hölder continuous. So even for simpler multiobjective gradient-like systems like $\dot{x}(t) = \proj_{C(x(t))}(0)$ it is not trivial to show uniqueness of trajectories in the general setting. The second problem is the implicit structure of the equation \eqref{eq:MTRIGS}. Therefore, we cannot use standard arguments like the Cauchy-Lipschitz theorem to derive the uniqueness of solutions. Note that the asymptotic analysis performed in this paper applies to any trajectory solution $x(\cdot)$ of \eqref{eq:MTRIGS}, which reduces the importance of the uniqueness statement.
\end{myremark}

\subsection{Preparatory results for the asymptotic analysis}
\label{subsec:energy_estimates}

In this subsection, we derive some properties that all trajectory solution $x(\cdot)$ of the system \eqref{eq:MTRIGS} share. 

\begin{myprop}
\label{prop:MTRIGS_proj_variational}
    Let $x(\cdot)$ be a trajectory solution of \eqref{eq:MTRIGS}. Then, for all $i = 1,\dots, m$ and almost all $t \ge t_0$ it holds
    \begin{align*}
        \left\langle \nabla f_i(x(t)) + \frac{\beta}{t^p} x(t) + \ddot{x}(t) + \frac{\alpha}{t^q} \dot{x}(t) , \dot{x}(t) \right\rangle \le 0,
    \end{align*}
    and therefore
    \begin{align*}
        \left\langle \nabla f_i(x(t)) + \frac{\beta}{t^p} x(t) + \ddot{x}(t) , \dot{x}(t) \right\rangle \le - \frac{\alpha}{t^q} \lVert \dot{x}(t) \Vert^2. 
    \end{align*}
\end{myprop}
\begin{proof}
    According to Definition \ref{def:sol_CP}, each solution $x(\cdot)$ satisfies 
    \begin{align*}
        - \frac{\alpha}{t^q} \dot{x}(t) = \proj_{C(x(t)) + \frac{\beta}{t^p} x(t) + \ddot{x}(t)}{(0)},
    \end{align*}
    for almost all $t \ge t_0$. From the variational characterization of the projection, it follows that 
    \begin{align*}
        \left\langle \nabla f_i(x(t)) + \frac{\beta}{t^p} x(t) + \ddot{x}(t) + \frac{\alpha}{t^q} \dot{x}(t), \frac{\alpha}{t^q}\dot{x}(t) \right\rangle \le 0, 
    \end{align*}
    for almost all $t \ge t_0$ and all $i = 1,\dots, m$, which leads to the desired inequality.
\end{proof}
In the next proposition, we define component-wise a multiobjective energy function and show that its components fulfill a decay property along each trajectory solution.
\begin{myprop}
\label{prop:ineq_W_i_t}
    Let $x(\cdot)$ be a trajectory solution of \eqref{eq:MTRIGS}. For all $i=1,\dots,m$, we define the energy function
    \begin{align}
    \label{eq:def_W_i_t}
        \W_i:[t_0, +\infty) \to \R, \quad t \mapsto f_i(x(t)) + \frac{\beta}{2 t^p}\lVert x(t) \rVert^2 + \frac{1}{2}\lVert \dot{x}(t) \rVert^2.
    \end{align}
    Then, for all $i=1,\dots,m$ and almost all $t \ge t_0$ it holds
    \begin{align*}
        \frac{d}{dt}\W_i(t) \le -\frac{p \beta}{2 t^{p+1}} \lVert x(t) \rVert^2 - \frac{\alpha}{t^q} \lVert \dot{x}(t) \Vert^2 \le 0.
    \end{align*}
    Further, for $a \in \R_{+}^m$ defined as $a_i := \frac{\beta}{2t_0^p}\lVert x(t_0) \Vert^2 + \frac{1}{2}\lVert \dot{x}(t_0) \rVert^2$ for $i = 1,\dots,m$, it holds
    \begin{align*}
        x(t) \in \L(F, F(x(t_0)) + a) \quad \text{ for all }\ t \ge t_0.
    \end{align*}
\end{myprop}
\begin{proof}
   According to Definition \ref{def:sol_CP}, the velocity $\dot{x}(\cdot)$ of a trajectory solution is differentiable almost everywhere. For all $i=1,\dots,m$ and almost all $t \ge t_0$ it holds
    \begin{align*}
            \frac{d}{dt}\W_i(t) = & \ \langle \nabla f_i(x(t)) , \dot{x}(t) \rangle - \frac{p \beta}{2 t^{p+1}}\lVert x(t) \rVert^2 + \frac{\beta}{t^p} \langle x(t), \dot{x}(t) \rangle + \langle \dot{x}(t), \ddot{x}(t) \rangle \nonumber \\
            = & \ - \frac{p \beta}{2 t^{p+1}} \lVert x(t) \rVert^2 + \left \langle \nabla f_i(x(t)) + \frac{\beta}{t^p} x(t) + \ddot{x}(t), \dot{x}(t) \right \rangle \\
            \le & \ - \frac{p \beta}{2 t^{p+1}} \lVert x(t) \rVert^2 - \frac{\alpha}{t^q} \lVert \dot{x}(t) \Vert^2 \leq 0,
    \end{align*}
where the penultimate inequality follows from Proposition \ref{prop:MTRIGS_proj_variational}. The last statement of the proposition follows using the monotonicity of each $\W_i$ for $i=1, \dots , m,$ on $[t_0, +\infty)$.
\end{proof}

Since for almost all $t \ge t_0$, $\proj_{C(x(t)) + \frac{\beta}{t^p} x(t) + \ddot{x}(t)}(0)$ belongs to $C(x(t)) + \frac{\beta}{t^p} x(t) + \ddot{x}(t)$, there exists $\theta(t) \in \Delta^m \coloneqq \left\lbrace \theta \in \R_{+}^m \, : \, \sum_{i=1}^m \theta_i = 1 \right\rbrace$ such that
\begin{align}
\label{eq:proj_convex_combination}
    -\frac{\alpha}{t^q} \dot{x}(t) = \proj_{C(x(t)) + \frac{\beta}{t^p} x(t) + \ddot{x}(t)}(0) = \sum_{i=1}^m \theta_i(t) \nabla f_i(x(t)) + \frac{\beta}{t^p} x(t) + \ddot{x}(t).
\end{align}
In the following proposition, we show that there exists a measurable function $\theta(\cdot)$ satisfying \eqref{eq:proj_convex_combination}.
\begin{myprop}
    \label{prop:theta_measurable}
   Let $x(\cdot)$ be a trajectory solution of \eqref{eq:MTRIGS}. Then, there exists a measurable function
    \begin{align*}
        \theta: [t_0, +\infty) \to \Delta^m, \quad t \mapsto \theta(t),
   \end{align*}
which satisfies for almost all $t \ge t_0$
    \begin{align}
    \label{eq:nice_theta_equality}
        -\frac{\alpha}{t^q} \dot{x}(t) = \proj_{C(x(t)) + \frac{\beta}{t^p} x(t) + \ddot{x}(t)}(0) = \sum_{i=1}^m \theta_i(t) \nabla f_i(x(t)) + \frac{\beta}{t^p} x(t) + \ddot{x}(t).
    \end{align}
\end{myprop}
\begin{proof}
    The proof follows the lines of the proof of Lemma 4.3 in \cite{Sonntag2024b}, where a similar result was shown for the system \eqref{eq:MAVD_intro}. For almost all $t \ge t_0$, there exists  $\theta(t) \in \Delta^m$ such that
    \begin{align}
    \label{eq:theta_measurable}
        \theta(t) \in \argmin_{\theta \in \Delta^m} j(t, \theta), \text{ where } j(t, \theta) \coloneqq \left\lVert \sum_{i = 1}^m \theta_i \nabla f_i(x(t)) + \frac{\beta}{t^p} x(t) + \ddot{x}(t) \right\rVert^2.
    \end{align}
The existence of a measurable selection $\theta : [t_0, +\infty) \to \Delta^m, t \mapsto \theta(t) \in \argmin_{\theta \in \Delta^m} j(t, \theta)$ can be verified using \cite[Theorem 14.37]{Rockafellar2009}. To this end, we have to show that $j(\cdot, \cdot)$ is a Carathéodory integrand, i.e., $j(\cdot,\theta)$ is measurable for all $\theta$ and $j(t,\cdot)$ is continuous for all $t \geq t_0$. The second condition is obviously satisfied. Since $x(\cdot)$ is a trajectory solution of \eqref{eq:MTRIGS} in the sense of Definition \ref{def:sol_CP}, $\ddot{x}(\cdot)$ is (Bochner) measurable. Hence, for all $\theta \in \Delta^m$, $j(\theta, \cdot)$ is measurable as a composition of measurable and continuous functions. This demonstrates that the first condition is also satisfied.
\end{proof}

By using the weight function $\theta(\cdot)$ we can give a further variational characterization of a trajectory solution of \eqref{eq:MTRIGS}.
\begin{myprop}
    \label{lem:variational_bound_d_dt_f(x(t))}
   Let $x(\cdot)$ be a trajectory solution of \eqref{eq:MTRIGS} and $\theta: [t_0, +\infty) \to \Delta^m$ the corresponding measurable weight function given by Proposition  \ref{prop:theta_measurable}. Then, for all $i = 1,\dots, m$ and almost all $t \ge t_0$ it holds
    \begin{align*}
        \langle \nabla f_i(x(t)), \dot{x}(t) \rangle \le \left\langle \sum_{i=1}^m \theta_i(t) \nabla f_i(x(t)), \dot{x}(t)\right\rangle.
    \end{align*}
    \end{myprop}
    \begin{proof}
    By Proposition \ref{prop:MTRIGS_proj_variational}, we have for all $i=1,\dots,m$ and almost all $t \ge t_0$
    \begin{align}
        \label{eq:d_dt_W_i_t_b1}
        \left \langle \nabla f_i(x(t)) + \frac{\beta}{t^p} x(t) + \ddot{x}(t) + \frac{\alpha}{t^q} \dot{x}(t), \dot{x}(t) \right \rangle \le 0,
    \end{align}
which, combined with \eqref{eq:nice_theta_equality}, yields
    \begin{align*}
        \langle \nabla f_i(x(t)) , \dot{x}(t) \rangle \le \left\langle \sum_{i=1}^m \theta_i(t) \nabla f_i(x(t)), \dot{x}(t) \right\rangle.
    \end{align*}
\end{proof}

We conclude this section with the following proposition.

\begin{myprop}
\label{lem:boundedness_x}
 Let $x(\cdot)$ be a trajectory solution of \eqref{eq:MTRIGS}. Then, the following statements are true:
    \begin{enumerate}
        \item[i\emph{)}] $\dot{x}(\cdot)$ is bounded;
        \item[ii\emph{)}] if $x(\cdot)$ is bounded, then $\ddot{x}(\cdot)$ is essentially bounded.
    \end{enumerate}
\end{myprop}

\begin{proof} \emph{i}) According to Proposition \ref{prop:ineq_W_i_t}, we have for all $i = 1,\dots,m$ and all $t \ge t_0$
\begin{align*}
    \frac{1}{2}\lVert \dot{x}(t) \rVert^2 \le \W_i(t) \le \W_i(t_0),
\end{align*}
which proves the first statement.

\emph{ii}) If $x(\cdot)$ is bounded, then $\nabla f_i(x(\cdot))$ is also bounded for all $i = 1,\dots,m,$ as a consequence of the Lipschitz continuity of the gradients. According to \eqref{eq:MTRIGS}, we have for almost all $t \geq t_0$
\begin{align*}
    \ddot{x}(t) + \frac{\alpha}{t^q} \dot{x}(t) = \proj_{C(x(t)) + \frac{\beta}{t^p} x(t)}(-\ddot{x}(t)),
\end{align*}
hence,
\begin{align}
\label{eq:bound_ddot_x}
    \left \lVert \ddot{x}(t) \right \rVert \le \frac{\alpha}{t^q} \left \lVert \dot{x}(t) \right \rVert  + \left \lVert  \proj_{C(x(t)) + \frac{\beta}{t^p} x(t)}(-\ddot{x}(t)) \right \rVert.
\end{align}
Since all expressions on the right hand side of \eqref{eq:bound_ddot_x} are bounded on $[t_0, +\infty)$,  $\ddot{x}(\cdot)$ is essentially bounded.
\end{proof}


\section{Asymptotic analysis}
\label{sec:convergence_anaylsis}

In this section, we study the asymptotic behavior of the trajectory solutions to \eqref{eq:MTRIGS}. The convergence rates for the merit function values and the convergence of the trajectory depend heavily on the parameters $p \in (0,2], q\in(0,1]$ and $\alpha, \beta > 0$. The results in this section extend those in \cite{Laszlo2023} from the single objective to the multiobjective framework. The following energy functions are the key to the asymptotic analysis of \eqref{eq:MTRIGS}.

\begin{mydef}
\label{def:mo_energy_functions5}
Let $x(\cdot)$ be a trajectory solution of \eqref{eq:MTRIGS}, $r \in [q,1]$ and $z \in \H$. Let $\gamma: [t_0, + \infty) \to [0, + \infty)$ and $\xi: [t_0, + \infty) \to \R$ be continuously differentiable functions. We define for $i = 1,\dots,m$
\begin{align*}
        \G^r_{i, \gamma, \xi, z}(t) \coloneqq & \ t^{2r} \left( f_{t,i}(x(t)) - f_{t,i}(z) \right) + \frac{1}{2}\lVert \gamma(t)(x(t) - z) + t^r \dot{x}(t) \rVert^2 + \frac{\xi(t)}{2} \lVert x(t) - z \rVert^2
\end{align*}
and
\begin{align*}
        \G^r_{\gamma, \xi, z}(t) \coloneqq & \ \,t^{2r} \min_{i=1,\dots,m}\left( f_{t,i}(x(t)) - f_{t,i}(z) \right) + \frac{1}{2}\lVert \gamma(t)(x(t) - z) + t^r \dot{x}(t) \rVert^2 + \frac{\xi(t)}{2} \lVert x(t) - z \rVert^2.
    \end{align*}
For $z(t) \coloneqq \argmin_{z \in \H} \max_{i=1,\dots,m} f_{t,i}(z) - f_{t,i}(x(t))$ for $t \geq t_0$, we define
 \begin{align*}
        \G^r_{\gamma, \xi}:[t_0, + \infty) \to \R, \quad t \mapsto \G^r_{\gamma, \xi,z(t)}(t) = & \ \,t^{2r} \min_{i=1,\dots,m}\left( f_{t,i}(x(t)) - f_{t,i}(z(t)) \right) \\
        & + \frac{1}{2}\lVert \gamma(t)(x(t) - z(t)) + t^r \dot{x}(t) \rVert^2 + \frac{\xi(t)}{2} \lVert x(t) - z(t) \rVert^2.\\
        = & \ \,t^{2r} \varphi_t(x(t)) \\
        & + \frac{1}{2}\lVert \gamma(t)(x(t) - z(t)) + t^r \dot{x}(t) \rVert^2 + \frac{\xi(t)}{2} \lVert x(t) - z(t) \rVert^2.
    \end{align*}
\end{mydef}

The functions $\gamma(\cdot)$ and $\xi(\cdot)$ will be specified at a later point in the analysis. In the next proposition, we derive estimates for the derivatives of the energy functions introduced above.

\begin{myprop}
\label{prop:ineq_d_dt_E_z_t2}
Let $x(\cdot)$ be a trajectory solution of \eqref{eq:MTRIGS}, $r \in [q,1]$ and $z \in \H$. Let $\gamma: [t_0, + \infty) \to [0, + \infty)$ and $\xi: [t_0, + \infty) \to \R$ be continuously differentiable functions.
\begin{enumerate}[i)]
    \item For all $i = 1,\dots,m$, the function $\G^r_{i, \gamma, \xi, z}(\cdot)$ is absolutely continuous on every interval $[t_0, T]$ for $T \geq t_0$, differentiable almost everywhere on $[t_0, +\infty)$, and its derivative satisfies for almost all $t \in [t_0, +\infty)$ 
    \begin{align}
    \label{eq:equation_deriv_G_i_z_t2}
    \begin{split}
        \frac{d}{dt}\G^r_{i, \gamma, \xi, z}(t) \le & \ 2rt^{2r - 1}\left( f_{t,i}(x(t)) - f_{t,i}(z) \right) - t^r \gamma(t) \min_{i=1,\dots,m} \left( f_{t,i}(x(t)) - f_{t,i}(z) \right) + \frac{p \beta t^{2r}}{2t^{p+1}} \lVert z \rVert^2\\
    & + \left (\gamma(t)(\gamma(t) + rt^{r-1} - \alpha t^{r - q}) + t^r \gamma'(t) + \xi(t) \right) \left\langle x(t) - z , \dot{x}(t)\right\rangle  \\
    & +  \left( \gamma(t)\gamma'(t) + \frac{\xi'(t)}{2} - \gamma(t) t^r \frac{\beta }{2 t^p} \right) \lVert x(t) - z \Vert^2 + t^r(\gamma(t) + rt^{r-1} - \alpha t^{r - q}) \lVert \dot{x}(t) \rVert^2.
    \end{split}
    \end{align}
    \item The function $\G^r_{\gamma, \xi, z}(\cdot)$ is absolutely continuous on every interval $[t_0, T]$ for $T \geq t_0$, differentiable almost everywhere on $[t_0, +\infty)$, and its derivative satisfies for almost all $t \in [t_0, +\infty)$ 
    \begin{align}
    \label{eq:equation_deriv_G_z_t2}
    \begin{split}
        \frac{d}{dt}\G^r_{\gamma, \xi, z}(t) \le & \ \left( 2rt^{2r - 1} - t^r \gamma(t) \right) \min_{i=1,\dots,m} \left( f_{t,i}(x(t)) - f_{t,i}(z) \right) + \frac{p \beta t^{2r}}{2t^{p+1}} \lVert z \rVert^2\\
     & + \left( \gamma(t)(\gamma(t) + rt^{r-1} - \alpha t^{r - q}) + t^r \gamma'(t) + \xi(t) \right) \left\langle x(t) - z , \dot{x}(t)\right\rangle \\
    & +  \left( \gamma(t)\gamma'(t) + \frac{\xi'(t)}{2} - \gamma(t) t^r \frac{\beta }{2 t^p} \right) \lVert x(t) - z \Vert^2  + t^r(\gamma(t) + rt^{r-1} - \alpha t^{r - q}) \lVert \dot{x}(t) \rVert^2.
    \end{split}
    \end{align}
\end{enumerate}
\end{myprop}

\begin{proof} 
    Fix an arbitrary $i \in \{1,\dots,m\}$. It is obvious that $\G^r_{i, \gamma, \xi, z}(\cdot)$ is absolutely continuous on every interval $[t_0, T]$ for $T \geq t_0$ and therefore differentiable almost everywhere on $[t_0, +\infty)$. Let $t \ge t_0$ be a point at which $\G^r_{i,\gamma, z}(\cdot)$ is differentiable. By the chain rule, it holds that
    \begin{align*}
        \frac{d}{dt} \G^r_{i, \gamma, \xi, z}(t) = & \ 2rt^{2r - 1}\left( f_{t,i}(x(t)) - f_{t,i}(z) \right) + t^{2r} \langle \nabla f_{t,i}(x(t)) , \dot{x}(t) \rangle - \frac{p \beta t^{2r}}{2 t^{p+1}} \lVert x(t) \rVert^2 + \frac{p \beta t^{2r}}{2 t^{p+1}} \lVert z \rVert^2 \\
        + & \left\langle \gamma(t)(x(t) - z) + t^r \dot{x}(t) , (\gamma(t) + rt^{r-1}) \dot{x}(t) + \gamma'(t)(x(t) - z) + t^{r} \ddot{x}(t) \right\rangle \\
        + & \ \xi(t) \langle x(t) - z , \dot{x}(t) \rangle + \frac{\xi'(t)}{2} \lVert x(t) - z \rVert^2.
    \end{align*}
Let $\theta(\cdot)$ be the measurable weight function given by Proposition \ref{prop:theta_measurable}. By Proposition \ref{lem:variational_bound_d_dt_f(x(t))}, we have
    \begin{align}
        \label{eq:inequality_10_1}
        \begin{split}
            \frac{d}{dt} \G^r_{i, \gamma, \xi, z}(t) \le & \ 2rt^{2r - 1}\left( f_{t,i}(x(t)) - f_{t,i}(z) \right) + t^{2r} \left\langle \sum_{i=1}^m \theta_i(t) \nabla f_{t,i}(x(t)) , \dot{x}(t) \right\rangle  + \frac{p \beta t^{2r}}{2t^{p+1}} \lVert z \rVert^2\\
            & + \left\langle \gamma(t)(x(t) - z) + t^r \dot{x}(t) , (\gamma(t) + rt^{r-1}) \dot{x}(t) + \gamma'(t)(x(t) - z) + t^{r} \ddot{x}(t) \right\rangle \\
            & + \xi(t) \langle x(t) - z , \dot{x}(t) \rangle + \frac{\xi'(t)}{2} \lVert x(t) - z \rVert^2.
        \end{split}
    \end{align}
    Using \eqref{eq:nice_theta_equality}, we write
    \begin{align*}
        t^r \ddot{x}(t) =  - \alpha t^{r - q} \dot{x}(t) - t^r \sum_{i=1}^m \theta_i(t) \nabla f_{t,i}(x(t)),
    \end{align*}
which we use to evaluate
    \begin{align}
        \label{eq:inequality_10_2}
        \begin{split}
            & \left\langle \gamma(t)(x(t) - z) + t^r \dot{x}(t) , (\gamma(t) + rt^{r-1}) \dot{x}(t) + \gamma'(t)(x(t) - z) + t^{r} \ddot{x}(t) \right\rangle \\
            = & \left\langle \gamma(t)(x(t) - z) + t^r \dot{x}(t) , (\gamma(t) + rt^{r-1} -\alpha t^{r - q}) \dot{x}(t) + \gamma'(t)(x(t) - z) - t^r \sum_{i=1}^m \theta_i(t) \nabla f_{t,i}(x(t)) \right\rangle \\
            = & \ \gamma(t)(\gamma(t) + rt^{r-1} - \alpha t^{r - q})\left\langle x(t) - z , \dot{x}(t)\right\rangle + \gamma(t)\gamma'(t) \lVert x(t) - z \Vert^2 - t^r\gamma(t) \left\langle x(t) - z , \sum_{i=1}^m \theta_i(t) \nabla f_{t,i}(x(t)) \right\rangle\\
            & \ + t^r(\gamma(t) + rt^{r-1} - \alpha t^{r - q}) \lVert \dot{x}(t) \rVert^2 + t^r \gamma'(t) \langle \dot{x}(t) , x(t) - z \rangle - t^{2r} \left\langle \sum_{i=1}^m \theta_i(t) \nabla f_{t,i}(x(t)), \dot{x}(t) \right\rangle \\
            = & \left[ \gamma(t)(\gamma(t) + r t^{r-1} - \alpha t^{r - q}) + t^r \gamma'(t) \right]\left\langle x(t) - z , \dot{x}(t)\right\rangle + \gamma(t)\gamma'(t) \lVert x(t) - z \Vert^2 \\
            & - t^r\gamma(t) \left\langle x(t) - z , \sum_{i=1}^m \theta_i(t) \nabla f_{t,i}(x(t)) \right\rangle + t^r(\gamma(t) + rt^{r-1} - \alpha t^{r - q}) \lVert \dot{x}(t) \rVert^2  - t^{2r} \left\langle \sum_{i=1}^m \theta_i(t) \nabla f_{t,i}(x(t)), \dot{x}(t) \right\rangle.
        \end{split}
    \end{align}
    We combine \eqref{eq:inequality_10_1} and \eqref{eq:inequality_10_2} to derive 
    \begin{align}
    \label{eq:inequality_10_234}
        \begin{split}
         \frac{d}{dt} \G^r_{i, \gamma, \xi, z}(t) \le & \ 2rt^{2r - 1}\left( f_{t,i}(x(t)) - f_{t,i}(z) \right) + t^{2r} \left\langle \sum_{i=1}^m \theta_i(t) \nabla f_{t,i}(x(t)) , \dot{x}(t) \right\rangle + \frac{p \beta t^{2r}}{2t^{p+1}} \lVert z \rVert^2 \\
        & + \left (\gamma(t)(\gamma(t) + rt^{r-1} - \alpha t^{r - q}) + t^r \gamma'(t) \right) \left\langle x(t) - z , \dot{x}(t)\right\rangle +  \gamma(t)\gamma'(t) \lVert x(t) - z \Vert^2 \\
        & - t^r \gamma(t) \left\langle x(t) - z , \sum_{i=1}^m \theta_i(t) \nabla f_{t,i}(x(t)) \right\rangle + t^r(\gamma(t) + rt^{r-1} - \alpha t^{r - q}) \lVert \dot{x}(t) \rVert^2\\
        & - t^{2r} \left\langle \sum_{i=1}^m \theta_i(t) \nabla f_{t,i}(x(t)), \dot{x}(t) \right\rangle + \xi(t) \langle x(t) - z , \dot{x}(t) \rangle + \frac{\xi'(t)}{2} \lVert x(t) - z \rVert^2 \\
        = & \ 2rt^{2r - 1}\left( f_{t,i}(x(t)) - f_{t,i}(z) \right) + \frac{p \beta t^{2r}}{2t^{p+1}} \lVert z \rVert^2\\
        & + \left( \gamma(t)(\gamma(t) + rt^{r-1} - \alpha t^{r - q}) + t^r \gamma'(t) + \xi(t) \right) \left\langle x(t) - z , \dot{x}(t)\right\rangle +  \left( \gamma(t)\gamma'(t) + \frac{\xi'(t)}{2}  \right) \lVert x(t) - z \Vert^2 \\
        & + t^r\gamma(t) \left\langle z - x(t) , \sum_{i=1}^m \theta_i(t) \nabla f_{t,i}(x(t)) \right\rangle + t^r(\gamma(t) + rt^{r-1} - \alpha t^{r - q}) \lVert \dot{x}(t) \rVert^2.
        \end{split}
    \end{align}
    We use the strong convexity of $x \mapsto \sum_{i=1}^m \theta_i(t) (f_{t,i}(x) - f_{t,i}(z))$ to derive
    \begin{align}
    \label{eq:inequality_10_235}
        \begin{split}
        \left\langle z - x(t) , \sum_{i=1}^m \theta_i(t) \nabla f_{t,i}(x(t)) \right\rangle \le & \ \sum_{i=1}^m \theta_i(t) \left( f_{t,i}(z) - f_{t,i}(x(t)) \right) - \frac{\beta}{2t^p}\lVert x(t) - z \rVert^2 \\
        \le &  \ -\min_{i=1,\dots,m} f_{t,i}(x(t)) - f_{t,i}(z) - \frac{\beta}{2t^p}\lVert x(t) - z \rVert^2.
        \end{split}
    \end{align}
Plugging \eqref{eq:inequality_10_234} into \eqref{eq:inequality_10_235} gives
    \begin{align*}
    \frac{d}{dt} \G^r_{i, \gamma, \xi, z}(t) \le & \ 2rt^{2r - 1}\left( f_{t,i}(x(t)) - f_{t,i}(z) \right) - t^r \gamma(t) \min_{i=1,\dots,m} \left( f_{t,i}(x(t)) - f_{t,i}(z) \right) - \gamma(t) t^r \frac{\beta }{2 t^p}\lVert x(t) - z \rVert^2 \\
        & +  \left( \gamma(t)(\gamma(t) + rt^{r-1} - \alpha t^{r - q}) + t^r \gamma'(t) + \xi(t) \right) \left\langle x(t) - z , \dot{x}(t)\right\rangle +  \left( \gamma(t)\gamma'(t) + \frac{\xi'(t)}{2}  \right) \lVert x(t) - z \Vert^2 \\
        & + t^r(\gamma(t) + rt^{r-1} - \alpha t^{r - q}) \lVert \dot{x}(t) \rVert^2 + \frac{p \beta t^{2r}}{2t^{p+1}} \lVert z \rVert^2,
    \end{align*}
    concluding part \emph{i}). Statement \emph{ii}) follows immediately from \emph{i}) and Lemma \ref{lem:diff_a_e}.
\end{proof}

For given $\lambda > 0$ and $r \in [q,1]$, we choose in the first part of the convergence analysis
    \begin{align*}
        \gamma:[t_0, +\infty) \to [0,+\infty), \, t\mapsto \gamma(t) \coloneqq \lambda, \quad \text{ and } \quad
        \xi:[t_0, +\infty) \to \R, \, t\mapsto\xi(t) \coloneqq \lambda\left( rt^{r-1} + \alpha t^{r-q} - 2\lambda \right).
    \end{align*}
For this choice of the two parameter functions, we rename the energy functions as follows:
\begin{align*}
        \E^r_{i,\lambda, z}:[t_0, + \infty) \to \R, \quad \E_{i, \lambda, z}(t)  :=  \G^r_{i,\gamma,\xi,z}(t) \coloneqq & \ t^{2r} (f_{t,i}(x(t)) - f_{t,i}(z)) + \frac{1}{2}\left\lVert \lambda (x(t) - z) + t^{r} \dot{x}(t) \right\rVert^2 \\
        & + \frac{\lambda}{2}\left( rt^{r-1} + \alpha t^{r-q} - 2 \lambda \right) \lVert x(t) - z \rVert^2,
\end{align*}
for $i=1, ..., m$,
\begin{align*}
        \E^r_{\lambda, z}:[t_0, + \infty) \to \R, \quad \E^r_{\lambda, z}(t) :=  \G^r_{\gamma,\xi,z}(t) = & \ t^{2r} \min_{i=1,\dots,m}(f_{t,i}(x(t)) - f_{t,i}(z)) + \frac{1}{2}\left\lVert \lambda (x(t) - z) + t^{r} \dot{x}(t) \right\rVert^2 \\
        & + \frac{\lambda}{2}\left( rt^{r-1} + \alpha t^{r-q} - 2 \lambda \right) \lVert x(t) - z \rVert^2,
    \end{align*}
and
 \begin{align*}
        \E^r_{\lambda}:[t_0, + \infty) \to \R, \quad \E^r_{\lambda}(t) := \G^r_{\gamma,\xi}(t)= & \ t^{2r} \min_{i=1,\dots,m}(f_{t,i}(x(t)) - f_{t,i}(z(t))) + \frac{1}{2}\left\lVert \lambda (x(t) - z(t)) + t^{r} \dot{x}(t) \right\rVert^2 \\
        & + \frac{\lambda}{2}\left( rt^{r-1} + \alpha t^{r-q} - 2 \lambda \right) \lVert x(t) - z(t) \rVert^2 \\
        = & \ t^{2r} \varphi_t(x(t)) + \frac{1}{2}\left\lVert \lambda (x(t) - z(t)) + t^{r} \dot{x}(t) \right\rVert^2 \\
        & + \frac{\lambda}{2}\left( rt^{r-1} + \alpha t^{r-q} - 2 \lambda \right) \lVert x(t) - z(t) \rVert^2,
    \end{align*}
where $z(t) \coloneqq \argmin_{z \in \H} \max_{i=1,\dots,m} f_{t,i}(z) - f_{t,i}(x(t))$ for $t \geq t_0$. In the following, we formulate a proposition on $\E^r_{i,\lambda, z}(\cdot)$ and $\E^r_{\lambda, z}(\cdot)$ similar to Proposition \ref{prop:ineq_d_dt_E_z_t2}. 

\begin{myprop}
\label{prop:ineq_d_dt_E_z_t}
Let $x(\cdot)$ be a trajectory solution of \eqref{eq:MTRIGS}, $\lambda > 0$, $r \in [q,1]$ and $z \in \H$.
\begin{enumerate}[i)]
    \item For all $i = 1,\dots,m$, the function $\E^r_{i,\lambda, z}(\cdot)$ is absolutely continuous on every interval $[t_0, T]$ for $T \geq t_0$, differentiable almost everywhere on $[t_0, +\infty)$, and its derivative satisfies for almost all $t \in [t_0, +\infty)$ 
    \begin{align}
    \label{eq:equation_deriv_E_i_z_t2}
    \begin{split}
        \frac{d}{dt}\E^r_{i,\lambda,z}(t) \le & \ 2r t^{2r-1}(f_{t,i}(x(t)) - f_{t,i}(z) - \lambda t^r \min_{i=1,\dots,m}(f_{t,i}(x(t)) - f_{t,i}(z)) + \frac{p \beta t^{2r}}{2t^{p+1}} \lVert z \rVert^2\\
        & + \lambda \left( 2 r t^{r-1} - \lambda \right) \langle x(t) - z , \dot{x}(t) \rangle + t^{r}\left( \lambda + r t^{r - 1} - \alpha t^{r - q} \right) \lVert \dot{x}(t) \rVert^2 \\
        & + \frac{\lambda}{2}\left( r(r-1) t^{r - 2} + \alpha (r-q) t^{r - q - 1} - \beta t^{r - p} \right) \lVert x(t) - z \rVert^2.
    \end{split}
    \end{align}
    \item The functions $\E^r_{\lambda, z}(\cdot)$ is absolutely continuous on every interval $[t_0, T]$ for $T \geq t_0$, differentiable almost everywhere on $[t_0, +\infty)$, and its derivative satisfies for almost all $t \in [t_0, +\infty)$ 
    \begin{align}
    \label{eq:equation_deriv_E_z_t2}
    \begin{split}
        \frac{d}{dt}\E^r_{\lambda, z}(t) \le & \, \left( 2r t^{2r-1} - \lambda t^r \right)\min_{i=1,\dots,m}(f_{t,i}(x(t)) - f_{t,i}(z)) + + \frac{p \beta t^{2r}}{2t^{p+1}} \lVert z \rVert^2\\ 
        & + \lambda \left( 2 r t^{r-1} - \lambda \right) \langle x(t) - z , \dot{x}(t) \rangle + t^{r}\left( \lambda + r t^{r - 1} - \alpha t^{r - q} \right) \lVert \dot{x}(t) \rVert^2 \\
        & + \frac{\lambda}{2}\left( r(r-1) t^{r - 2} + \alpha (r-q) t^{r - q - 1} - \beta t^{r - p} \right) \lVert x(t) - z \rVert^2.
    \end{split}
    \end{align}
\end{enumerate}
\end{myprop}

\begin{proof}
    The proof follows immediately by Proposition \ref{prop:ineq_d_dt_E_z_t2} using $\gamma'(t) = 0$ and $\xi'(t) = \lambda(r(r-1)t^{r-2} + \alpha(r-q)t^{r-q-1})$ for $t \ge t_0$.
\end{proof}

\begin{mylemma}
\label{lem:d_dt_E(t)_mu(t)E(t)}
    Let $q \in (0,1)$, $x(\cdot)$ be a trajectory solution of \eqref{eq:MTRIGS}, $\lambda > 0$, $r \in [q,1)$, and $z \in \H$. Define $\mu_r:[t_0, +\infty) \to \R,\, \mu_r(t) \coloneqq \frac{\lambda}{t^r} - \frac{2r}{t}$. Then, for almost all $t \ge t_1 \coloneqq \max\left(\left( \frac{2r}{\lambda} \right)^{\frac{1}{1-r}}, t_0 \right)$, it holds
    \begin{align}
    \label{eq:d_dt_E(t)_mu(t)E(t)}
        \frac{d}{dt}\E^r_{\lambda, z}(t) + \mu_r(t) \E^r_{\lambda, z}(t) \le t^{r}\left( \frac{3}{2}\lambda - \alpha t^{r - q} \right) \lVert \dot{x}(t) \rVert^2 + \frac{p \beta t^{2r}}{2t^{p+1}} \lVert z \rVert^2 + \frac{\lambda}{2} \left[ \frac{3 \lambda r}{t} - \frac{\lambda^2}{t^r} + \frac{ \lambda \alpha}{t^q} - \frac{\beta}{t^{p-r}}\right] \lVert x(t) - z \rVert^2.
    \end{align}
\end{mylemma}

\begin{proof}
For all $t \ge t_0$ it holds
    \begin{align}
        \label{eq:general_bound_E_z_t}
        \begin{split}
            \E^r_{\lambda}(t) = & \ t^{2r}\min_{i=1,\dots,m}(f_{t,i}(x(t)) - f_{t,i}(z)) + \frac{\lambda^2}{2} \lVert x(t) - z \rVert^2 + \lambda t^r \langle x(t) - z , \dot{x}(t) \rangle \\
            & + \frac{t^{2r}}{2} \lVert \dot{x}(t) \rVert^2 + \frac{\lambda}{2} \left( rt^{r-1} + \alpha t^{r - q} - 2 \lambda \right) \lVert x(t) - z \rVert^2 \\
            = & \ t^{2r}\min_{i=1,\dots,m}(f_{t,i}(x(t)) - f_{t,i}(z)) + \frac{\lambda}{2} \left( rt^{r-1} + \alpha t^{r - q} - \lambda \right) \lVert x(t) - z \rVert^2 \\
            & + \lambda t^r \langle x(t) - z , \dot{x}(t) \rangle + \frac{t^{2r}}{2} \lVert \dot{x}(t) \rVert^2.
        \end{split}
    \end{align}
    Note that $\mu_r(t) \ge 0$ for all $t \ge \left( \frac{2r}{\lambda} \right)^{\frac{1}{1-r}}$. Then, combining \eqref{eq:equation_deriv_E_z_t2} and \eqref{eq:general_bound_E_z_t}, it yields for almost all $t \ge t_1$
    \begin{align*}
    \frac{d}{dt}\E^r_{\lambda, z}(t) + \mu_r(t) \E^r_{\lambda, z}(t) \le & \ \left( 2r t^{2r-1} - \lambda t^r \right)\min_{i=1,\dots,m}(f_{t,i}(x(t)) - f_{t,i}(z)) + t^{r}\left( \lambda + r t^{r - 1} - \alpha t^{r - q} \right) \lVert \dot{x}(t) \rVert^2 \\
        & + \frac{\lambda}{2}\left( r(r-1) t^{r - 2} + \alpha (r-q) t^{r - q - 1} - \beta t^{r - p} \right) \lVert x(t) - z \rVert^2 \\
        & + \lambda \left( 2r t^{r-1} - \lambda \right) \langle x(t) - z , \dot{x}(t) \rangle  + \frac{p \beta t^{2r}}{2t^{p+1}} \lVert z \rVert^2\\
        & + \left( \lambda t^{r} - 2r t^{2r - 1 } \right)\min_{i=1,\dots,m}(f_{t,i}(x(t)) - f_{t,i}(z)) \\
        & + \frac{\lambda}{2} \left[ 
        \frac{3 \lambda r}{t} + \frac{\lambda \alpha}{t^q} - \frac{\lambda^2}{t^r} - \frac{2r^2}{t^{2 - r}} - \frac{2r \alpha}{t^{1 - r + q}} \right] \lVert x(t) - z \rVert^2\\
        & + \lambda \left( \lambda - 2r t^{r - 1} \right) \langle x(t) - z , \dot{x}(t) \rangle + \frac{1}{2} \left( \lambda t^r - 2r t^{2r - 1} \right) \lVert \dot{x}(t) \rVert^2\\
        = & \ t^{r}\left( \frac{3}{2}\lambda - \alpha t^{r - q} \right) \lVert \dot{x}(t) \rVert^2 + \frac{p \beta t^{2r}}{2t^{p+1}} \lVert z \rVert^2\\
        & + \frac{\lambda}{2} \left[ -\frac{r(r + 1)}{t^{2 - r}} - \frac{\alpha (r + q)}{t^{1 - r + q}} + \frac{3 \lambda r}{t} + \frac{ \lambda \alpha}{t^q} - \frac{\lambda^2}{t^r} - \beta t^{r - p} \right] \lVert x(t) - z \rVert^2 \\
        \le & \ t^{r}\left( \frac{3}{2}\lambda - \alpha t^{r - q} \right) \lVert \dot{x}(t) \rVert^2 + \frac{p \beta t^{2r}}{2t^{p+1}} \lVert z \rVert^2 + \frac{\lambda}{2} \left[ \frac{3 \lambda r}{t} - \frac{\lambda^2}{t^r} + \frac{ \lambda \alpha}{t^q} - \frac{\beta}{t^{p-r}}\right] \lVert x(t) - z \rVert^2 .
\end{align*}
\end{proof}

The result above can be extended to the case $q \in (0,1]$ and $r = 1$ for $\lambda \ge 2$ as we state in the following lemma. 

\begin{mylemma}
\label{lem:d_dt_E(t)_mu(t)E(t)_r=1}
    Let $q \in (0,1]$, $x(\cdot)$ be a trajectory solution of \eqref{eq:MTRIGS}, $\lambda \ge 2$, $r = 1$ and $z \in \H$. Define $\mu_1:[t_0, +\infty) \to \R,\, t \mapsto \mu_1(t) \coloneqq \frac{\lambda - 2}{t}$. Then, for almost all $t \ge t_0$, it holds
    \begin{align}
    \label{eq:d_dt_E(t)_mu(t)E(t)_r=1}
    \begin{split}
        \frac{d}{dt}\E^1_{\lambda, z}(t) + \mu_1(t) \E^1_{\lambda, z}(t) \le & \ t \left( \frac{3}{2}\lambda - \alpha t^{1 - q} \right) \lVert \dot{x}(t) \rVert^2 + \frac{p \beta }{2t^{p-1}} \lVert z \rVert^2 \\  
        & + \frac{\lambda}{2} \left[ \frac{(1 - \lambda)(\lambda - 2)}{t} + \frac{\alpha (\lambda - (1 + q)) }{t^q} - \frac{\beta}{t^{p-1}}\right] \lVert x(t) - z \rVert^2.
    \end{split}
    \end{align}
\end{mylemma}

\begin{proof}
    The proof is analogous to that of Lemma \ref{lem:d_dt_E(t)_mu(t)E(t)}.
\end{proof}

\subsection{The case $p \in (0,2]$ and $q < \frac{p}{2}$ : convergence rates}\label{subsec41}

In Theorem \ref{thm:main_convergence_special_r=q} we derive convergence rates for the merit function along trajectory solutions of \eqref{eq:MTRIGS} when $q \in (0,1)$ is such that $p \in (0,2]$ and $q < \frac{p}{2}$.

\begin{theorem}
\label{thm:main_convergence_special_r=q}
    Let $p \in (0,2]$ with $q < \frac{p}{2}$, $x(\cdot)$ be a bounded trajectory solution of \eqref{eq:MTRIGS}, and $z(t) := \argmin_{z \in \H} \max_{i=1,\dots,m} f_{t,i}(z) - f_{t,i}(x(t))$ for $t \geq t_0$. Then, we have the following convergence rates as $t \rightarrow +\infty$:
    \begin{enumerate}[i)]
        \item $\E^q_\lambda(t) = \mathcal{O}\left( 1 \right)$ for $0 < \lambda < \frac{\alpha}{2}$;
        \item $\varphi_t(x(t)) = \mathcal{O}\left( t^{-2q} \right);$
        \item $\varphi(x(t)) = \mathcal{O}\left( t^{-2q} \right);$
        \item $\lVert x(t) - z(t) \Vert = \mathcal{O}\left( 1 \right);$
        \item $\lVert \dot{x}(t) \Vert  = \mathcal{O}\left( t^{- q} \right).$
    \end{enumerate}
\end{theorem}

\begin{proof}
    \emph{i)} Let $0 < \lambda < \frac{\alpha}{2}$ and $z \in \H$ fixed. We derive a bound for the energy function $\E^q_{\lambda, z}(\cdot)$ by considering inequality \eqref{eq:d_dt_E(t)_mu(t)E(t)} with $r = q$, i.e., for almost all $t \ge \max\left(\left( \frac{2q}{\lambda} \right)^{\frac{1}{1-q}}, t_0 \right)$
    \begin{align}
    \label{eq:bound_d_dt_E_z(t)+_mu(t)E_z(t)_prep}
    \begin{split}
        \frac{d}{dt}\E^q_{\lambda, z}(t) + \mu_q(t) \E^q_{\lambda, z}(t) \le t^{q}\left( \frac{3}{2}\lambda - \alpha \right) \lVert \dot{x}(t) \rVert^2 + \frac{p \beta}{2}t^{2q - p - 1} \lVert z \rVert^2 + \frac{\lambda}{2} \left[ \frac{3 \lambda q}{t} - \frac{\lambda^2}{t^q} + \frac{ \lambda \alpha}{t^q} - \frac{\beta}{t^{p-q}}\right] \lVert x(t) - z \rVert^2.
    \end{split}
    \end{align}
   From here, we derive for almost all $t \ge \max\left(\left( \frac{2q}{\lambda} \right)^{\frac{1}{1-q}}, t_0, 1 \right)$
    \begin{align*}
        \frac{d}{dt}\E^q_{\lambda, z}(t) + \mu_q(t) \E^q_{\lambda, z}(t) \le & \ \frac{p \beta}{2} t^{2q - p -1} \lVert z \rVert^2 + \frac{\lambda^2(3 + \alpha - \lambda)}{2 t^q}\lVert x(t) - z \rVert^2\\
         \le & \ \frac{p \beta}{2} t^{2q - p - 1} \lVert z \rVert^2 + \lambda^2(3 + \alpha - \lambda) t^{-q} \left( \lVert z \rVert^2 + \lVert x(t) \rVert^2 \right).
    \end{align*}
    Since $x( \cdot )$ is bounded and $q < \frac{p}{2} \leq 1$, there exist $t_2 \ge \max\left(\left( \frac{2q}{\lambda} \right)^{\frac{1}{1-q}}, t_0, 1 \right)$ and $c, M > 0$ such that for almost all $t \ge t_2$
    \begin{align}
    \label{eq:bound_d_dt_E_z(t)+_mu(t)E_z(t)_32}
        \frac{d}{dt}\E^q_{\lambda, z}(t) + \mu_q(t) \E^q_{\lambda, z}(t) \le c \left( M + \lVert z \Vert^2 \right) t^{-q}.
    \end{align}
    We define the function
    \begin{align}
    \label{eq:def_M(t)}
        \M_q:[t_2, + \infty)\to \R, \quad t \mapsto \M_q(t) \coloneqq \exp\left(\int_{t_2}^t \mu_q(s) ds\right) = \exp \left( \int_{t_2}^t \frac{\lambda}{s^q} - \frac{2q}{s} ds \right) = C_{\M_q} \frac{\exp\left( \frac{\lambda}{1 - q} t^{1 - q}\right)}{t^{2q}},
    \end{align}
    with $C_{\M_q} = \frac{t_2^{2q}}{\exp\left( \frac{\lambda}{1 - q} t_2^{1 - q}\right)} > 0$. The function $\M_q(\cdot)$ is constructed such that $\frac{d}{dt} \M_q(t) = \M_q(t) \mu_q(t)$ and hence
    \begin{align}
    \label{eq:d_dt(M(t)E_z(t))}
        \frac{d}{dt}\left( \M_q(t) \E^q_{\lambda, z}(t)\right) = \M_q(t) \left( \frac{d}{dt}\E^q_{\lambda, z}(t) + \mu_q(t) \E^q_{\lambda, z}(t) \right) \ \mbox{for almost all} \ t \geq t_2.
    \end{align}
    The relations \eqref{eq:d_dt(M(t)E_z(t))} and \eqref{eq:bound_d_dt_E_z(t)+_mu(t)E_z(t)_32} give for almost all $t \ge t_2$
    \begin{align}
    \label{eq:d_dt(M(t)E_z(t))2}
        \frac{d}{dt}\left( \M_q(t) \E^q_{\lambda, z}\right) \le  c\, \M_q(t) \left( M + \lVert z \Vert^2 \right)  t^{-q}.
    \end{align}
    We integrate \eqref{eq:d_dt(M(t)E_z(t))2} from $t_2$ to $t \geq t_2$ to get
    \begin{align*}
        \M_q(t) \E^q_{\lambda, z}(t) - \M_q(t_2) \E^q_{\lambda, z}(t_2) \le c \left(M + \lVert z \rVert^2 \right) \int_{t_2}^t \M_q(s) s^{-q} ds,
    \end{align*}
   thus, for all $t \geq t_2$ it holds
    \begin{align}
    \label{eq:E_z_t_bound2}
        \E^q_{\lambda, z}(t) \le \frac{\M_q(t_2)\E^q_{\lambda, z}(t_2)}{\M_q(t)} + c \left(M + \lVert z \rVert^2 \right) \frac{C_{\M_q}}{\M_q(t)} \int_{t_2}^t \exp\left( \frac{\lambda}{1-q} s^{1 - q}\right)  s^{- 3q} ds.
    \end{align}
    The inequality above holds for all $z \in \H$ and all $t \ge t_2$. For all $t \ge t_2$, we choose 
    $$z :=z(t) = \argmin_{z \in \H} \max_{i=1,\dots,m} f_{t,i}(z) - f_{t,i}(x(t)),$$ 
    which, since $\E^q_{\lambda}(t) = \E^q_{\lambda, z(t)}(t)$, yields
    \begin{align*}
        \E^q_\lambda(t) \le \frac{\M_q(t_2)\E^q_{\lambda, z(t)}(t_2)}{\M_q(t)} + c \left(M + \lVert z(t) \rVert^2 \right) \frac{C_{\M_q}}{\M_q(t)} \int_{t_2}^t \exp\left( \frac{\lambda}{1-q} s^{1 - q}\right) s^{- 3q} ds.
    \end{align*}
    By Proposition \ref{prop:z_bounded}, $z(\cdot)$ is bounded, and hence there exist constants $C_1, C_2 > 0$ such that for all $t \geq t_2$
    \begin{align}
    \label{eq:E_z_t_bound23}
        \E^q_\lambda(t) \le \frac{C_1}{\M_q(t)} + \frac{C_2}{\M_q(t)} \int_{t_2}^t \exp\left( \frac{\lambda}{1-q} s^{1 - q}\right) s^{- 3q} ds.
    \end{align}
    We apply Lemma \ref{lem:int_poly_exp_bound} to the integral in \eqref{eq:E_z_t_bound23} to derive the asymptotic bound
    \begin{align*}
        \int_{t_2}^t \exp\left( \frac{\lambda}{1-q} s^{1 - q}\right)s^{-3q} ds = \mathcal{O}\left( t^{-2q} \exp\left( \frac{\lambda}{1 - q} t^{1 - q}\right)\right) \quad \text{ as } t \to + \infty,
    \end{align*}
hence
    \begin{align}
    \label{eq:intergral_bound12}
        \frac{C_2}{\M_q(t)} \int_{t_2}^t \exp\left( \frac{\lambda}{1 - q} s^{1 - q}\right) s^{ - 2q } ds = \mathcal{O}\left( 1 \right) \quad \text{ as } t \to + \infty.
    \end{align}
We conclude from \eqref{eq:E_z_t_bound23} and \eqref{eq:intergral_bound12} that
    \begin{align}
    \label{eq:asymptotic_bound_E(t)2}
        \E^q_\lambda(t) = \mathcal{O}\left( 1 \right) \quad \text{ as } t \to + \infty,
    \end{align}
    proving statement \emph{i)}. From here, we can prove the remaining four statements of the theorem.
    
    \emph{ii)} By the choice of $0 < \lambda < \frac{\alpha}{2}$, we have for all $t \ge t_0$
    \begin{align*}
        qt^{q-1} + \alpha - 2\lambda \ge 0.
    \end{align*}
    Then, by the definition of $\E^q_\lambda(\cdot)$ we have for all $t \ge t_0$
    \begin{align*}  
        t^{2q} \varphi_t(x(t)) \le \E^q_\lambda(t),
    \end{align*}
        which, according to \eqref{eq:asymptotic_bound_E(t)2}, gives
    \begin{align*}
        \varphi_t(x(t)) = \mathcal{O}\left( t^{-2q} \right) \quad \text{ as } t \to + \infty.
    \end{align*}
    
    \emph{iii)} Using Proposition \ref{prop:basic_properties_tikh} and \emph{ii)} yields
    \begin{align*}
        \varphi(x(t)) \le \varphi_t(x(t)) + \frac{\beta R^2}{2t^p} = \mathcal{O}\left( t^{-2q} \right) \quad \text{ as } t \to + \infty.
    \end{align*}
    
    \emph{iv)} Since for all $t \ge t_0$
    \begin{align*}
        q t^{q - 1} + \alpha - 2 \lambda \ge \alpha - 2\lambda > 0,
    \end{align*}
it holds
    \begin{align*}
        \frac{\lambda}{2} (\alpha - 2\lambda) \lVert x(t) - z(t) \rVert^2 \le \E^q_\lambda(t).
    \end{align*}
    This estimate together with \eqref{eq:asymptotic_bound_E(t)2} implies that
    \begin{align}
    \label{eq:x(t)_z(t)_bound12}
        \lVert x(t) - z(t) \rVert = \mathcal{O}\left( 1 \right) \quad \text{ as } t \to + \infty.
    \end{align}
    
    \emph{v)} From \emph{i)} and \emph{iv)}, we have
    \begin{align*}
        \frac{t^{2q}}{2}\lVert \dot{x}(t) \rVert^2 \le & \ \lVert \lambda(x(t) - z(t)) + t^q \dot{x}(t) \rVert^2 + \lambda^2 \lVert x(t) - z(t) \rVert^2 \\
        \le & \ 2 \E^q_\lambda(t) + \lambda^2 \lVert x(t) - z(t) \rVert^2 = \mathcal{O}\left( 1 \right) \quad \text{ as }\, t \to + \infty.
    \end{align*}
    From here, we conclude 
    \begin{align*}
        \lVert \dot{x}(t) \rVert = \mathcal{O}\left( t^{- q} \right) \quad \text{ as }\, t \to + \infty.
    \end{align*}
\end{proof}


\subsection{The case $q \in (0,1)$ and $p < q + 1$ : convergence rates and strong convergence of the trajectories}

In this section, we perform the asymptotic analysis for \eqref{eq:MTRIGS} in case $p < q + 1$.

\begin{theorem}
\label{thm:main_convergence_special_p<q+1}
    Let $q \in (0,1)$ and $p < q + 1$, $x(\cdot)$ be a trajectory solution of \eqref{eq:MTRIGS}, and $z(t) := \argmin_{z \in \H} \max_{i=1,\dots,m} f_{t,i}(z) - f_{t,i}(x(t))$ for $t \geq t_0$. Then, for $r \in [q, 1) \cap [p - q, 1)$, we have the following convergence rates as $t \rightarrow +\infty$:
    \begin{enumerate}[i)]
        \item $\E^r_{\lambda}(t) = \mathcal{O}\left( t^{3r - (p + 1)} \right)$ for $\lambda \in \left( 0 ,\frac{2\alpha}{3} \right] \cap \left( 0 , \frac{\beta}{\alpha} \right]$;
        \item $\varphi_t(x(t)) = \mathcal{O}\left( t^{r - (p + 1)} \right)$;
        \item $\varphi(x(t)) = \mathcal{O}\left( t^{-p} \right);$
        \item $\lVert x(t) - z(t) \Vert = \mathcal{O}\left( t^{\frac{r - 1}{2}} \right)$;
        \item $\lVert \dot{x}(t) \Vert  = \mathcal{O}\left(t^{\frac{r - (p + 1)}{2}} \right)$.
    \end{enumerate}
\end{theorem}

\begin{proof}
    \emph{i)} Let $r \in [q, 1) \cap [p - q, 1)$ and $z \in \H$ fixed. From \eqref{eq:d_dt_E(t)_mu(t)E(t)}, we have for almost all $t \ge \max\left(\left( \frac{2r}{\lambda} \right)^{\frac{1}{1-r}}, t_0 \right)$
    \begin{align}
    \label{eq:d_dt_E(t)_mu(t)E(t)_p<q+1}
        \frac{d}{dt}\E^r_{\lambda, z}(t) + \mu_r(t) \E^r_{\lambda, z}(t) \le t^{r}\left( \frac{3}{2}\lambda - \alpha t^{r - q} \right) \lVert \dot{x}(t) \rVert^2 + \frac{p \beta t^{2r}}{2t^{p+1}} \lVert z \rVert^2 + \frac{\lambda}{2} \left[ \frac{3 \lambda r}{t} - \frac{\lambda^2}{t^r} + \frac{ \lambda \alpha}{t^q} - \frac{\beta}{t^{p-r}}\right] \lVert x(t) - z \rVert^2.
    \end{align}
    Since $r < 1$, and $p - r \le q$, $\lambda \le \frac{\beta}{\alpha}$, and $r - q \ge 0$, $\lambda \le \frac{2\alpha}{3}$ there exists $t_2 \ge \max\left(\left( \frac{2r}{\lambda} \right)^{\frac{1}{1-r}}, t_0 \right)$ such that for almost all $t \ge t_2$
    \begin{align}
    \label{eq:bound_d_dt_E_z(t)+mu(t)E_z(t)_32}
        \frac{d}{dt}\E^r_{\lambda, z}(t) + \mu_q(t) \E^r_{\lambda, z}(t) \le \frac{p \beta t^{2r}}{2t^{p+1}} \lVert z \rVert^2.
    \end{align}
    As before, we define the function
    \begin{align}
    \label{eq:def_M(t)2}
        \M_r:[t_2, + \infty)\to \R, \quad t \mapsto \M_r(t) \coloneqq \exp\left(\int_{t_2}^t \mu_r(s) ds\right) = \exp \left( \int_{t_1}^t \frac{\lambda}{s^r} - \frac{2r}{s} ds \right) = C_{\M_r} \frac{\exp\left( \frac{\lambda}{1 - r} t^{1 - r}\right)}{t^{2r}},
    \end{align}
    with $C_{\M_r} = \frac{t_2^{2r}}{\exp\left( \frac{\lambda}{1 - r} t_2^{1 - r}\right)} > 0$. The function $\M_r(\cdot)$ is constructed such that $\frac{d}{dt} \M_r(t) = \M_r(t) \mu_r(t)$ and hence
    \begin{align}
    \label{eq:d_dt(M(t)E_z(t))2c}
        \frac{d}{dt}\left( \M_r(t) \E^r_{\lambda, z}(t)\right) = \M_r(t) \left( \frac{d}{dt}\E^r_{\lambda, z}(t) + \mu_r(t) \E^r_{\lambda, z}(t) \right) \ \mbox{for almost all} \ t \geq t_2.
    \end{align}
    The relations \eqref{eq:d_dt(M(t)E_z(t))2c} and \eqref{eq:bound_d_dt_E_z(t)+mu(t)E_z(t)_32} give for almost all $t \ge t_2$
    \begin{align}
    \label{eq:d_dt(M(t)E_z(t))3}
        \frac{d}{dt}\left( \M_r(t) \E^r_{\lambda, z}(t)\right) \le  \, \frac{p\beta}{2} \lVert z \rVert^2 \M_r(t)  t^{2r - (p + 1)} ,
    \end{align}
    We integrate \eqref{eq:d_dt(M(t)E_z(t))3} from $t_2$ to $t \geq t_2$ to get
    \begin{align*}
        \M_r(t) \E^r_{\lambda, z}(t) - \M_r(t_2) \E^r_{\lambda, z}(t_2) \le  \frac{p \beta }{2} \lVert z \Vert^2 \int_{t_2}^t \M_r(s) s^{2r - (p + 1)} ds,
    \end{align*}
   thus, for all $t \geq t_2$ it holds
    \begin{align}
    \label{eq:E_z_t_bound}
        \E^r_{\lambda, z}(t) \le \frac{\M_r(t_2)\E^r_{\lambda, z}(t_2)}{\M_r(t)} + \frac{p \beta }{2} \lVert z \Vert^2 \frac{C_{\M_r}}{\M_r(t)} \int_{t_2}^t \exp\left( \frac{\lambda}{1-r} s^{1 - r}\right)  s^{- (p + 1)} ds.
    \end{align}
    The inequality above holds for all $z \in \H$ and all $t \ge t_2$. For all $t \ge t_2$, we choose 
    $$z :=z(t) = \argmin_{z \in \H} \max_{i=1,\dots,m} f_{t,i}(z) - f_{t,i}(x(t)),$$ 
    which, since $\E^r_{\lambda}(t) = \E^r_{\lambda, z(t)}(t)$, yields
    \begin{align*}
        \E^r_{\lambda}(t) \le \frac{\M_r(t_2)\E^r_{\lambda, z(t)}(t_2)}{\M_r(t)} + \frac{p \beta }{2} \lVert z(t) \Vert^2 \frac{C_{\M_r}}{\M_r(t)} \int_{t_2}^t \exp\left( \frac{\lambda}{1-r} s^{1 - r}\right)  s^{- (p + 1)} ds.
    \end{align*}
    By Proposition \ref{prop:z_bounded}, $z(\cdot)$ is bounded, hence there exist constants $C_1, C_2 > 0$ such that for all $t \geq t_2$
    \begin{align}
    \label{eq:E_z_t_bound3}
        \E^r_{\lambda}(t) \le \frac{C_1}{\M_r(t)} + \frac{C_2}{\M_r(t)} \int_{t_2}^t \exp\left( \frac{\lambda}{1-r} s^{1 - r}\right) s^{- (p + 1)} ds.
    \end{align}
    We apply Lemma \ref{lem:int_poly_exp_bound} to the integral in \eqref{eq:E_z_t_bound3} to derive the asymptotic bound
    \begin{align*}
        \int_{t_2}^t \exp\left( \frac{\lambda}{1 - r} s^{1 - r}\right) s^{- (p + 1)} ds = \mathcal{O}\left( t^{r - (p + 1)}  \exp\left( \frac{\lambda}{1 - r} t^{1 - r}\right)\right) \quad \text{ as } t \to + \infty,
    \end{align*}
    hence
    \begin{align}
    \label{eq:intergral_bound123}
        \frac{C_2}{\M_r(t)} \int_{t_2}^t \exp\left( \frac{\lambda}{1 - r} s^{1 - r}\right) s^{- (p + 1)} ds = \mathcal{O}\left( t^{3r - (p + 1)} \right) \quad \text{ as } t \to + \infty.
    \end{align}
    We conclude from \eqref{eq:E_z_t_bound3} and \eqref{eq:intergral_bound123} that
    \begin{align}
    \label{eq:asymptotic_bound_E(t)3}
        \E^r_\lambda(t) = \mathcal{O}\left( t^{3r - (p + 1)} \right) \quad \text{ as } t \to + \infty,
    \end{align}
    proving statement \emph{i)}. From here, we can prove the other four statements of the theorem.\\

    \emph{ii)} If $r > q$, for $t \ge \left( \frac{2 \lambda}{\alpha} \right)^{\frac{1}{r - q}}$ we have $ rt^{r-1} + \alpha t^{r-q} - 2\lambda \ge 0$ and hence 
    \begin{align}
    \label{eq:varphi_t_x(t)_bound}
        t^{2r} \varphi_t(x(t)) \le \E^r_\lambda(t).
    \end{align}
    For the case $r = q$ the argument follows in a similar manner. We apply part \emph{i)} for $\lambda \in \left( 0 ,\frac{\alpha}{2} \right) \cap \left( 0 , \frac{\beta}{\alpha} \right] \subseteq \left( 0 ,\frac{2\alpha}{3} \right] \cap \left( 0 , \frac{\beta}{\alpha} \right]$. Then $qt^{q-1} + \alpha - 2\lambda \ge 0$ for all $t \ge t_0$ and hence
     \begin{align}
    \label{eq:varphi_t_x(t)_bound2}
        t^{2q} \varphi_t(x(t)) \le \E^q_\lambda(t).
    \end{align}
    Both cases, together with \eqref{eq:asymptotic_bound_E(t)3}, imply that for all $r \in [q, 1) \cap [p - q, 1)$
    \begin{align*}
        \varphi_t(x(t)) = \mathcal{O}\left( t^{r - (p + 1)} \right) \quad \text{ as } t \to + \infty.
    \end{align*}

    \emph{iii)} Using Proposition \ref{prop:basic_properties_tikh} and \emph{ii)} yields
    \begin{align*}
        \varphi(x(t)) \le \varphi_t(x(t)) + \frac{\beta R^2}{2t^p} = \mathcal{O}\left( t^{-p} \right) \quad \text{ as } t \to + \infty.
    \end{align*}

    \emph{iv)} By Proposition \ref{prop:basic_properties_tikh}, we have for all $t \ge t_0$
    \begin{align*}
        \lVert x(t) - z(t) \rVert^2 \le \frac{2 t^p}{\beta} \varphi_t(x(t)),
    \end{align*}
    and hence by \emph{ii)} we get
    \begin{align}
    \label{eq:x(t)_z(t)_bound223}
        \lVert x(t) - z(t) \rVert = \mathcal{O}\left( t^\frac{r - 1}{2} \right) \quad \text{ as } t \to + \infty.
    \end{align}

    \emph{v)} From the above considerations, we have
    \begin{align*}
        \frac{t^{2r}}{2}\lVert \dot{x}(t) \rVert^2 \le & \ \lVert \lambda(x(t) - z(t)) + t^r \dot{x}(t) \rVert^2 + \lambda^2 \lVert x(t) - z(t) \rVert^2 \\
        \le & \ 2 \E^r_\lambda(t) + \lambda^2 \lVert x(t) - z(t) \rVert^2 = \mathcal{O}\left( t^{3r - (p + 1)} \right) \quad \text{ as }\, t \to + \infty.
    \end{align*}
    From here, we conclude 
    \begin{align*}
        \lVert \dot{x}(t) \rVert = \mathcal{O}\left( t^{\frac{r - (p + 1)}{2}} \right) \quad \text{ as }\, t \to + \infty.
    \end{align*}
\end{proof}

For this parameter settings, alongside establishing convergence rates, we demonstrate that the bounded trajectory solutions of \eqref{eq:MTRIGS} strongly converge to a weak Pareto optimal point of \eqref{eq:MOP}. Notably, this point is also the element of minimum norm within the lower level set of the objective function with respect to its value at the weak Pareto optimal point.

\begin{theorem}
\label{thm:strong_convergence}
    Let $q \in (0,1)$, $p < q + 1$, and $x(\cdot)$ be a bounded trajectory solution of \eqref{eq:MTRIGS}. Then, $x(t)$ converges strongly to a weak Pareto optimal point $x^*$ of \eqref{eq:MOP} as $t \to +\infty$, which is the element of minimum norm in $\bigcap_{i=1}^m {\cal L}(f_i, f_i(x^*))$.
\end{theorem}

\begin{proof}
To prove the strong convergence of the trajectory solution $x(\cdot)$ we use Theorem \ref{thm:moo_tikh_conv}, which states that $z(\cdot)$ converges strongly, in combination with Theorem \ref{thm:main_convergence_special_p<q+1}  \emph{iv)}, which states that $\lVert x(t) - z(t) \rVert \to 0$ as $t \to + \infty$.  Since $x(\cdot)$ is bounded, it holds $\inf_{t \ge t_0} f_i(x(t)) > - \infty$ for $i=1, \ldots, m$, and so
    \begin{align*}
        \inf_{t \ge t_0} \W_i(t) = \inf_{t \ge t_0} \left( f_i(x(t)) + \frac{\beta}{2t^p}\lVert x(t) \rVert^2 + \frac{1}{2}\lVert \dot{x}(t) \rVert^2 \right) \ge \inf_{t \ge t_0} f_i(x(t)) > - \infty,
    \end{align*}
where $\W_i(\cdot)$ is the function introduced in \eqref{eq:def_W_i_t}. By Proposition \ref{prop:ineq_W_i_t}, the function $\W_i(\cdot)$ is monotonically decreasing and therefore, $\lim_{t\to + \infty} \W_i(t)$ exists for $i=1, \ldots, m$. According to Theorem \ref{thm:main_convergence_special_p<q+1}, $\dot{x}(t) \to 0$, hence $\frac{\beta}{2t^p}\lVert x(t) \Vert^2 + \frac{1}{2}\lVert \dot{x}(t) \Vert^2 \to 0$ as $t \to  +\infty$. Thus, for $i = 1,\dots,m$,
    \begin{align*}
        \lim_{t \to +\infty} f_i(x(t)) = \lim_{t \to +\infty} \W_i(t) =  \inf_{t \ge t_0} \W_i(t) > - \infty.
    \end{align*}
    We denote by $f^* \coloneqq \lim_{t \to +\infty} f(x(t)) = \lim_{t \to +\infty} \left(f_1(x(t)), \dots, f_m(x(t))\right) \in \R^m$. We use Theorem \ref{thm:moo_tikh_conv} with $q(t) := f(x(t))$ to conclude
    \begin{align*}
        z(t) \to x^*:=\proj_{S(f^*)}(0) \ \mbox{as} \ t \rightarrow +\infty,
    \end{align*}
    where $z(t) := \argmin_{z \in \H} \max_{i=1,\dots,m} f_{t,i}(z) - f_{t,i}(x(t))$ and $S(f^*) := \argmin_{z \in \H} \max_{i=1,\dots,m} \left( f_i(z) - f_i^* \right)$. According to Theorem \ref{thm:main_convergence_special_p<q+1}, we have $\lVert x(t) - z(t) \rVert \to 0$, hence
    \begin{align*}
        x(t) \to x^* \ \mbox{as} \ t \rightarrow +\infty.
    \end{align*}
Since $\varphi(x(t)) \to 0$ as $t \to +\infty$, it yields $\varphi(x^*) = 0$, thus $x^*$ is a weak Pareto optimal point of \eqref{eq:MOP}.
    By continuity, $f^* = f(x^*)$ and, since $x^*$ is a weak Pareto optimal solution of \eqref{eq:MOP}, it holds $S(f^*) = \bigcap_{i=1}^m {\cal L}(f_i, f_i(x^*))$.
\end{proof}


\subsection{The case $p\in (0,2]$  and $q = 1$}

In this subsection, we consider the boundary case $q = 1$, allowing $p$ to be chosen in $(0,2]$. The assumption we make for $\alpha$ is consistent with that made in the setting of inertial dynamics with vanishing damping in the single objective case, see \cite{Su2016, Attouch2018}.

\begin{theorem}
\label{thm:main_convergence_special_p<q+1_r=1}
    Let $p \in (0,2],\, q = 1$ and $\alpha \ge 3$, $x(\cdot)$ be a bounded trajectory solution of \eqref{eq:MTRIGS}, and $z(t) := \argmin_{z \in \H} \max_{i=1,\dots,m} f_{t,i}(z) - f_{t,i}(x(t))$ for $t \geq t_0$. Then, we have the following convergence rates as $t \rightarrow +\infty$:
    \begin{enumerate}[i)]
        \item $\E^1_{\lambda}(t) = \mathcal{O}\left( t^{2 - p} \right)$ for $\lambda \in \left[2 , \frac{2\alpha}{3} \right]$;
        \item $\varphi_t(x(t)) = \mathcal{O}\left( t^{- p} \right)$;
        \item $\varphi(x(t)) = \mathcal{O}\left( t^{-p} \right);$
        \item $\lVert x(t) - z(t) \Vert = \mathcal{O}\left( 1 \right)$;
        \item $\lVert \dot{x}(t) \Vert  = \mathcal{O}\left( t^{-\frac{p}{2}} \right)$.
    \end{enumerate}
\end{theorem}

\begin{proof}
    \emph{i)} Let $r = q = 1$ and $z \in \H$ fixed. We consider the  energy function $\E^r_{\lambda, z}(\cdot)$. From inequality \eqref{eq:d_dt_E(t)_mu(t)E(t)_r=1} we get for almost all $t \ge t_0$
    \begin{align}
    \label{eq:d_dt_E(t)_mu(t)E(t)_r=1_2}
        \frac{d}{dt}\E^1_{\lambda, z}(t) + \mu_1(t) \E^1_{\lambda, z}(t) \le & \ t \left( \frac{3}{2}\lambda - \alpha \right) \lVert \dot{x}(t) \rVert^2 + \frac{p \beta}{2t^{p - 1}} \lVert z \rVert^2 + \frac{\lambda}{2} \left[ \frac{\alpha(\lambda - 2)}{t} - \frac{\beta}{t^{p-1}}\right] \lVert x(t) - z \rVert^2.
    \end{align}
    Since $p - 1 \le 1$, $\lambda \le \frac{2\alpha}{3}$ and $x(\cdot)$ is bounded, there exist $t_1 \ge t_0$ and $M, c > 0$ such that for almost all $t \ge t_1$
    \begin{align}
    \label{eq:bound_d_dt_E_z(t)+mu(t)E_z(t)_32_r=1}
        \frac{d}{dt}\E^1_{\lambda, z}(t) + \mu_1(t) \E^1_{\lambda, z}(t) \le \frac{c}{2t^{p-1}} \left( M + \lVert z \rVert^2 \right).
    \end{align}
    As before, we define the function
    \begin{align}
    \label{eq:def_M(t)2_r=1}
        \M_1:[t_1, + \infty)\to \R, \quad t \mapsto \M_1(t) \coloneqq \exp\left(\int_{t_1}^t \mu_1(s) ds\right) = \exp \left( \int_{t_1}^t \frac{\lambda - 2}{s} ds \right) = C_{\M_1} t^{\lambda - 2},
    \end{align}
    with $C_{\M_1} = t_1^{2 - \lambda}$. The function $\M_1(\cdot)$ is constructed such that $\frac{d}{dt} \M_1(t) = \M_1(t) \mu_1(t)$, hence
    \begin{align}
    \label{eq:d_dt(M(t)E_z(t))2c_r=1}
        \frac{d}{dt}\left( \M_1(t) \E^1_{\lambda, z}(t)\right) = \M_1(t) \left( \frac{d}{dt}\E^1_{\lambda, z}(t) + \mu_1(t) \E^1_{\lambda, z}(t) \right) \ \mbox{for almost all} \ t \geq t_1.
    \end{align}
   The relations \eqref{eq:d_dt(M(t)E_z(t))2c_r=1} and \eqref{eq:bound_d_dt_E_z(t)+mu(t)E_z(t)_32_r=1} give for almost all $t \ge t_1$
    \begin{align}
    \label{eq:d_dt(M(t)E_z(t))3_r=1}
        \frac{d}{dt}\left( \M_1(t) \E^1_{\lambda, z}(t)\right) \le  \, \frac{c}{2}\left( M + \lVert z \rVert^2 \right) \M_1(t)  t^{1 - p}.
    \end{align}
    We integrate \eqref{eq:d_dt(M(t)E_z(t))3_r=1} from $t_1$ to $t \geq t_1$ to get
    \begin{align*}
        \M_1(t) \E^1_{\lambda, z}(t) - \M_1(t_1) \E^1_{\lambda, z}(t_1) \le  \frac{c }{2} \left( M + \lVert z \rVert^2 \right) \int_{t_1}^t \M_1(s) s^{1 - p} ds,
    \end{align*}
   thus, for all $t \geq t_1$ it holds
    \begin{align}
    \label{eq:E_z_t_bound_r=1}
        \E^1_{\lambda, z}(t) \le \frac{\M_1(t_1)\E^1_{\lambda, z}(t_1)}{\M_1(t)} + \frac{c }{2} \left( M + \lVert z \rVert^2 \right) \frac{C_{\M_1}}{\M_1(t)} \int_{t_1}^t s^{\lambda - (p + 1)} ds.
    \end{align}
    The inequality above holds for all $z \in \H$ and all $t \ge t_1$. For all $t \ge t_1$, we choose 
    $$z :=z(t) = \argmin_{z \in \H} \max_{i=1,\dots,m} f_{t,i}(z) - f_{t,i}(x(t)),$$ 
    which, since $\E^1_{\lambda}(t) = \E^1_{\lambda, z(t)}(t)$, yields
    \begin{align*}
        \E^1_{\lambda}(t) \le \frac{\M_1(t_1)\E^1_{\lambda, z(t)}(t_1)}{C_{\M_1} t^{\lambda - 2}} + \frac{c }{2 t^{\lambda - 2}}\left( M + \lVert z(t) \rVert^2 \right) \left[ \frac{t^{\lambda - p}}{\lambda - p} - \frac{t_1^{\lambda - p}}{\lambda - p}\right].
    \end{align*}
    By Proposition \ref{prop:z_bounded}, $z(\cdot)$ is bounded, which means that there exist constants $C_1, C_2 > 0$ such that for all $t \geq t_1$
    \begin{align}
    \label{eq:E_z_t_bound3_r=1}
        \E^1_{\lambda}(t) \le C_1 + C_2 t^{2 - p},
    \end{align}
    hence
    \begin{align}
    \label{eq:asymptotic_bound_E(t)3_r=1}
        \E^1_\lambda(t) = \mathcal{O}\left( t^{2 - p} \right) \quad \text{ as } t \to + \infty,
    \end{align}
    proving statement \emph{i)}. From here, the remaining four statements of the theorem follow as in the proof of Theorem \ref{thm:main_convergence_special_p<q+1}.
\end{proof}

\begin{myremark}
\label{rem:boundedness_pin(0,2],q=1,alpha>=3}
    If we choose $\lambda = 2$ in the proof of Theorem \ref{thm:main_convergence_special_p<q+1_r=1} we do not need to assume the boundedness of $x(\cdot)$ to conclude \eqref{eq:bound_d_dt_E_z(t)+mu(t)E_z(t)_32_r=1} from \eqref{eq:d_dt_E(t)_mu(t)E(t)_r=1_2}. This implies that in the case $q = 1$ and $\alpha \ge 3$ the bound $\lVert x(t) - z(t) \rVert = \mathcal{O}(1)$ as $t \to + \infty$ follows without the boundedness assumption on $x(\cdot)$.
\end{myremark}


\subsection{The case $p \in (0,2]$ and $q + 1 < p$ : weak convergence of the trajectories}

In this section, we show that in the case $p \in (0,2]$ and $q + 1 < p$ the bounded trajectory solutions of \eqref{eq:MTRIGS} converge weakly to a weak Pareto optimal point of \eqref{eq:MOP}. To this end, we make use of Opial's Lemma and the energy function from Definition \ref{def:mo_energy_functions5} with $\gamma(\cdot)$ and $\xi(\cdot)$ to be specified later.
The convergence rates derived in Subsection \ref{subsec41} are valid in this setting.

\begin{theorem}
    \label{thm:intergral_estimate}
    Let $p \in (0,2)$, $q + 1 < p$, and $x(\cdot)$ be a trajectory solution of \eqref{eq:MTRIGS}. Then, for $r \in \left[q, \frac{q+1}{2} \right]$, we have
    \begin{align*}
        \int_{t_0}^{+ \infty} s^{2r - q}\lVert \dot{x}(s) \rVert^2 ds < + \infty.
    \end{align*}
\end{theorem}

\begin{proof}
Let $z \in \H$ fixed. Define
    \begin{align*}
        \gamma:[t_0, + \infty) \to \R, \quad t \mapsto \gamma(t) = 2rt^{r - 1}.
    \end{align*}
With this choice, inequality  \eqref{eq:equation_deriv_G_z_t2} reads for almost all $t \geq t_0$ 
    \begin{align}
    \label{eq:d_dt_Grz(t)_xi}
    \begin{split}
        \frac{d}{dt} \G^r_{\gamma, \xi, z}(t) \le & \ \frac{p \beta t^{2r}}{2t^{p+1}} \lVert z \rVert^2 + \left( 2rt^{r-1}(2rt^{r-1} + rt^{r-1} - \alpha t^{r - q}) + 2r(r-1)t^{2r - 2} + \xi(t) \right) \left\langle x(t) - z , \dot{x}(t)\right\rangle \\
        & +  \left( 4r^2(r-1)t^{2r - 3} + \frac{\xi'(t)}{2} - \beta r t^{2r-1 - p} \right) \lVert x(t) - z \Vert^2 + t^r(2rt^{r - 1} + rt^{r-1} - \alpha t^{r - q}) \lVert \dot{x}(t) \rVert^2 \\
        = & \ \frac{p \beta t^{2r}}{2t^{p+1}} \lVert z \rVert^2 + \left( 2rt^{r-1}(3rt^{r-1} - \alpha t^{r - q}) + 2r(r-1)t^{2r - 2} + \xi(t) \right) \left\langle x(t) - z , \dot{x}(t)\right\rangle \\
        & +  \left( 4r^2(r-1)t^{2r - 3} + \frac{\xi'(t)}{2} - \beta r t^{2r-1 - p} \right) \lVert x(t) - z \Vert^2 + t^r(3rt^{r - 1} - \alpha t^{r - q}) \lVert \dot{x}(t) \rVert^2.
    \end{split}
    \end{align}
Now we choose
    \begin{align*}
        \xi:[t_0, +\infty) \to \R, \quad \xi(t) \coloneqq & \ 2rt^{r-1}(\alpha t^{r-q} - 3rt^{r-1}) + 2r(1 - r)t^{2(r - 1)} =  2\alpha rt^{2r - q -1} + 2r(1 - 4r)t^{2(r - 1)}, 
    \end{align*}
and notice that $\xi'(t) = 2\alpha r (2r - q -1) t^{2r - q -2} + 4r(r-1)(1 - 4r)t^{2r - 3}$ for all $t \geq t_0$.
With this choice, inequality \eqref{eq:d_dt_Grz(t)_xi} simplifies for almost all $t \geq t_0$ to
    \begin{align}
    \label{eq:d_dt_Grz(t)}
    \begin{split}
        \frac{d}{dt} \G^r_{\gamma, \xi, z}(t) \le & \ \frac{p \beta t^{2r}}{2t^{p+1}} \lVert z \rVert^2 + \left( 2r(r-1)(1 - 2r)t^{2r-3} + \alpha r(2r - q- 1)t^{2r-q-2} - \beta rt^{2r - 1- p} \right) \lVert x(t) - z \Vert^2 \\
        & + t^r(3rt^{r - 1} - \alpha t^{r - q}) \lVert \dot{x}(t) \rVert^2.
    \end{split}
    \end{align}
    Since $r \le \frac{q+ 1}{2}$, we conclude from \eqref{eq:d_dt_Grz(t)} that for almost all $t \ge \max \left(\left(\frac{\max(2(r-1)(1 - 2r), 0)}{\beta} \right)^{\frac{1}{2-p}}, t_0 \right)$
    \begin{align}
    \label{eq:d_dt_Grz(t)2}
         \frac{d}{dt} \G^r_{\gamma, \xi, z}(t) \le & t^r(3rt^{r - 1} - \alpha t^{r - q}) \lVert \dot{x}(t) \rVert^2 + \frac{p \beta t^{2r}}{2t^{p+1}} \lVert z \rVert^2.
    \end{align}
    Hence, there exist $t_1 \ge \max \left(\left(\frac{\max(2(r-1)(1 - 2r), 0)}{\beta} \right)^{\frac{1}{2-p}}, t_0 \right)$ and $a, b > 0$ such that for almost all $t \ge t_1$
    \begin{align*}
        \frac{d}{dt} \G^r_{\gamma, \xi, z}(t) \le & -a t^{2r - q}\lVert \dot{x}(t) \rVert^2 + b t^{2r - p - 1} \lVert z \rVert^2,
    \end{align*}
    therefore
    \begin{align*}
        \G^r_{\gamma, \xi, z}(t) - \G^r_{\gamma, \xi, z}(t_1) \le & -a \int_{t_1}^t s^{2r - q}\lVert \dot{x}(s) \rVert^2 ds + b \lVert z \rVert^2 \int_{t_1}^t s^{2r - p - 1} ds \quad \forall t \geq t_1.
    \end{align*}
    Since this holds for all $z \in \H$, we conclude
    \begin{align*}
        \G^r_{\lambda, \xi}(t) - \G^r_{\lambda, \xi, z(t)}(t_1) \le & -a \int_{t_1}^t s^{2r - q}\lVert \dot{x}(s) \rVert^2 ds + b \lVert z(t) \rVert^2 \int_{t_1}^t s^{2r - p - 1} ds \quad \forall t \geq t_1.
    \end{align*}
    For $t \ge \left( \frac{\max(1 - 4r, 0)}{\alpha} \right)^{\frac{1}{1-q}}$, it holds that $\xi(t) \ge 0$ and hence $\G_{\lambda, \xi}^r(t) \ge 0$. Then, for all $t \ge \max\left( \frac{\max(1 - 4r, 0)}{\alpha} , t_1\right)$ 
    \begin{align*}
        a \int_{t_1}^t s^{2r - q}\lVert \dot{x}(s) \rVert^2 ds \le \G^r_{\lambda, \xi, z(t)}(t_1) + b \lVert z(t) \rVert^2 \int_{t_1}^t s^{2r - p - 1} ds.
    \end{align*}
    Since $z(\cdot)$ is bounded by Proposition \ref{prop:z_bounded} and $2r - p - 1 < -1$, the right hand side of the previous inequality is uniformly bounded for all $t \ge \max\left( \left( \frac{1 - 4r}{\alpha} \right)^{\frac{1}{1-q}} , t_1 \right)$, hence
    \begin{align*}
        \int_{t_0}^{+ \infty} s^{2r - q}\lVert \dot{x}(s) \rVert^2 ds < + \infty.
    \end{align*}
\end{proof}

Next, we discuss the boundary case $p = 2$. To derive weak convergence, we need an additional condition on the parameter $\beta > 0$.

\begin{theorem}
\label{thm:intergral_estimate2}
    Let $p = 2,\, q \in (0,1)$, $\beta \ge q(1-q)$, and $x(\cdot)$ be a bounded trajectory solution of \eqref{eq:MTRIGS}. Then, for $r \in \left[q, \frac{1+q}{2}\right]$, we have
     \begin{align}
     \label{eq:integral_estimate}
        \int_{t_0}^{+ \infty} s^{2r - q} \, \lVert \dot{x}(s) \rVert^2 ds < + \infty.
    \end{align}
\end{theorem}

\begin{proof}
    The proof follows analogously to the proof of Theorem \ref{thm:intergral_estimate}, with the difference that in order to conclude \eqref{eq:d_dt_Grz(t)2} from \eqref{eq:d_dt_Grz(t)} the additional inequality
    \begin{align}
    \label{eq:beta_r}
        2(r-1)(1 - 2r) \le \beta,
    \end{align}
    is necessary. Since $r := \frac{q+1}{2}$ satisfies \eqref{eq:beta_r}, it holds  
    \begin{align}
        \int_{t_0}^{+ \infty} s \, \lVert \dot{x}(s) \rVert^2 ds < + \infty,
    \end{align}
which implies that \eqref{eq:integral_estimate} holds for all $r \in \left[ q , \frac{q +1 }{2}\right]$.
\end{proof}

\begin{myremark}
\label{rem:intergral_estimate}
 In both regimes, namely, for $p \in (0,2)$ and $ q + 1<p$, and for $p = 2,\, q\in (0,1)$ and $\beta \ge q(1-q)$, choosing $r := \frac{1 + q}{2}$ we obtain the following integral estimate, which describes the convergence behavior of the velocity of the trajectory
     \begin{align*}
        \int_{t_0}^{+ \infty} s \, \lVert \dot{x}(s) \rVert^2 ds < + \infty.
    \end{align*}
\end{myremark}

We use the integral estimates given in Theorem \ref{thm:intergral_estimate} and in Theorem \ref{thm:intergral_estimate2} to prove the weak convergence of the trajectory solution using Opial's Lemma (see Lemma \ref{lem:opial}). The following two results prove that the first condition in Opial's Lemma is satisfied, while the final weak convergence statement is shown in Theorem \ref{thm:weak_convergence}.

\begin{mylemma}
    \label{lem:conv_f_i(x(t))}
    Let $p \in (0,2]$. Let $q \in (0,1)$, or $q = 1$ and $\alpha \ge 3$, and $x(\cdot)$ be a bounded trajectory solution of \eqref{eq:MTRIGS}. Let $\W_i(\cdot), i=1, ..., m,$ be the energy function defined in Proposition \ref{prop:ineq_W_i_t}. Then, for all $i = 1,\dots, m$, the limit
    \begin{align*}
        f_i^{\infty} \coloneqq \lim_{t \to + \infty} f_i(x(t)) = \lim_{t \to + \infty} \W_i(t) = \inf_{t \ge t_0} \W_i(t) \in \R
    \end{align*}
    exists.
\end{mylemma}

\begin{proof}
    Let $i \in \{1, \dots, m\}$ be fixed. Since $x(\cdot)$ is bounded, $\inf_{t \ge t_0} f_i(x(t)) \in \R$ holds, therefore
    \begin{align}
    \label{eq:inf_W_i_t}
        \inf_{t \ge t_0} \W_i(t) = \inf_{t \ge t_0} \left( f_i(x(t)) + \frac{\beta}{2t^p}\lVert x(t) \rVert^2 + \frac{1}{2}\lVert \dot{x}(t) \rVert^2 \right) \ge \inf_{t \ge t_0} f_i(x(t)) \in \R.
    \end{align}
    By Proposition \ref{prop:ineq_W_i_t}, $\W_i(\cdot)$ is monotonically decreasing, thus
    \begin{align}
    \label{eq:inf_W_i_t2}
        \lim_{t \to + \infty} \W_i(t) = \inf_{t \ge t_0} \W_i(t) > - \infty.
    \end{align}
    By Theorem \ref{thm:main_convergence_special_r=q}, Theorem \ref{thm:main_convergence_special_p<q+1} and Theorem \ref{thm:main_convergence_special_p<q+1_r=1}, it holds $\dot{x}(t) \to 0$ as $t \to + \infty$. Hence, $\frac{\beta}{2t^p}\lVert x(t) \rVert^2 + \frac{1}{2}\lVert \dot{x}(t) \rVert^2 \to 0$ as $t \to + \infty$. Thus
    \begin{align}
    \label{eq:inf_W_i_t3}
        \lim_{t \to + \infty} f_i(x(t)) = \lim_{t \to + \infty} \W_i(t),
    \end{align}
which leads to the desired result.
\end{proof}

\begin{mylemma}
\label{lem:preparation_opial}
    Let $p \in (0,2),\, q \in (0,1)$ with $q + 1 < p$, or $p = 2,\, q \in (0,1)$ and $\beta \ge q(1-q)$, $x(\cdot)$ be a bounded trajectory solution of \eqref{eq:MTRIGS}, and assume that $$S \coloneqq \left\lbrace z \in \H \, : \, f_i(z) \le f_i^{\infty} \, \text{ for }\, i = 1,\dots, m \right\rbrace \neq \emptyset,$$ 
    with $f_i^{\infty} = \lim_{t \to \infty} f_i(x(t)) \in \R$. Then, for all $z \in S$, the limit $\lim_{t \to + \infty} \lVert x(t) - z \rVert$ exists.
\end{mylemma}

\begin{proof}
Let $z \in S$, and define the function
    \begin{align*}
        h_z:[t_0, + \infty) \to \R, \, z \mapsto h_z(t) \coloneqq \frac{1}{2} \lVert x(t) - z \rVert^2.
    \end{align*}
    For almost all $t \ge t_0$ it holds that
    \begin{align}
    \label{eq:deriv_h_z(t)}
    \begin{split}
        h_z'(t) = \langle x(t) - z , \dot{x}(t) \rangle \quad \text{ and } \quad h_z''(t) = \langle x(t) - z , \ddot{x}(t) \rangle + \lVert \dot{x}(t) \rVert^2.
    \end{split}
    \end{align}
    From \eqref{eq:deriv_h_z(t)} and \eqref{eq:nice_theta_equality}, we have for almost all $t \ge t_0$
    \begin{align}
    \label{eq:deriv_h_z(t)2}
    \begin{split}
        h''_z(t) + \frac{\alpha}{t^q} h_z'(t) = & \left\langle \ddot{x}(t) + \frac{\alpha}{t^q} \dot{x}(t), x(t) - z \right\rangle + \lVert \dot{x}(t) \rVert^2, \\
        = &  \left\langle - \sum_{i=1}^m \theta_i(t) \nabla f_i(x(t)) - \frac{\beta}{t^p}x(t), x(t) - z \right\rangle + \lVert \dot{x}(t) \rVert^2,
    \end{split}
    \end{align}
where $\theta(\cdot)$ be the measurable weight function given by Proposition \ref{prop:theta_measurable}. Since $z \in S$, we have for all $i = 1,\dots, m,$ and almost all $t \ge t_0$
    \begin{align*}
        f_i(x(t)) + \frac{\beta}{2t^p} \lVert x(t) \rVert^2 + \frac{1}{2} \lVert \dot{x}(t) \rVert^2 \ge & \ f_i(z) = f_i(z) + \frac{\beta}{2t^p} \lVert z \rVert^2 - \frac{\beta}{2t^p} \lVert z \rVert^2 \\
        \ge & \ f_i(x(t)) + \frac{\beta}{2t^p} \lVert x(t) \rVert^2 + \left\langle \nabla f_i(x(t)) + \frac{\beta}{t^p}x(t), z - x(t) \right\rangle - \frac{\beta}{2t^p} \lVert z \rVert^2,
    \end{align*}
hence 
    \begin{align}
    \label{eq:deriv_h_z(t)3}
        \left\langle \nabla f_i(x(t)) + \frac{\beta}{t^p}x(t), z - x(t) \right\rangle \le \frac{\beta}{2t^p} \lVert z \rVert^2 + \frac{1}{2} \lVert \dot{x}(t) \rVert^2.
    \end{align}
We define function $k: [t_0, + \infty) \to [0, +\infty),\, k(t) \coloneqq \frac{\beta}{2t^p} \lVert z \rVert^2 + \frac{3}{2} \lVert \dot{x}(t) \rVert^2$. By Theorem \ref{thm:intergral_estimate} and Theorem \ref{thm:intergral_estimate2}, we have $\left(t \mapsto t^q \lVert \dot{x}(t) \rVert^2\right) \in L^1\left( [t_0, +\infty) \right)$. On the other hand, since $ q + 1 < p$, we get $\left(t \mapsto\frac{\beta t^q}{2 t^p} \lVert z \rVert^2\right) \in L^1\left( [t_0, +\infty) \right)$,  consequently, $\left(t \mapsto t^q k(t) \right)\in L^1\left( [t_0, +\infty) \right)$. Combining \eqref{eq:deriv_h_z(t)2} and \eqref{eq:deriv_h_z(t)3} gives
    \begin{align*}
        h''_z(t) + \frac{\alpha}{t^q} h_z'(t) \le k(t) \quad \mbox{for almost all} \ t \geq t_0.
    \end{align*}
Now, we can use Lemma \ref{lem:converngence_h} to conclude that the limit
    \begin{align*}
        \lim_{t \to + \infty} \lVert x(t) - z \rVert \ \text{exists.} 
    \end{align*}
\end{proof}

\begin{theorem}
\label{thm:weak_convergence}
    Let $p \in (0,2)$ and $q + 1 < p$, or $p = 2,\, q \in (0,1)$ and $\beta \ge q(1-q)$, and $x(\cdot)$ be a bounded trajectory solution of \eqref{eq:MTRIGS}. Then $x(t)$ converges weakly to a weak Pareto optimal solution of \eqref{eq:MOP} as $t \to + \infty$, which belongs to $\bigcap_{i=1}^m {\cal L}(f_i, f_i^{\infty})$, where $f_i^{\infty} = \lim_{t \to + \infty} f_i(x(t))$ for $i = 1,\dots,m$.
\end{theorem}

\begin{proof}
    We define the set $S \coloneqq \left\lbrace z \in \H \, : \, f_i(z) \le f_i^{\infty} \, \text{ for }\, i = 1,\dots, m \right\rbrace$ as in Lemma \ref{lem:preparation_opial}. Since $x(\cdot)$ is bounded, it possesses a weak sequential cluster point $x^{\infty} \in \H$. This means that there exists a sequence $\{t_k\}_{k \ge 0}$ which converges to $+\infty$ with the property that $x(t_k)$ converges weakly to $x^{\infty}$ as $k \to + \infty$. The functions $f_i$ being weakly lower semicontinuous fulfill for all $i = 1,\dots,m$
    \begin{align*}
        f_i(x^{\infty}) \le \liminf_{k \to + \infty} f_i(x(t_k)) =  \lim_{k \to + \infty} f_i(x(t_k)) = f_i^{\infty},
    \end{align*}
    therefore $x^{\infty} \in S$. We conclude that $S$ is nonempty and all weak sequential cluster points of $x(\cdot)$ belong to $S$. On the other hand, according to Lemma \ref{lem:preparation_opial} we have that $\lim_{t \to + \infty} \lVert x(t) - z \rVert$ exists for all $z \in S$. We can use Opial's Lemma (Lemma \ref{lem:opial}) to conclude that $x(t)$ converges weakly to an element in $S$ for $t \to + \infty$. By Theorem \ref{thm:main_convergence_special_r=q}, $\varphi(x(t)) \to 0$ as $t \to + \infty$, therefore, since $\varphi(\cdot)$ is weakly lower semicontinuous, $\varphi(x^{\infty}) \le \liminf_{k \to + \infty} \varphi(x(t_k)) = 0$. By Theorem \ref{thm:varphi_properties}, $x^{\infty}$ is a weak Pareto optimal solution of \eqref{eq:MOP}. 
\end{proof}


\section{Numerical experiments}\label{sec:numerical_experiments}

In this section, we illustrate the typical behavior of the trajectory solution $x(\cdot)$ of \eqref{eq:MTRIGS} using two example problems. In the first example, presented in Subsection \ref{subsec:comparison_MTRIGS_MAVD}, we show that trajectory solutions $x(\cdot)$ of \eqref{eq:MTRIGS} converge to a weak Pareto optimal point $x^*$, which is the element of minimum norm in  $\bigcap_{i=1}^m {\cal L}(f_i, f_i(x^*))$, whereas those of \eqref{eq:MAVD_intro} may fail to exhibit this behavior. In Subsection \ref{subsec:experiment_p_q}, we analyze the sensitivity of trajectory solutions of \eqref{eq:MTRIGS} with respect to $q \in (0,1]$ and $p \in (0,2]$. We highlight how different parameter choices affect the decay of the merit function values $\varphi(x(t))$ and the asymptotic behavior of the distance $\lVert x(t) - z(t) \rVert$ to the generalized regularization path as $t \to +\infty$.

\begin{figure}[H]
    \centering
    \includegraphics[width=0.5\linewidth]{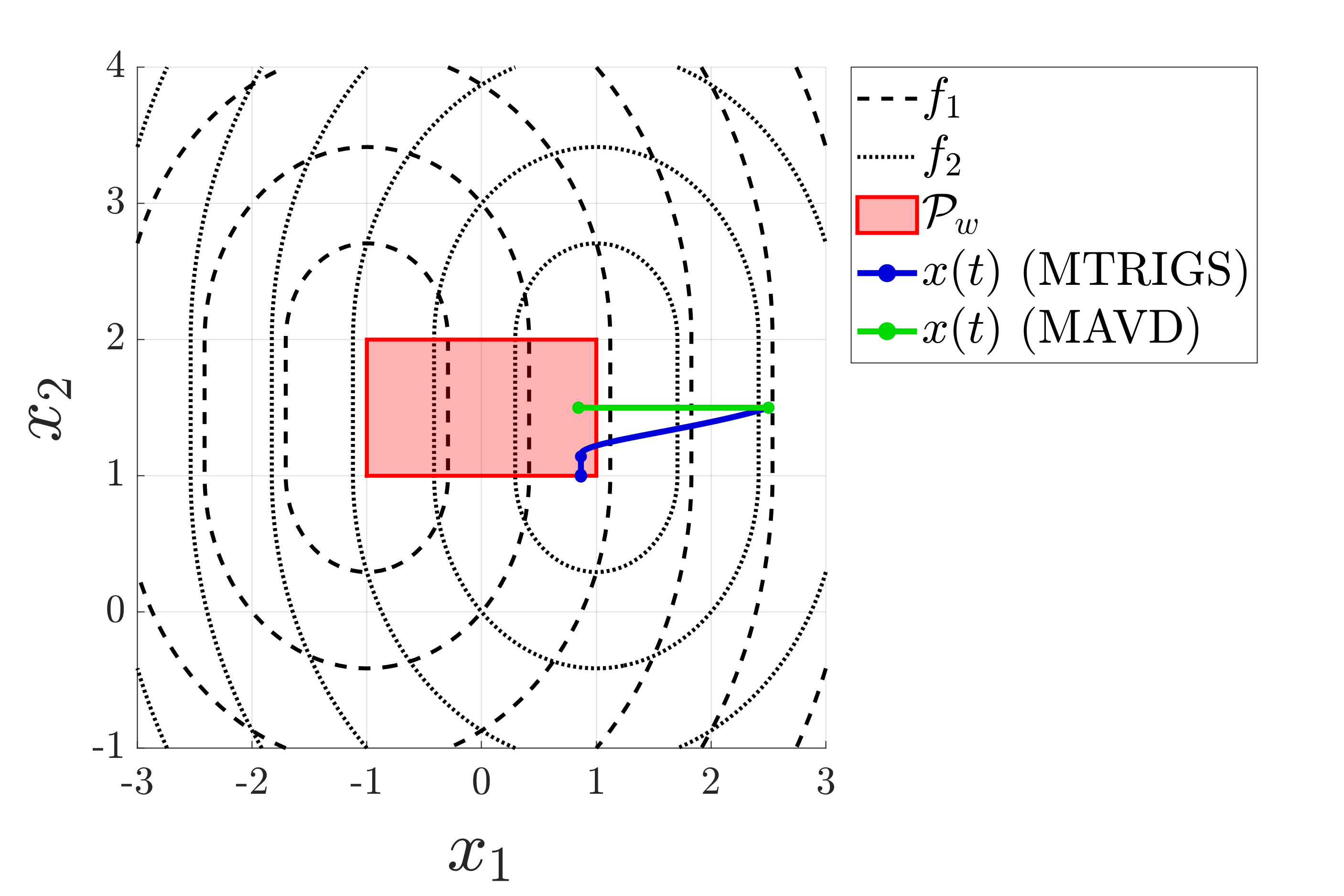}
    \caption{Contour plots of $f_1$ and $f_2$ defined in \eqref{eq:obj_ex1}, the weak Pareto set $\Pw$ of the problem \eqref{eq:MOP_ex1} and the trajectory solutions $x(\cdot)$ of \eqref{eq:MTRIGS} and \eqref{eq:MAVD_intro} with identical initial conditions, respectively.}
    \label{fig:num_exp1_traj}
\end{figure}

\subsection{Comparison of \eqref{eq:MTRIGS} with \eqref{eq:MAVD_intro}}
\label{subsec:comparison_MTRIGS_MAVD}

In the first example, we consider the following instance of \eqref{eq:MOP}. Define the sets
\begin{align*}
    S_1 \coloneqq \{-1\} \times [1,2] \subseteq \R^2 \quad \text{and}\quad S_2 \coloneqq \{1\} \times [1,2] \subseteq \R^2,
\end{align*}
and the functions
\begin{align}
\label{eq:obj_ex1}
    f_i:\R^2 \to \R, \quad x \mapsto f_i(x) \coloneqq \frac{1}{2} \dist(x, S_i)^2,\quad \text{for}\quad i=1,2,
\end{align}
which are both convex and continuously differentiable, and have Lipschitz continuous gradients. The weak Pareto set of the multiobjective optimization problem
\begin{align}
\label{eq:MOP_ex1}
\tag{MOP-Ex$_1$}
    \min_{x \in \R^2} \left[ \begin{array}{c}
         f_1(x)  \\
         f_2(x)
    \end{array} \right]
\end{align}
is given by
\begin{align*}
    \Pw = \conv\left( S_1 \cup S_2 \right) = [-1,1] \times [1,2].
\end{align*}
Let $z = (z_1, z_2)^{\top} \in \Pw$. Then, the element of minimum norm in $\bigcap_{i=1}^2 {\cal L}(f_i, f_i(z))$ is given by
\begin{align}
\label{eq:ex1_min_norm}
    \proj_{\bigcap_{i=1}^2 {\cal L}(f_i, f_i(z))}(0) = (z_1, 1).
\end{align}
We approximate a trajectory solution for \eqref{eq:MTRIGS} and \eqref{eq:MAVD_intro}, respectively, in the following context:
\begin{itemize}
    \item For \eqref{eq:MTRIGS}, we set $\alpha := 4$, $\beta := \frac{1}{2}$, $q := \frac{7}{8}$ and $p := \frac{7}{4}$;
    \item For \eqref{eq:MAVD_intro}, we set $\alpha := 4$;
    \item For both systems, we use as initial conditions $x(t_0) = (2.5,\,0.5)$ and $\dot{x}(t_0) = (0,0)$, where $t_0 = 1$;
    \item For both systems, we use an equidistant discretization in time, i.e., time steps $t_k := t_0 + kh$ with step size $h = \expnumber{1}{-2}$;
    \item For both systems, we approximate the first and second derivatives by $\dot{x}(t_k) = \frac{x(t_{k+1}) - x(t_k)}{h}$ and $\ddot{x}(t_k) = \frac{x(t_{k+1}) - 2x(t_k) + x(t_{k-1})}{h^2}$, respectively;
    \item For both systems, we consider the trajectory solutions for $t \in [1, 100]$.
\end{itemize}

Note that for \eqref{eq:MTRIGS} it holds that $p < q +1$. According to Theorem \ref{thm:main_convergence_special_p<q+1} and Theorem \ref{thm:strong_convergence}, we have convergence of the merit function values $\varphi(x(t)) \to 0$, convergence of the distance of the trajectory to the regularization path $\lVert x(t) - z(t) \rVert \to 0$ and strong convergence of the trajectory $x(t)$ to a weak Pareto optimal point as $t \to + \infty$.

\noindent Figure \ref{fig:num_exp1_traj} shows the contour plots of the objective function $f_1$ and $f_2$ defined in \eqref{eq:obj_ex1}, along with the weak Pareto set $\Pw$ highlighted in red in the decision space. The figure also displays the trajectory solutions of \eqref{eq:MTRIGS} and \eqref{eq:MAVD_intro} with identical initial conditions, respectively, which both converge to points in the weak Pareto set. Notably, the solution of \eqref{eq:MAVD_intro} evolves solely in the $x_1$-direction, whereas the Tikhonov regularization ensures that the solution of \eqref{eq:MTRIGS} converges to an element as specified by \eqref{eq:ex1_min_norm}.

\begin{figure}
    \begin{center}
        \begin{subfigure}[t]{.33\textwidth}
            \centering
            \includegraphics[width=\linewidth]{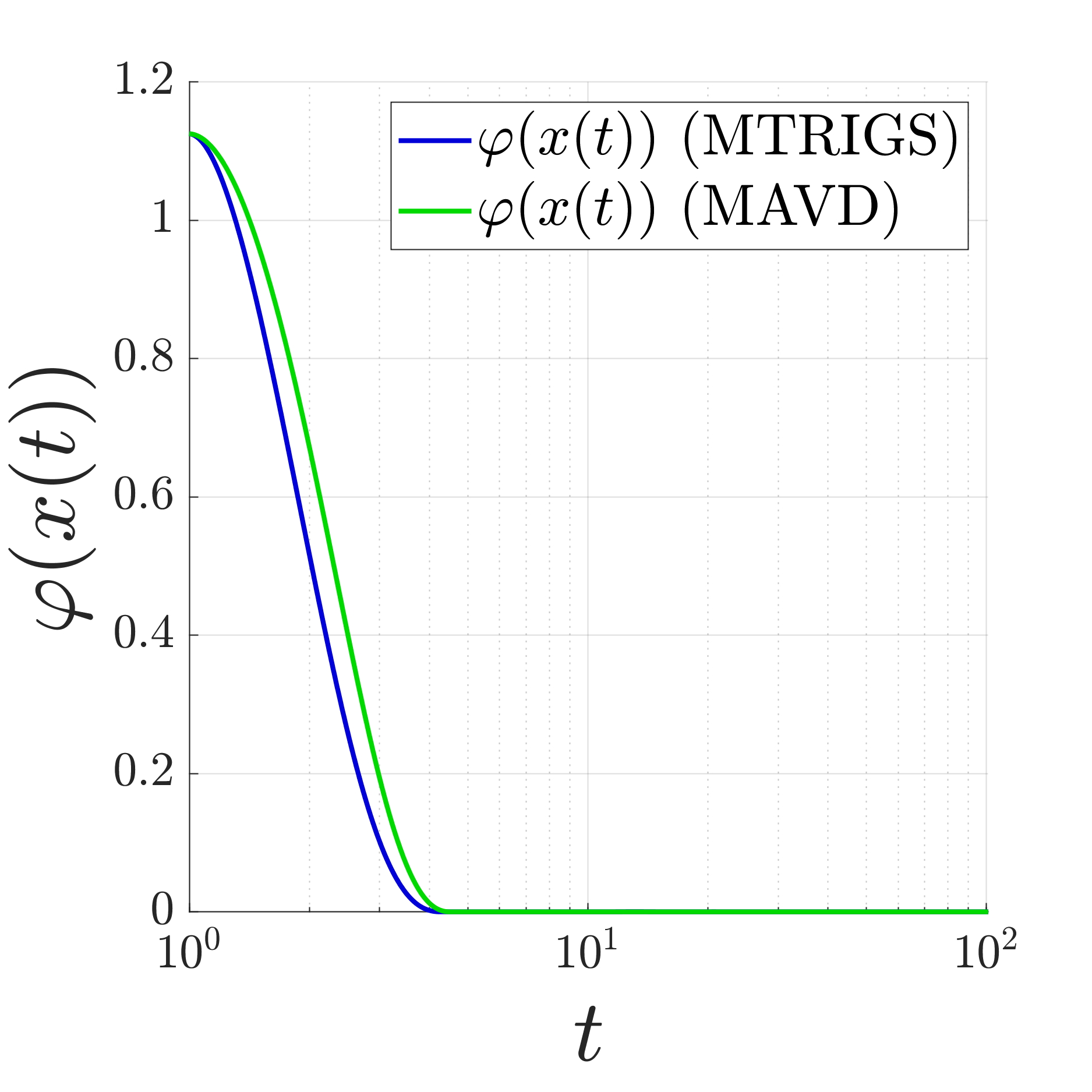}
            \caption{}
            \label{subfig:num_exp1_phi}
        \end{subfigure}
        \hspace{5mm}
        \begin{subfigure}[t]{.33\textwidth}
            \centering
            \includegraphics[width=\linewidth]{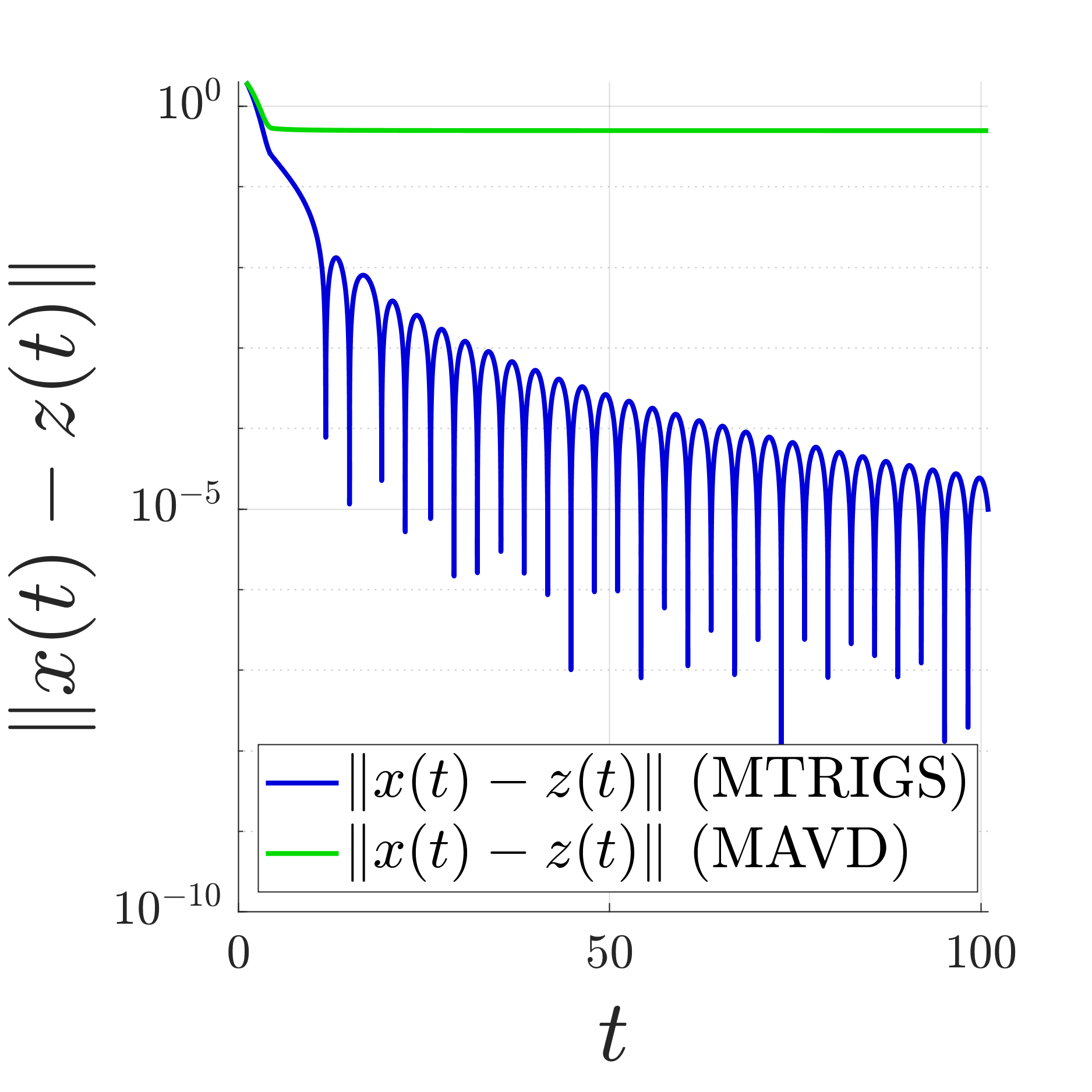}
            \caption{}
            \label{subfig:num_exp1_dist}
        \end{subfigure}    
    \end{center}
    \caption{The merit function values $\varphi(x(t))$ and the distance $\lVert x(t) - z(t) \rVert$ of the trajectory solutions to the generalized regularization path  for \eqref{eq:MTRIGS} and \eqref{eq:MAVD_intro} for the problem \eqref{eq:MOP_ex1}.}
    \label{fig:num_exp1_dist_phi}
\end{figure}

\noindent Figure \ref{fig:num_exp1_dist_phi} visualizes the behavior of the trajectory solutions of \eqref{eq:MTRIGS} and \eqref{eq:MAVD_intro} by showing, in two subfigures, the evolution of the merit function values and the distance of the trajectories to the generalized regularization paths. As already shown in Figure \ref{fig:num_exp1_traj}, the trajectories enter the weak Pareto set $\Pw$ after some time, implying that the merit function values $\varphi(x(t))$ vanish accordingly. This is illustrated in Subfigure \ref{subfig:num_exp1_phi}. Subfigure \ref{subfig:num_exp1_dist} depicts the distance between the trajectory and the generalized regularization path, i.e., $\lVert x(t) - z(t) \rVert$ for $t \in [1,100]$. For the solution of \eqref{eq:MAVD_intro}, this distance converges to a positive limit as $t \to +\infty$. In contrast, for the solution of \eqref{eq:MTRIGS}, the distance decays to zero at a sublinear rate, as predicted by Theorem \ref{thm:main_convergence_special_p<q+1}.

\subsection{The convergence behaviour of \eqref{eq:MTRIGS} for different values of $q \in (0,1]$ and $p \in (0,2]$}
\label{subsec:experiment_p_q}

The numerical experiments in this subsection demonstrate a similar influence of the parameters $q$ and $p$ in on the asymptotic behaviour of \eqref{eq:MTRIGS} as was observed in \cite{Laszlo2023} for the system \eqref{eq:TRIGS_intro} in the context of single objective optimization. Consider
\begin{align*}
    f_1:\R^4 \to \R, \quad & x \mapsto f_1(x) \coloneqq \frac{1}{2} (x_1 - 1)^2 + \frac{1}{2} (x_2 - 1)^2, \quad \text{and}\\[4pt]
    f_2:\R^4 \to \R, \quad & x \mapsto f_1(x) \coloneqq \frac{1}{2} (x_1 + 1)^2 + \frac{1}{2} (x_2 - 1)^2,
\end{align*}
which are both convex and continuously differentiable functions, and have Lipschitz continuous gradients. The weak Pareto set of the multiobjective optimization problem
\begin{align}
\label{eq:MOP_ex2}
\tag{MOP-Ex$_2$}
    \min_{x \in \R^4} \left[ \begin{array}{c}
        f_1(x)  \\
        f_2(x) 
    \end{array}\right]
\end{align}
is given by
\begin{align*}
    \Pw \coloneqq [-1, 1] \times \{1\} \times \R \times \R \subseteq \R^4.
\end{align*}
We approximate a trajectory solution for \eqref{eq:MTRIGS} in the following context:

\begin{itemize}
    \item We set $\alpha := 4$, $\beta := \frac{1}{2}$, and consider different values for $q \in (0,1]$ and $p \in (0,2]$;
    \item We use as initial conditions $x(t_0) = x_0$ and $\dot{x}(t_0) = 0$ with $t_0 = 1$ and $x_0 = (2,3,4,5)^{\top}$;
    \item We use an equidistant discretization in time, i.e., time steps $t_k := t_0 + kh$ with step size $h = \expnumber{1}{-3}$;
    \item We approximate the first and second derivative of $x(\cdot)$ in time by $\dot{x}(t_k) = \frac{x(t_{k+1}) - x(t_k)}{h}$ and $\ddot{x}(t_k) = \frac{x(t_{k+1}) - 2x(t_k) + x(t_{k-1})}{h^2}$ respectively;
    \item We consider the trajectory solutions for $t \in [1, 100]$.
\end{itemize}

\noindent We first fix $q = 0.8$ and vary the parameter $p$ over the set $\{ 0.25,\,0.75,\,1.25,\,1.75 \}$. Afterwards, we fix $p = 1.1$ and vary $q$ over the set $\{ 0.3,\,0.6,\,0.8,\,0.99 \}$.

\begin{figure}[H]
    \begin{center}        
        \begin{subfigure}[t]{.24\textwidth}
            \centering
            \includegraphics[width=\linewidth]{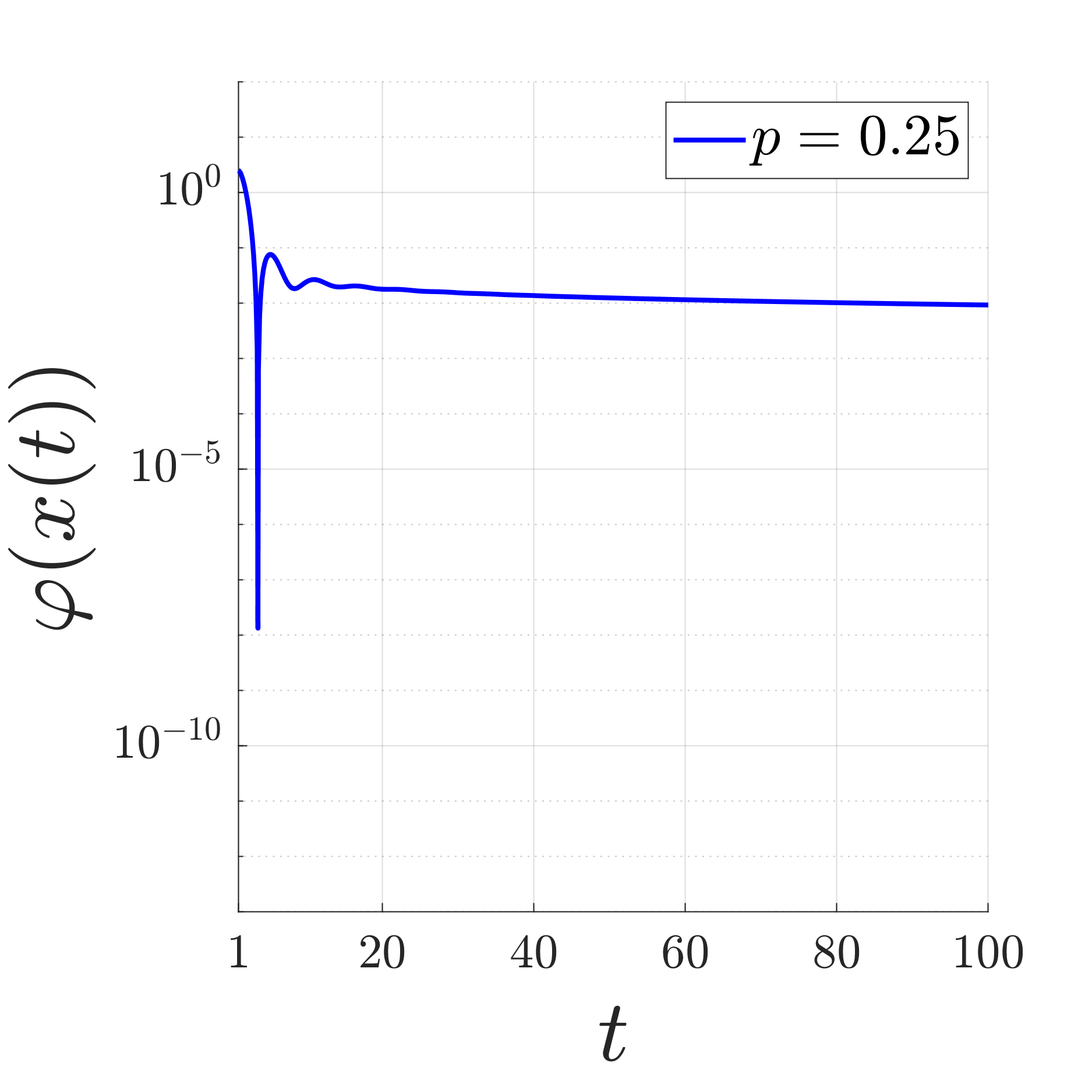}
            \caption{}
            \label{subfig:phi_p=.25}
        \end{subfigure}
        \begin{subfigure}[t]{.24\textwidth}
            \centering
            \includegraphics[width=\linewidth]{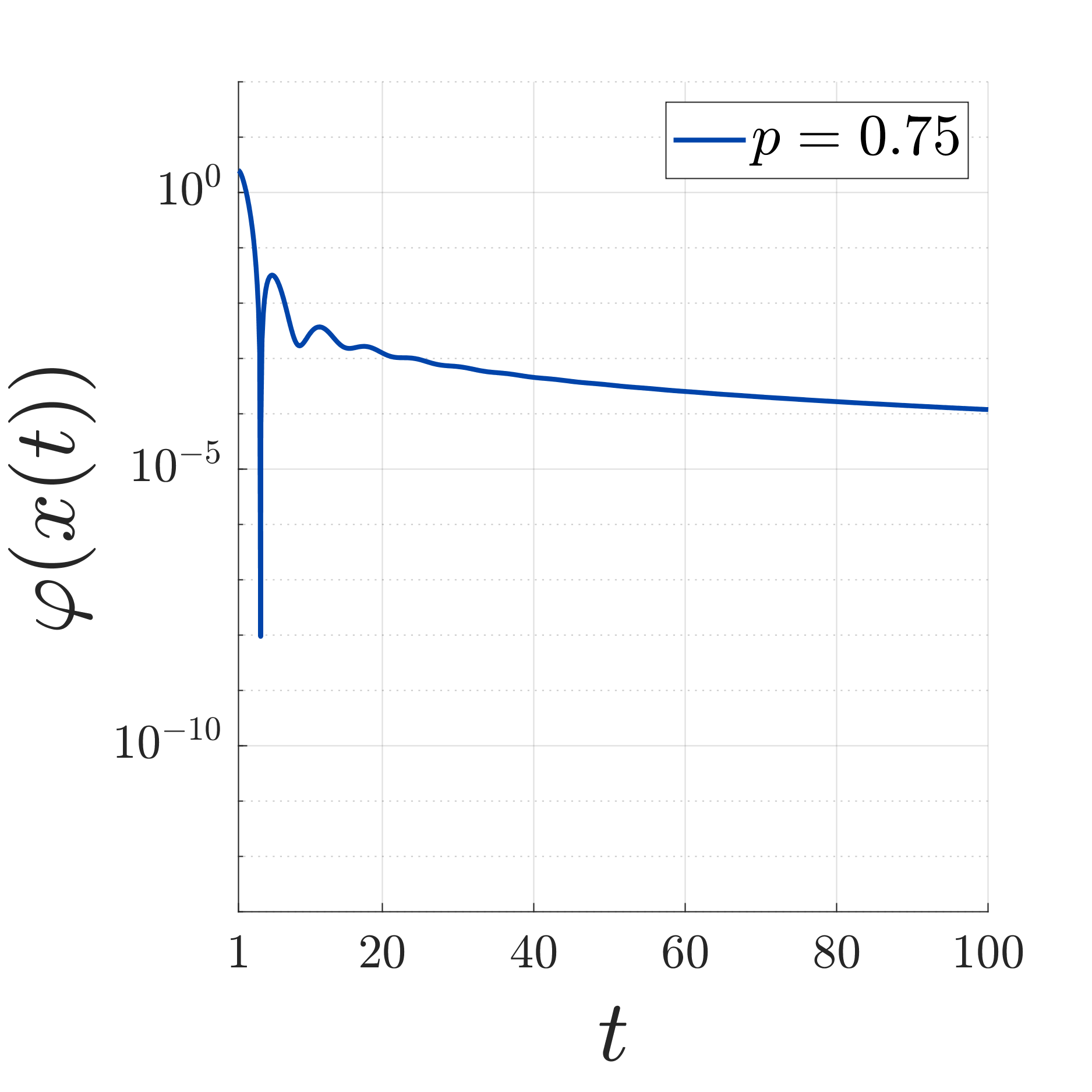}
            \caption{}
            \label{subfig:phi_p=.75}
        \end{subfigure}    
          \begin{subfigure}[t]{.24\textwidth}
            \centering
            \includegraphics[width=\linewidth]{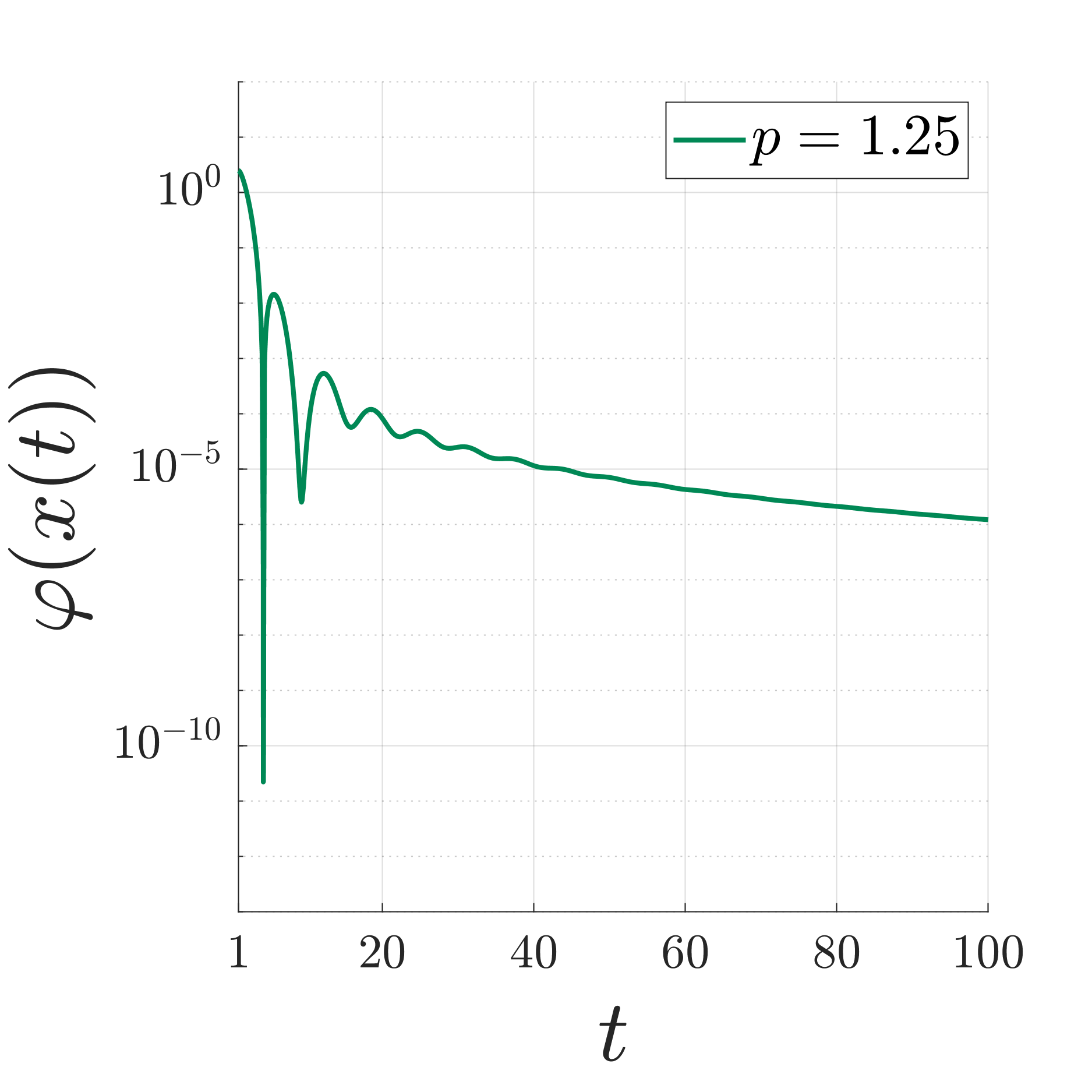}
            \caption{}
            \label{subfig:phi_p=1.25}
          \end{subfigure}
          \begin{subfigure}[t]{.24\textwidth}
            \centering
            \includegraphics[width=\linewidth]{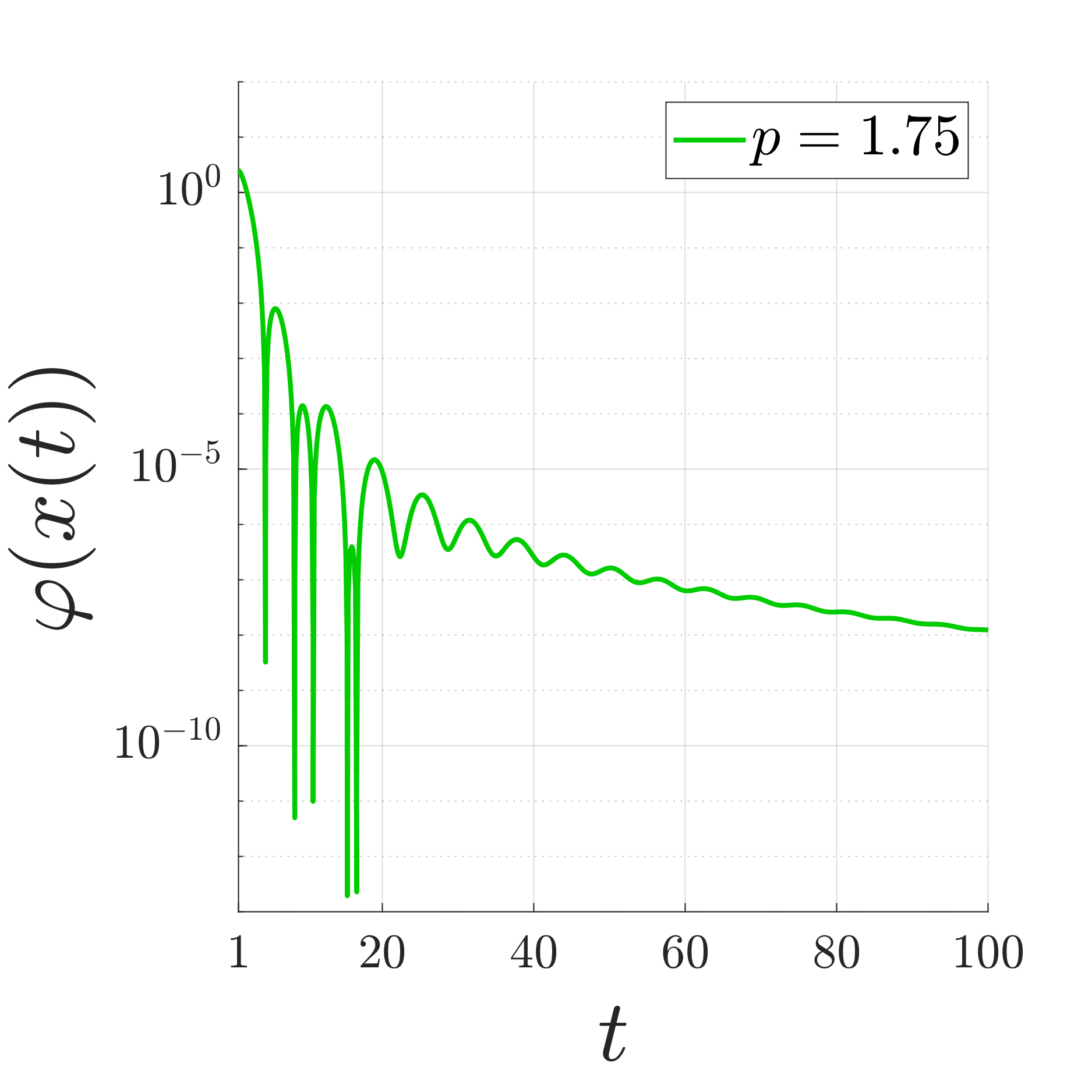}
            \caption{}
            \label{subfig:phi_p=1.75}
          \end{subfigure}
    \medskip
          \begin{subfigure}[t]{.24\textwidth}
            \centering
            \includegraphics[width=\linewidth]{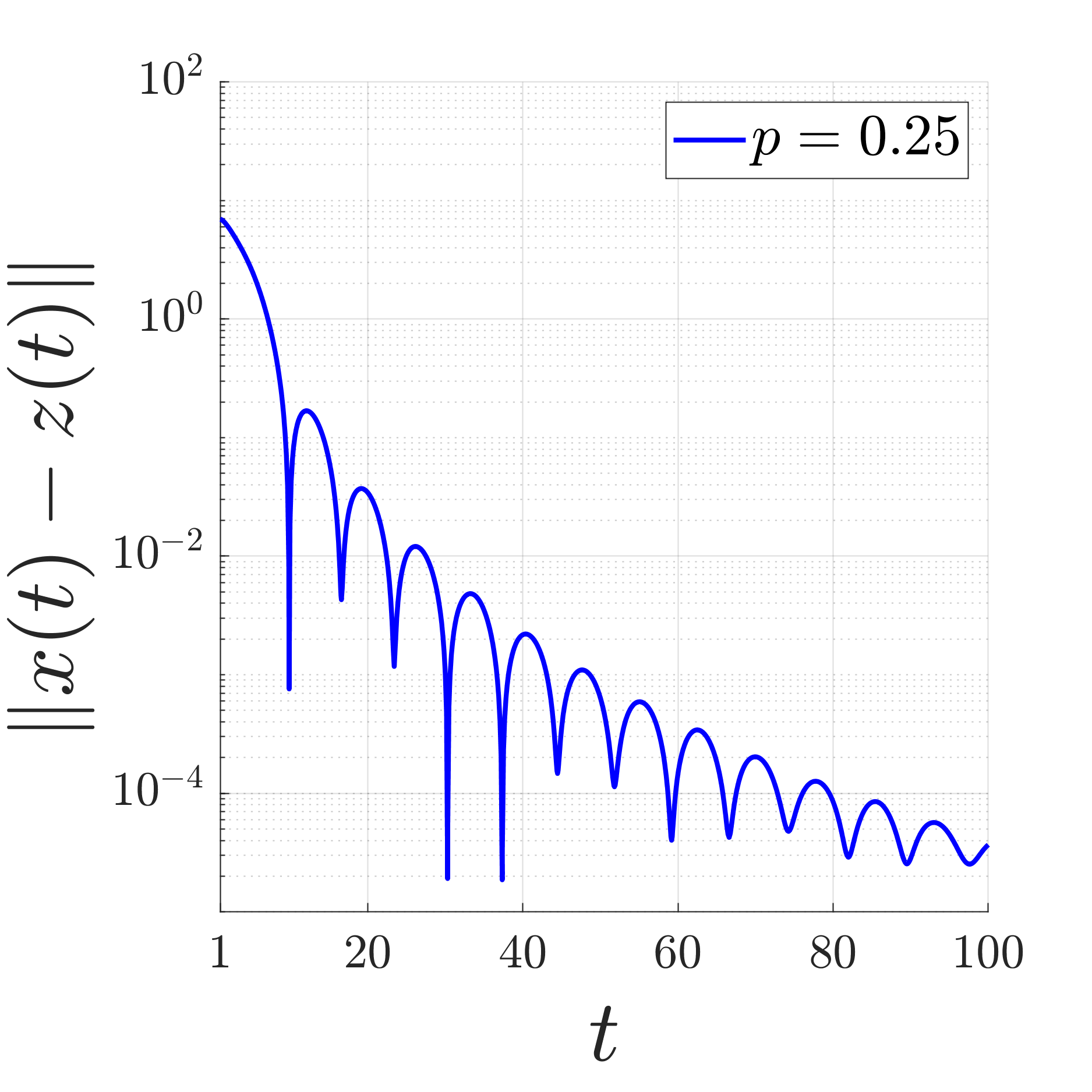}
            \caption{}
            \label{subfig:dist_p=.25}
        \end{subfigure}
        \begin{subfigure}[t]{.24\textwidth}
            \centering
            \includegraphics[width=\linewidth]{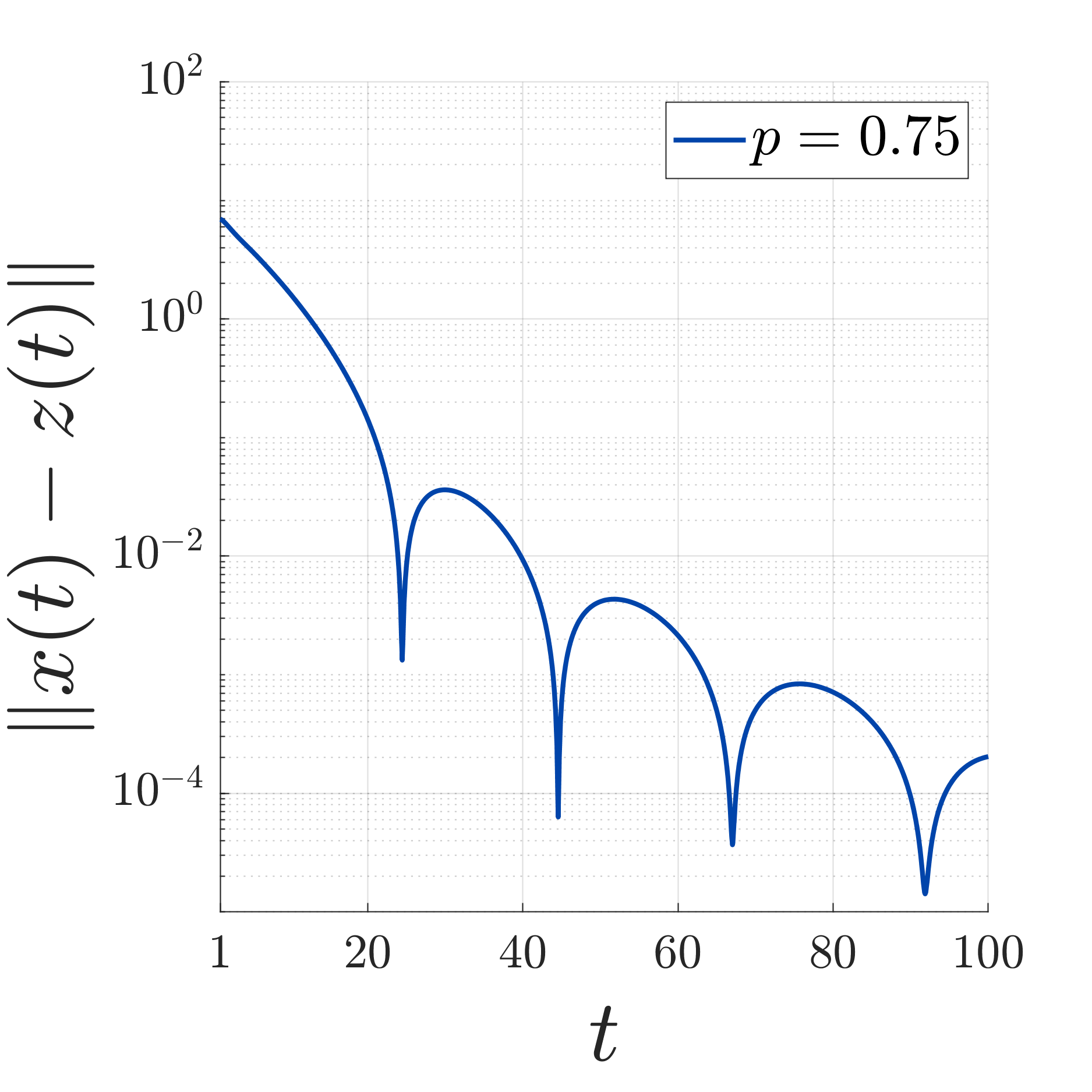}
            \caption{}
            \label{subfig:dist_p=.75}
        \end{subfigure}    
          \begin{subfigure}[t]{.24\textwidth}
            \centering
            \includegraphics[width=\linewidth]{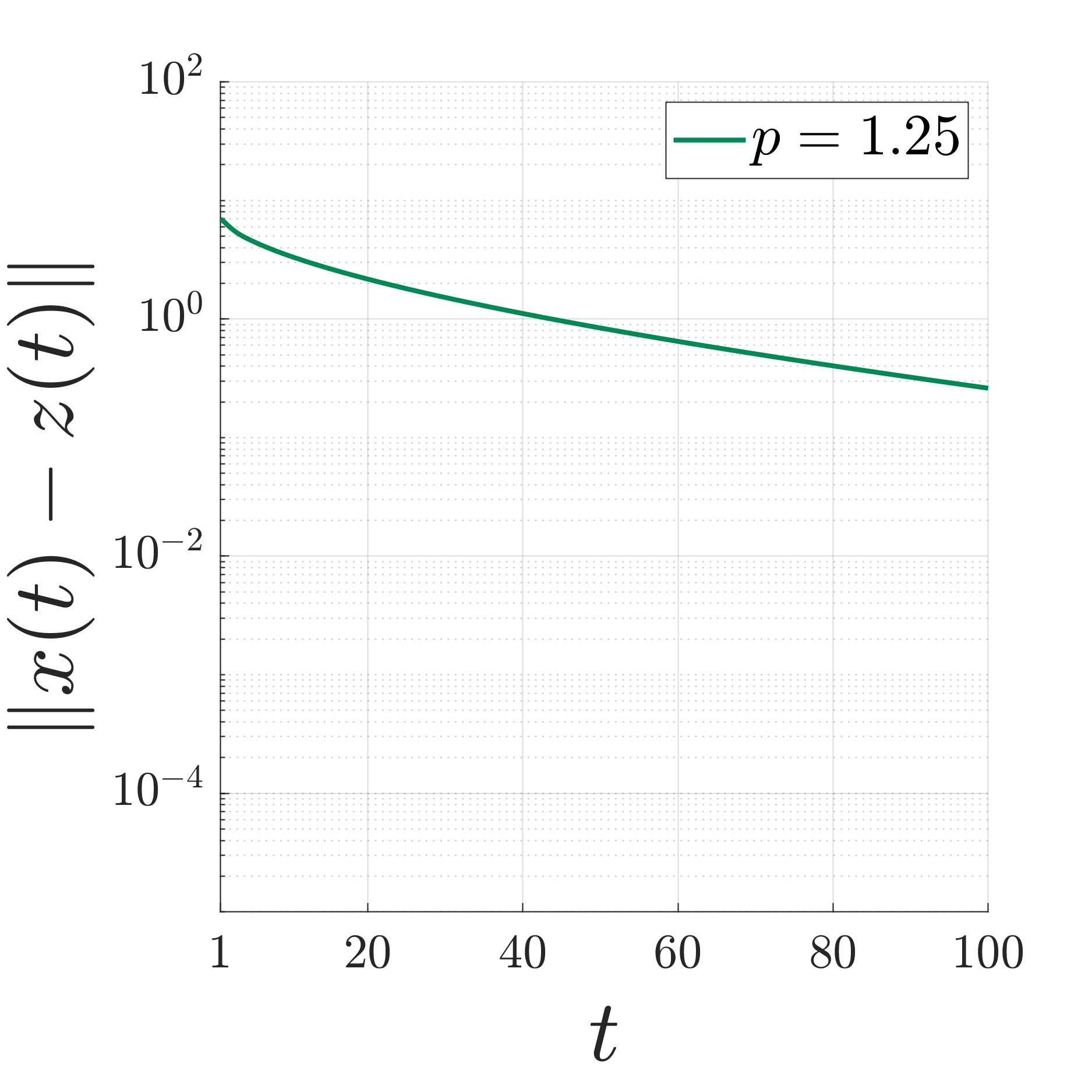}
            \caption{}
            \label{subfig:dist_p=1.25}
          \end{subfigure}
          \begin{subfigure}[t]{.24\textwidth}
            \centering
            \includegraphics[width=\linewidth]{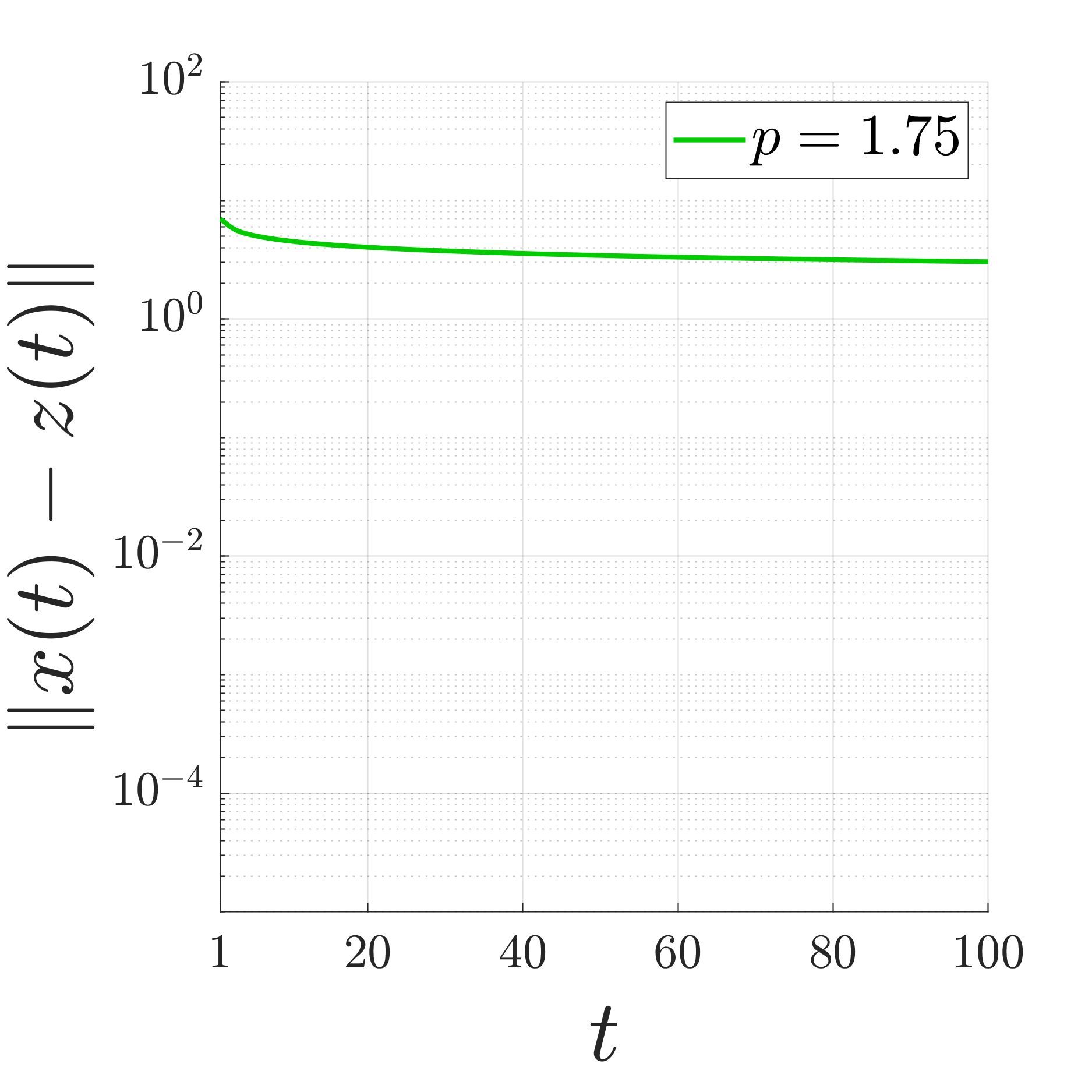}
            \caption{}
            \label{subfig:dist_p=1.75}
          \end{subfigure}
          \caption{The merit function values  $\varphi(x(t))$ and the distance $\lVert x(t) - z(t) \rVert$ of the trajectory to the generalized regularization path for $q = 0.8$ and $p \in \{0.25,\,0.75,\,1.25,\,1.75\}$.}
          \label{fig:phi_dist_p}
    \end{center}
\end{figure}

Figure \ref{fig:phi_dist_p} shows the evolution of the merit function values $\varphi(x(t))$ and of the distance  $\lVert x(t) - z(t)\rVert$ of the trajectory to the generalized regularization path for $q = 0.8$ and $p \in \{0.25,\,0.75,\,1.25,\,1.75\}$. The merit function values exhibit the fastest decay for the largest value of $p=1.75$. This behavior is expected, as higher values of $p$ cause the Tikhonov regularization parameter to decay more rapidly, thus exerting less influence and allowing the function values to converge more quickly. Conversely, the distance $\lVert x(t) - z(t)\rVert$ decays most rapidly for smaller values of $p$, where the regularization parameter vanishes more slowly and effectively guides the trajectory towards the regularization path.

\begin{figure}
    \begin{center}        
        \begin{subfigure}[t]{.24\textwidth}
            \centering
            \includegraphics[width=\linewidth]{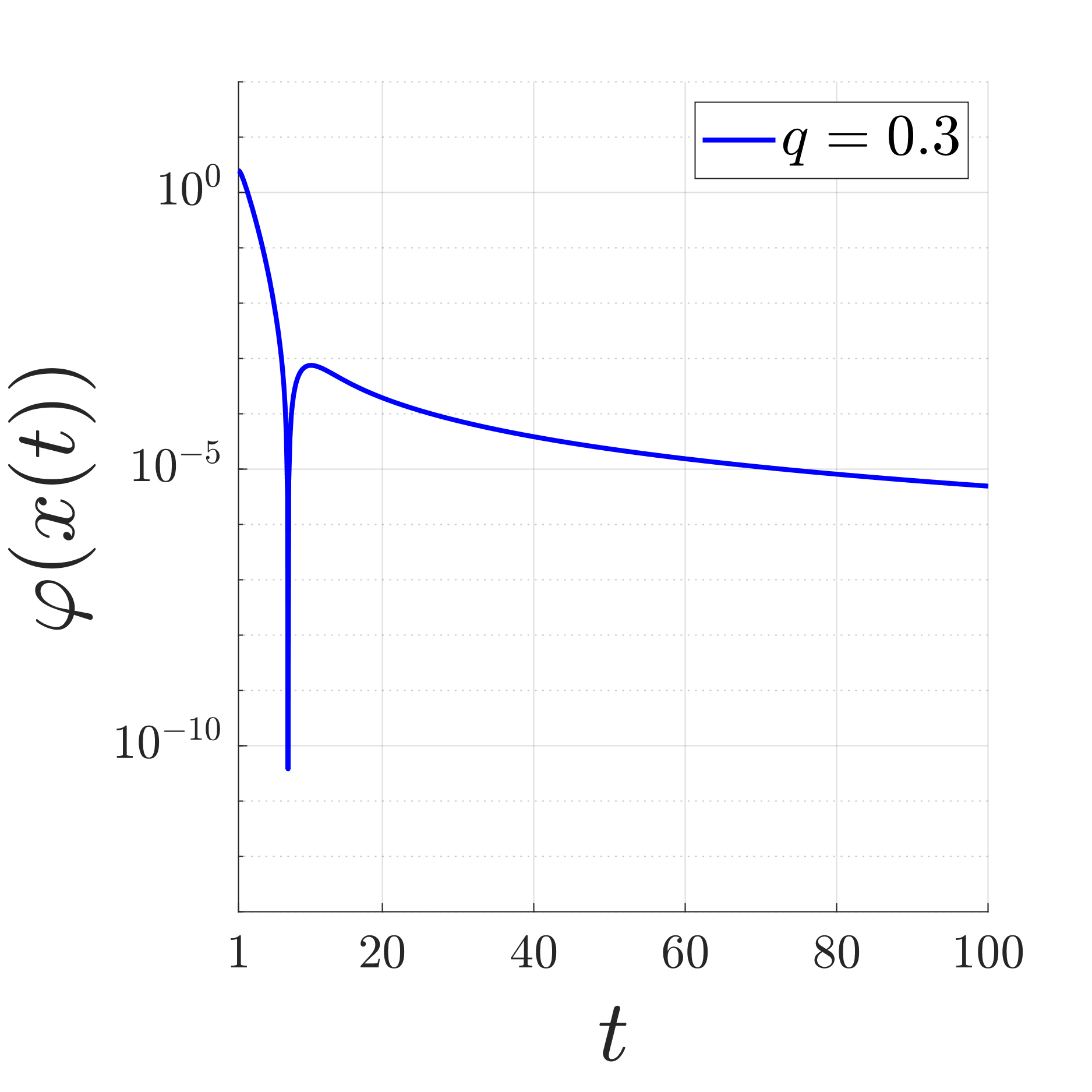}
            \caption{}
            \label{subfig:phi_q=.3}
        \end{subfigure}
        \begin{subfigure}[t]{.24\textwidth}
            \centering
            \includegraphics[width=\linewidth]{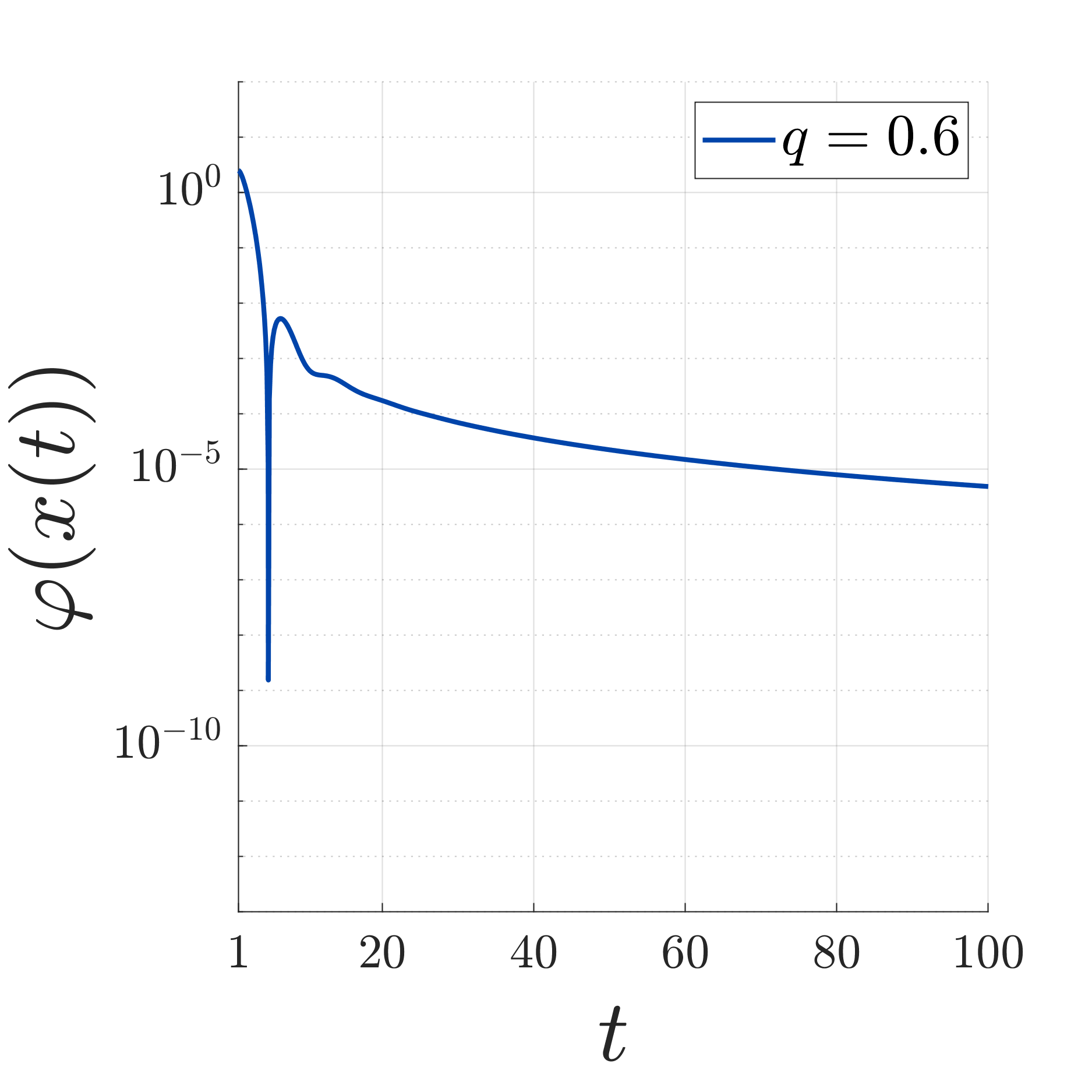}
            \caption{}
            \label{subfig:phi_q=.6}
        \end{subfigure}    
          \begin{subfigure}[t]{.24\textwidth}
            \centering
            \includegraphics[width=\linewidth]{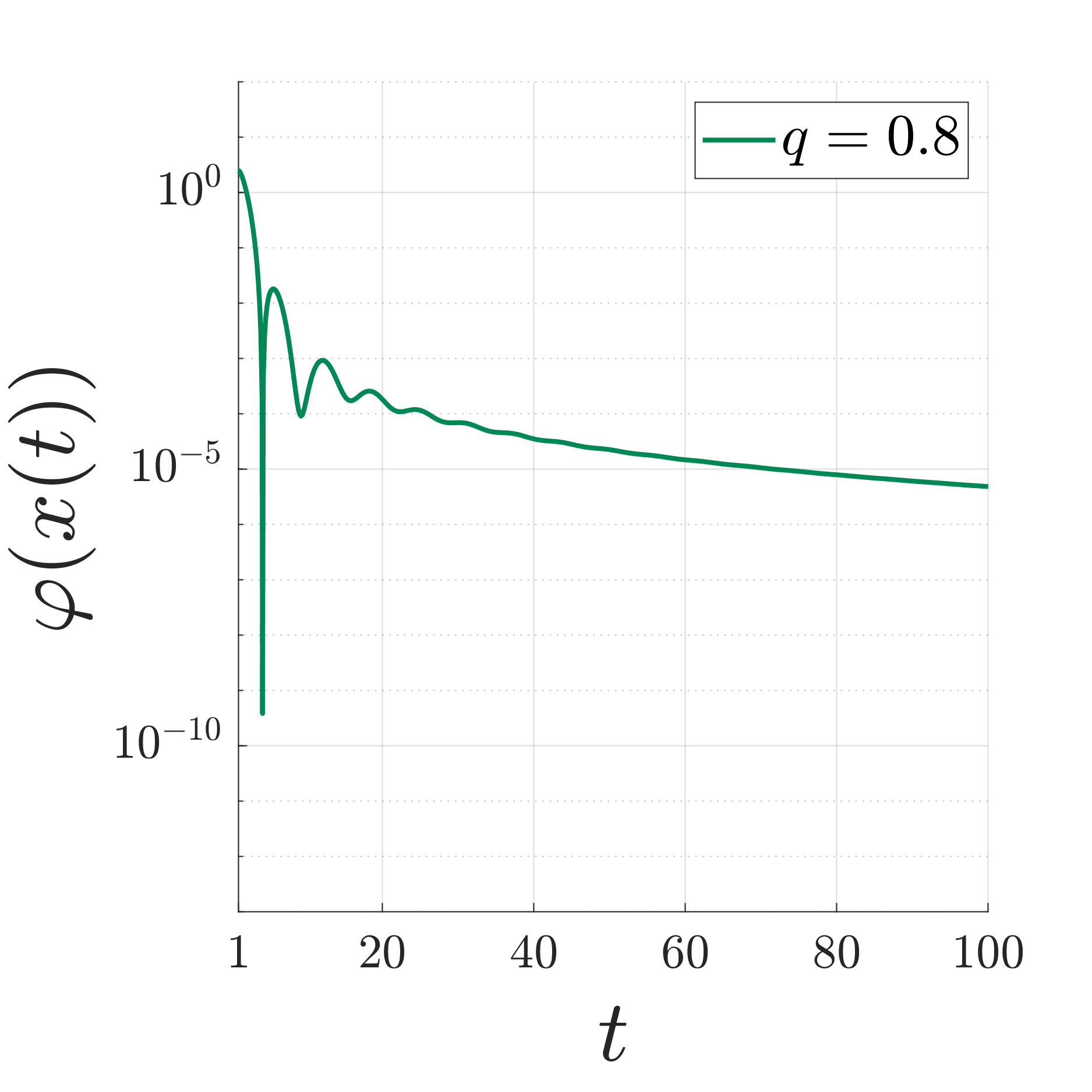}
            \caption{}
            \label{subfig:phi_q=.8}
          \end{subfigure}
          \begin{subfigure}[t]{.24\textwidth}
            \centering
            \includegraphics[width=\linewidth]{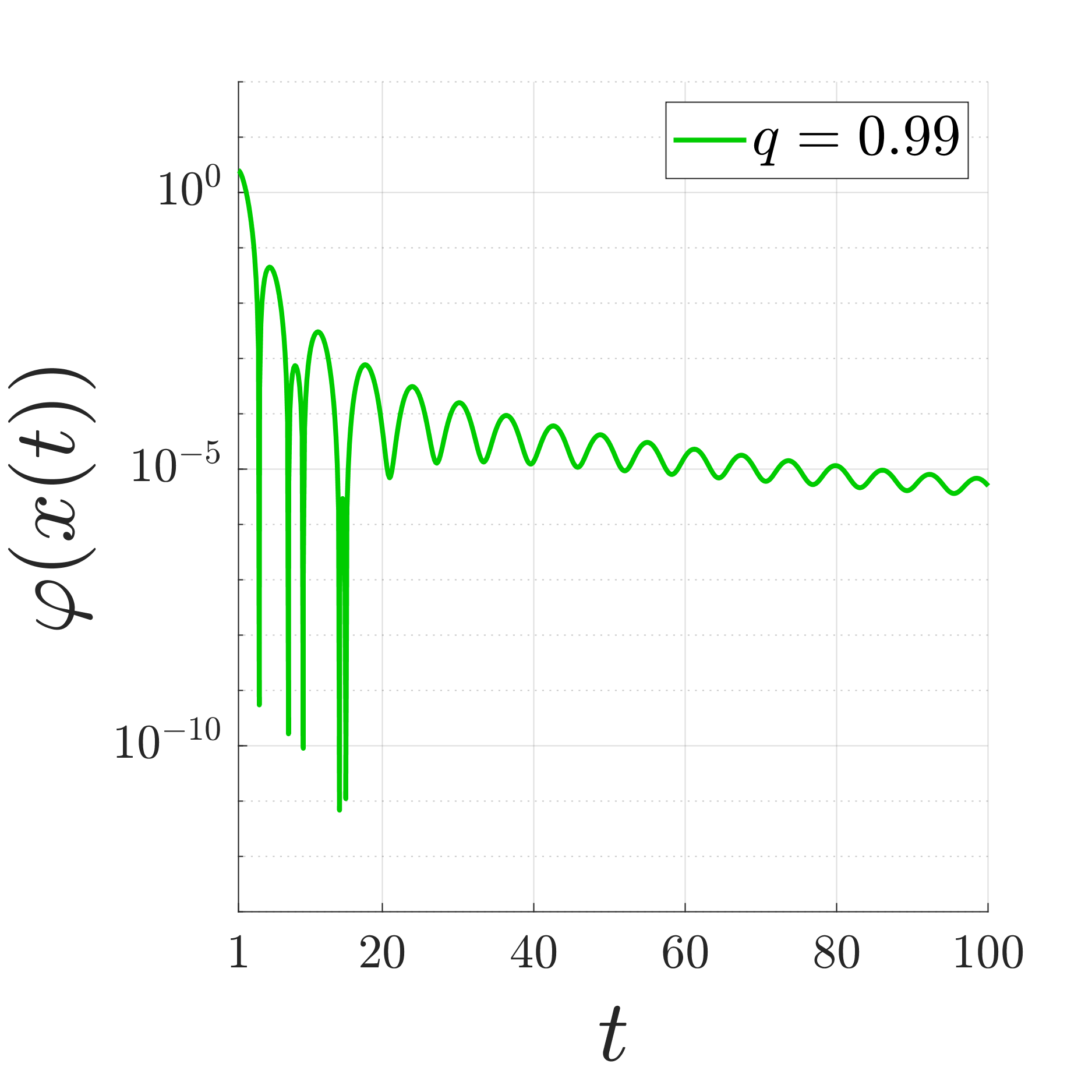}
            \caption{}
            \label{subfig:phi_q=.99}
          \end{subfigure}
    \medskip
          \begin{subfigure}[t]{.24\textwidth}
            \centering
            \includegraphics[width=\linewidth]{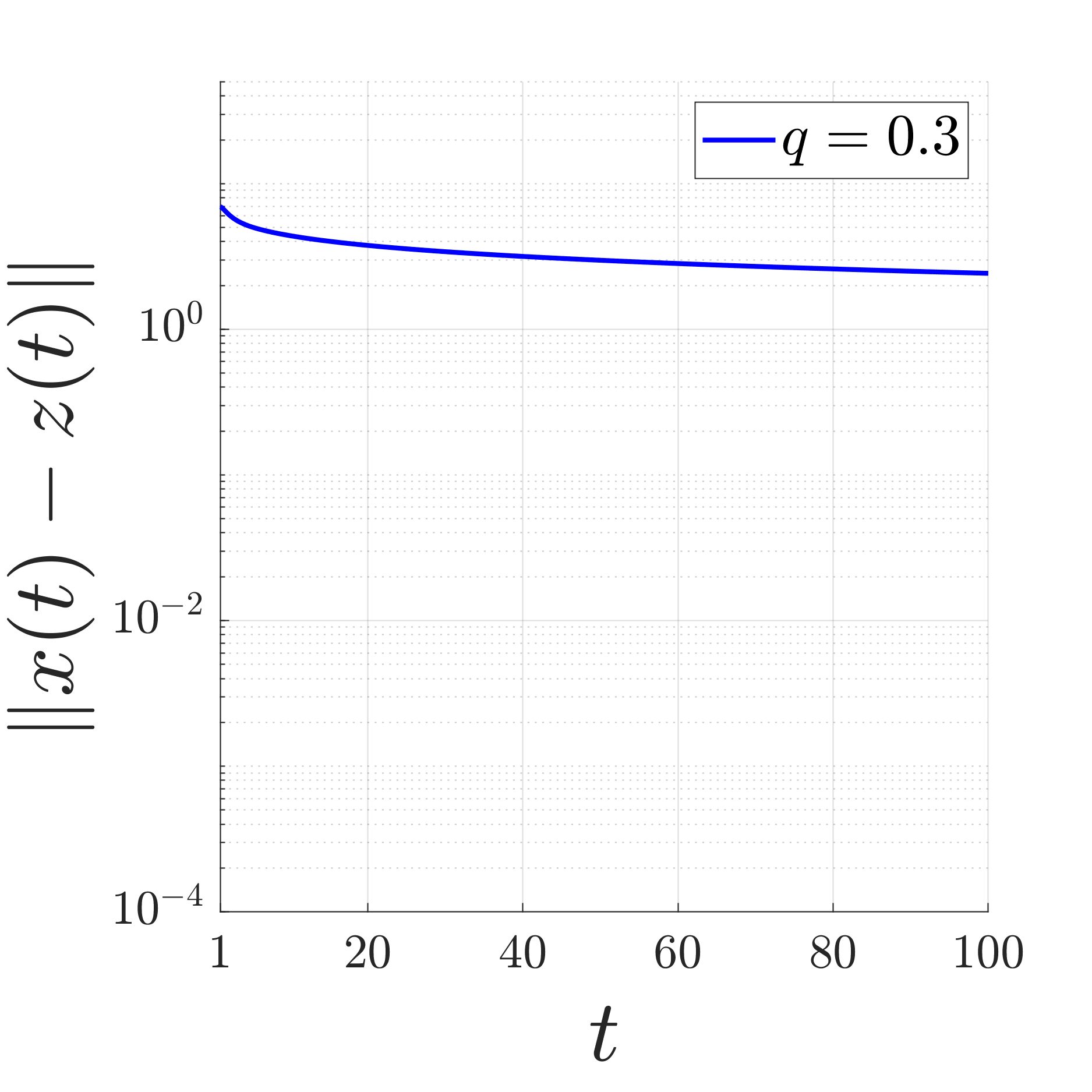}
            \caption{}
            \label{subfig:dist_q=.3}
        \end{subfigure}
        \begin{subfigure}[t]{.24\textwidth}
            \centering
            \includegraphics[width=\linewidth]{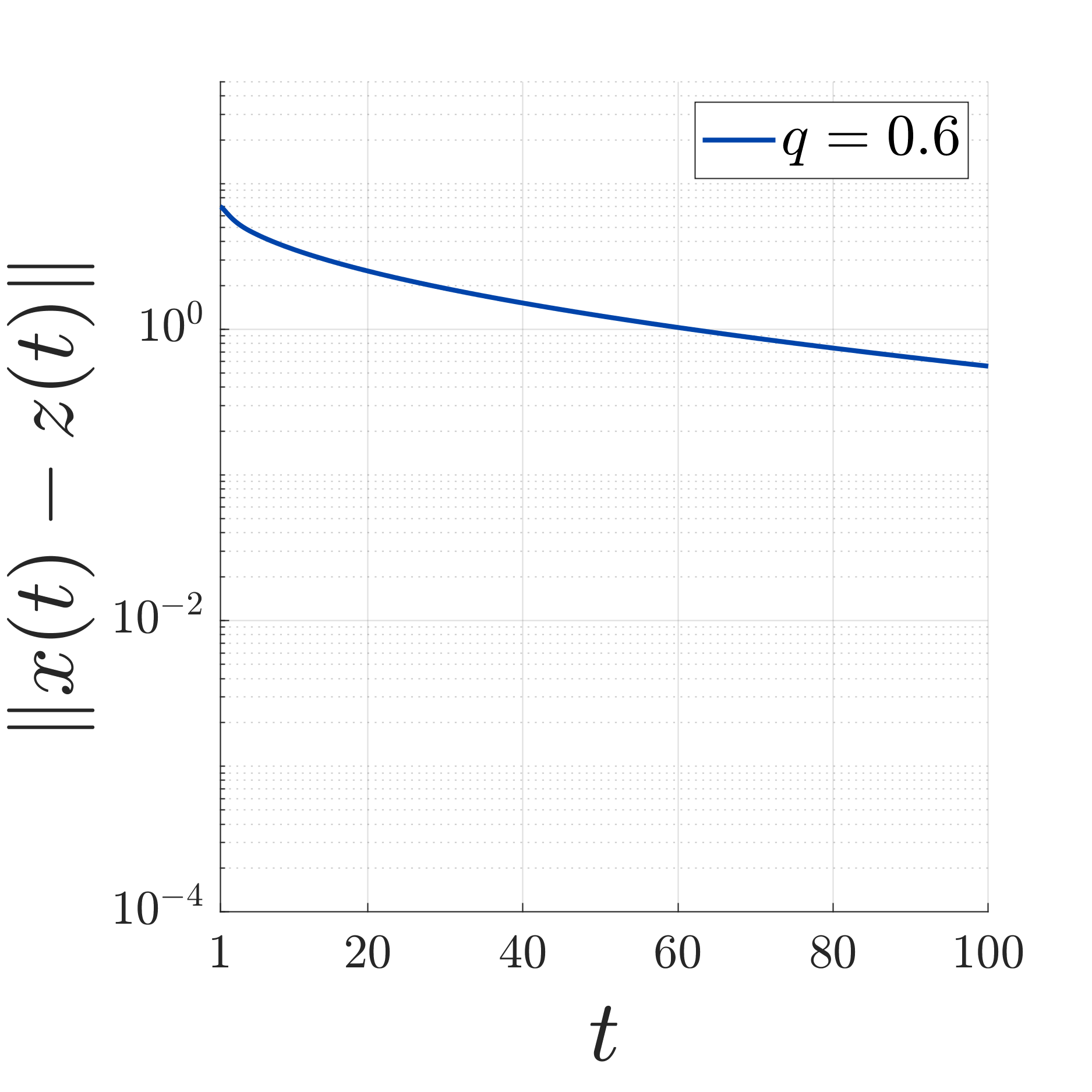}
            \caption{}
            \label{subfig:dist_q=.6}
        \end{subfigure}    
          \begin{subfigure}[t]{.24\textwidth}
            \centering
            \includegraphics[width=\linewidth]{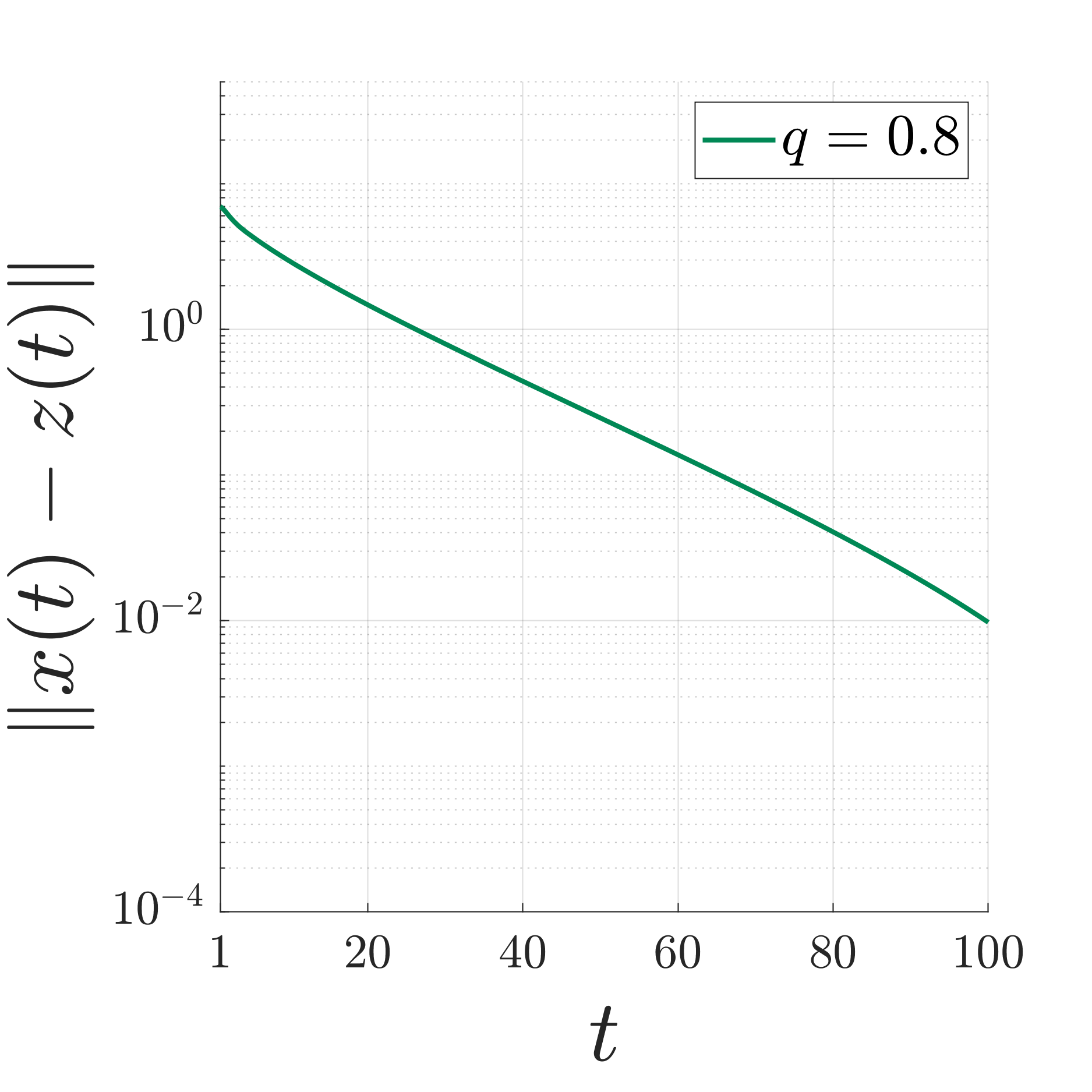}
            \caption{}
            \label{subfig:dist_q=.8}
          \end{subfigure}
          \begin{subfigure}[t]{.24\textwidth}
            \centering
            \includegraphics[width=\linewidth]{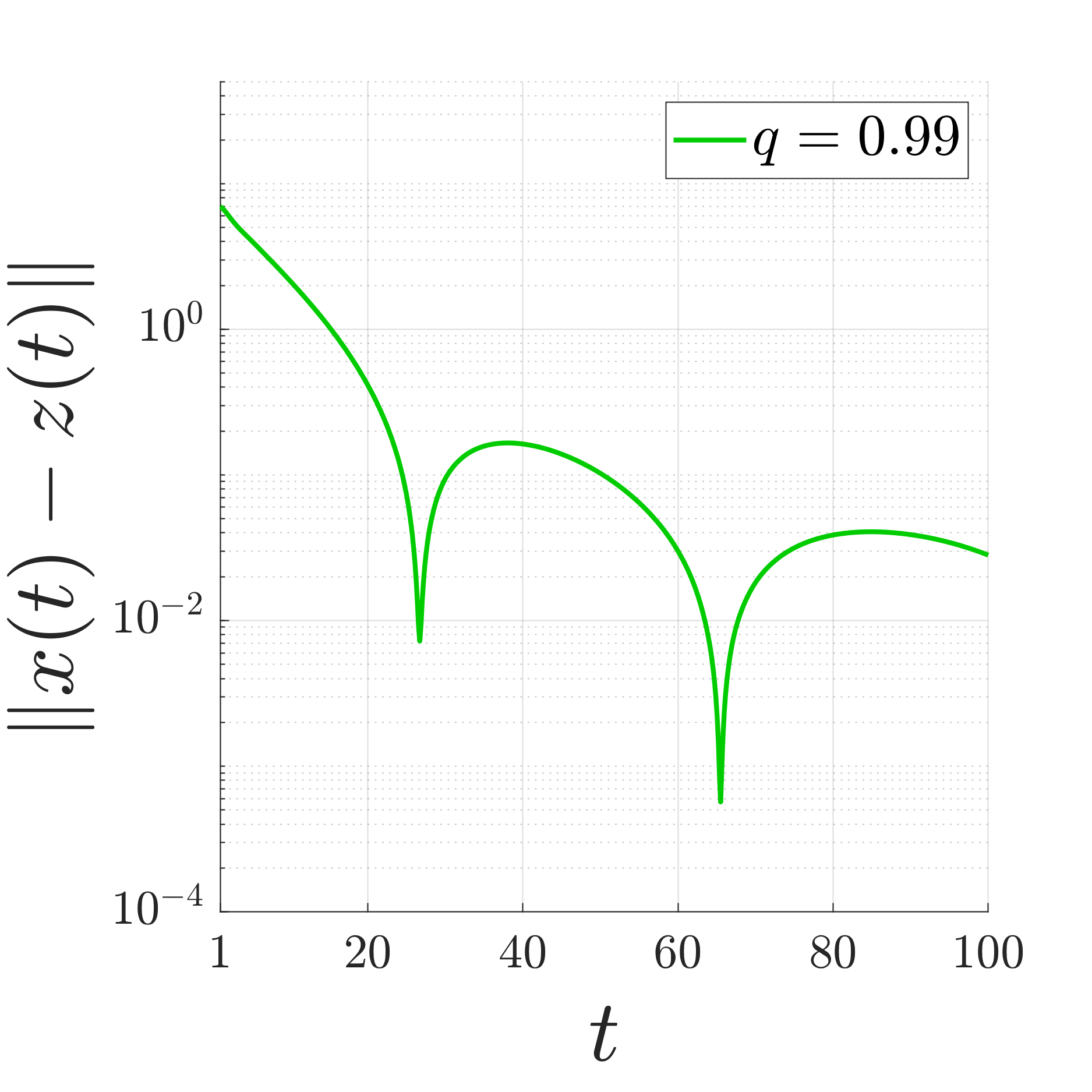}
            \caption{}
            \label{subfig:dist_q=.99}
          \end{subfigure}
          \caption{The merit function values $\varphi(x(t))$ and the distance $\lVert x(t) - z(t) \rVert$ of the trajectory to the generalized regularization path for $p = 1.1$ and $q \in \{ 0.3,\,0.6,\,0.8,\,0.99 \}$.}
          \label{fig:dist_q_p}
    \end{center}
\end{figure}

Figure \ref{fig:dist_q_p} shows the evolution of the merit function values $\varphi(x(t))$ and the distance $\lVert x(t) - z(t) \rVert$ of the trajectory to the generalized regularization path for $p = 1.1$ and $q \in \{ 0.3,\,0.6,\,0.8,\,0.99 \}$. The decay of the merit function values $\varphi(x(t))$ is generally insensitive to the choice of $q$; for all considered values of $q$, the convergence rate remains essentially the same. However, for larger values of $q$, the merit function exhibits more pronounced oscillations. This behavior is expected, as a larger value of $q$ implies a faster decay of the friction term $\frac{\alpha}{t^q}$, thereby reducing damping. In contrast, the decay of the distance  $\lVert x(t) - z(t) \rVert$ is strongly influenced by $q$, particularly for $q = 0.99$, where convergence is significantly faster. For the smallest value $q = 0.3$, the distance decreases only slowly, at a sublinear rate. These observations align with expectations: higher values of $q$ correspond to weaker friction, which allows the trajectory to approach the regularization path more rapidly in the early phase.


\section{Conclusion}\label{sec:conclusion}

In this paper, we propose a novel second-order dynamical system, \eqref{eq:MTRIGS}, tailored for multiobjective optimization problems. This system incorporates asymptotically vanishing damping and vanishing Tikhonov regularization. Leveraging existence theorems for differential inclusions, we establish the existence of solutions to this system in the finite dimensional setting. To analyze the asymptotic behavior of the trajectory solutions, we introduce a new regularization path for multiobjective optimization problems, derived from the Tikhonov regularization of an adaptive scalarization. Using this framework, we demonstrate the strong convergence of the trajectory solutions $x(\cdot)$ of \eqref{eq:MTRIGS} to the weak Pareto optimal point with minimal norm in a particular lower level set of the objective function. Furthermore, we recover fast convergence rates quantified in terms of a merit function. We investigate the qualitative behavior of the solution to \eqref{eq:MTRIGS} through multiple numerical experiments. These findings form the basis for developing inertial proximal point methods with vanishing Tikhonov regularization for multiobjective optimzsation problems, which yield fast convergence of function values and strong convergence of iterates. Future research directions include designing second-order gradient dynamics for multiobjective optimization problems with Hessian-driven damping, as well as addressing multiobjective problems with linear constraints using primal-dual dynamical systems.


\appendix
\section{Auxiliary lemmas}
\label{sec:ap_aux_lem}
In the first part of the appendix we introduce some auxiliary lemmas that we use in the asymptotic analysis of the trajectory solutions of \eqref{eq:MTRIGS}.

\begin{mylemma}
    For $i=1, \dots, m$, let $h_i:[t_0, +\infty) \to \R$ be absolutely continuous functions on every interval $[t_0, T]$ for $T \geq t_0$. Define $h:[t_0, +\infty) \to \R,\quad t \mapsto h(t) \coloneqq \min_{i=1,\dots, m}h_i(t)$. Then, the following statements are true:
    \begin{enumerate}
        \item[i\emph{)}] The function $h$ is absolutely continuous on every interval $[t_0, T]$ for $T \geq t_0$, and therefore differentiable at almost all $t \geq t_0$;
        \item[ii\emph{)}] For almost all $t \ge t_0$ there exists $i \in \{1,\dots, m\}$ such $h(t) = h_i(t)$ and $\frac{d}{dt}h(t) = \frac{d}{dt}h_i(t)$.
    \end{enumerate}
    \label{lem:diff_a_e}
\end{mylemma}
\begin{proof}~\\
\vspace{-5mm}
\begin{enumerate}
    \item[\emph{i})] The minimum of a family of finitely many absolutely continuous functions is absolutely continuous.

    \item[\emph{ii})]  Let $t \ge t_0$ be such that $h(\cdot)$ and $h_i(\cdot)$ are differentiable in $t$ for all $i =1,\dots,m$. Take an arbitrary sequence $\lbrace \tau_k \rbrace_{k \ge 0}$ with $\lim_{k \to + \infty} \tau_k = 0$. Then, there exists $i \in \{ 1, \dots, m\}$ and a subsequence $\lbrace k_l \rbrace_{l \ge 0} \subset \N$ with $h(t + \tau_{k_l}) = h_i(t + \tau_{k_l})$ for all $l \ge 0$. From the continuity of $h(\cdot)$ and $h_i(\cdot)$, it holds $h(t) = h_i(t)$. By the definition of the derivative, we get
    \begin{align*}
        \frac{d}{dt}h(t) = \lim_{l \to + \infty} \frac{h(t + \tau_{k_l}) - h(t)}{\tau_{k_l}} = \lim_{l \to + \infty} \frac{h_i(t + \tau_{k_l}) - h_i(t)}{\tau_{k_l}} = \frac{d}{dt}h_i(t).
    \end{align*}
\end{enumerate}
\end{proof}

\begin{mylemma}
\label{lem:int_poly_exp_bound}
    Let $\alpha, \beta, a, b > 0$ be given constants, and $t_0 > 0$. Then,
    \begin{align*}
        \int_{t_0}^t \alpha s^{-a}\exp(\beta s^b) ds = \mathcal{O}\left( t^{1 - (a+b)} \exp(\beta t^b) \right) \, \text{ as } t \to + \infty. 
    \end{align*}
\end{mylemma}

\begin{proof}
For $t \geq t_0$, we use integration by parts to get
    \begin{align}
        \int_{t_0}^t \alpha s^{-a} \exp \left( \beta s^{b} \right) ds = & \frac{\alpha}{\beta b} \int_{t_0}^t s^{1- (a + b)} \frac{d}{ds} \exp \left( \beta s^{b} \right) ds \nonumber\\
        = & \frac{\alpha}{\beta b} \left[ s^{1 - (a + b)} \exp \left( \beta s^{b} \right)\right]_{t_0}^t - \frac{1 - (a+b)}{\beta b} \int_{t_0}^t \alpha s^{-(a + b)} \exp \left( \beta s^b \right) ds.\label{eq:int_parts_b_1}
    \end{align}
    Since $b > 0$, there exists $t_1 \ge t_0$ such that for all $t \ge t_1$
    \begin{align}
    \label{eq:bound_s^b}
        \left\lvert \frac{1 - (a + b)}{\beta b} \right\rvert t^{-b} \le \frac{1}{2}.
    \end{align}
Define $C_1 \coloneqq \left\lvert \frac{1 - (a+b)}{\beta b} \right\rvert \int_{t_0}^{t_1} \alpha s^{-(a+b)} \exp \left( \beta s^b \right) ds$. Then, \eqref{eq:int_parts_b_1} and \eqref{eq:bound_s^b} yield for all $t \geq t_0$
    \begin{align*}
       \int_{t_0}^t \alpha s^{-a} \exp \left( \beta s^{b} \right) ds \le & \frac{\alpha} {\beta b} \left[ s^{1 - (a+b)} \exp \left( \beta s^b \right)\right]_{t_0}^{t} + C_1 + \left\lvert \frac{1 - (a + b)}{\beta b} \right\rvert \int_{t_1}^{t} \alpha  s^{- (a + b)} \exp \left( \beta s^b \right) ds\\
        \le & \frac{\alpha} {\beta b} \left[ s^{1 - (a+b)} \exp \left( \beta s^b \right)\right]_{t_0}^{t} + C_1 + \frac{1}{2} \int_{t_1}^{t} \alpha s^{- a} \exp \left( \beta s^b \right) ds \\
        \le & \frac{\alpha} {\beta b} \left[ s^{1 - (a+b)} \exp \left( \beta s^b \right)\right]_{t_0}^{t} + C_1 + \frac{1}{2} \int_{t_0}^{t} \alpha s^{- a} \exp \left( \beta s^b \right) ds,
    \end{align*}
hence
    \begin{align*}
        \int_{t_0}^t \alpha s^{-a} \exp \left( \beta s^{b} \right) ds \le \frac{2\alpha} {\beta b} \left[ s^{1 - (a+b)} \exp \left( \beta s^b \right)\right]_{t_0}^{t} + 2C_1.
    \end{align*}
    Defining $C_2 \coloneqq -\frac{2\alpha} {\beta b}(t_0)^{1 - (a+b)} \exp \left( \beta (t_0)^b \right) + 2 C_1$, we obtain for all $t \geq t_0$
    \begin{align*}
        \int_{t_0}^t \alpha s^{-a} \exp \left( \beta s^{b} \right) ds \le \frac{2\alpha}{\beta b}  t^{1 - (a+b)} \exp \left( \beta t^b \right) + C_2,
    \end{align*}
    and the asymptotic bound holds.
\end{proof}

To prove weak convergence of the trajectory solutions, we use the following continuous version of Opial's Lemma (see \cite[Lemma 5.7]{Attouch2018}).

\begin{mylemma}
Let $S \subseteq \H$ be a nonempty set and let $x:[t_0, +\infty) \to \H$ be a function satisfying the following conditions:
\begin{enumerate}
\item[(i)] For every $z\in S$, $\lim_{t\to + \infty} \lVert x(t) - z \rVert$ exists;
\item[(ii)] Every weak sequential cluster point of $x$ belongs to $S$.
\end{enumerate}
Then, $x(t)$ converges weakly to an element $x^{\infty} \in S$ as $t \to +\infty$.
\label{lem:opial}
\end{mylemma}

The following lemma is a modification of \cite[Lemma 16]{Laszlo2023}.

\begin{mylemma}
    \label{lem:converngence_h}
    Let $t_0 > 0$, $\alpha > 0$, $q \in (0,1)$, and $k:[t_0, +\infty) \to \R$ a nonnegative function such that
    \begin{align}
    \label{eq:integral_estimate_k(t)}
        \left( t \mapsto t^q k(t) \right) \in L^1\left( [t_0, + \infty) \right).
    \end{align}
    Let $h:[t_0, +\infty) \to \R$ be a continuously differentiable function that is bounded from below and possesses an absolutely continuous derivative $h'(\cdot)$. Further, assume $h(\cdot)$ satisfies
    \begin{align}
    \label{eq:h''(t)_h'(t)_k(t)}
        h''(t) + \frac{\alpha}{t^q} h'(t) \le k(t) \quad \text{for almost all} \ t \geq t_0.
    \end{align}
 Then, $\left( t \mapsto \left[ h'(t) \right]_+ \right) \in L^1\left( [t_0, + \infty) \right)$, where $\left[ h'(t) \right]_+$ denotes the positive part of $h'(t)$, and further $\lim_{t \to + \infty} h(t)$ exists.
\end{mylemma}

\begin{proof}
    Define the function
    \begin{align*}
        \M:[t_0, +\infty) \to \R,\, t\mapsto  \M(t) \coloneqq \exp\left( \int_{t_0}^t \frac{\alpha}{s^q} ds \right) = C_{\M} \exp \left( \frac{\alpha}{1 - q} t^{1 - q}\right),
    \end{align*}
    with $C_{\M} \coloneqq \exp \left( - \frac{\alpha}{1 - q} t_0^{1 - q}\right)$, and $b \coloneqq \frac{\alpha}{1 - q} > 0$. For $t \ge t_0$, using integration by parts, we have
    \begin{align}
    \label{eq:converngence_h_1}
        C_{\M} \int_{t}^{+\infty} \frac{ds}{\M(s)}  = & \int_{t}^{+\infty} \exp\left( - b s^{1 - q} \right) ds = - \frac{1}{\alpha} \int_{t}^{+\infty}  s^q \frac{d}{ds} \exp\left( - b s^{1 - q} \right) ds \nonumber\\
        = & - \frac{1}{\alpha} \left( \left[ s^q \exp\left( - b s^{1 - q} \right) \right]_{t}^{+\infty} - \int_{t}^{+\infty} qs^{q-1} \exp\left( - b s^{1 - q} \right) ds \right) \\
        = & \frac{t^q}{\alpha} \exp\left( -b t^{1-q} \right) + \frac{q}{\alpha} \int_{t}^{+\infty} s^{q-1} \exp(-b s^{1-q}) ds. \nonumber
    \end{align}
    As $q - 1 < 0$, there exists $t_1 \ge t_0$ such that for all $t \ge t_1$ the inequality $\frac{q}{\alpha} t^{q-1} \le \frac{1}{2}$ holds and hence
    \begin{align}
    \label{eq:converngence_h_3}
        \frac{q}{\alpha} \int_{t}^{+\infty} s^{q-1} \exp(-b s^{1-q}) ds \le \frac{1}{2} \int_{t}^{+\infty}\exp(-b s^{1-q}) ds.
    \end{align}
    Combining \eqref{eq:converngence_h_1} and \eqref{eq:converngence_h_3}, we conclude that for all $t \ge t_1$
    \begin{align}
    \label{eq:converngence_h_4}
        C_{\M} \int_{t}^{+\infty} \frac{ds}{\M(s)} = \int_{t}^{+\infty} \exp\left( - b s^{1 - q} \right) ds \le \frac{2 t^q}{\alpha} \exp\left( -b t^{1-q} \right).
    \end{align}
    Using the definition of $\M(\cdot)$, equality \eqref{eq:converngence_h_4} yields for all $t \ge t_1$
    \begin{align}
    \label{eq:converngence_h_45}
        \left( \int_{t}^{+\infty} \frac{ds}{\M(s)} \right) \M(t) = \left( \int_{t}^{+\infty} \exp\left( - b s^{1 - q} \right) \right) \exp\left( b t^{1-q} \right) \le \frac{2t^q}{\alpha}.
    \end{align}
    We multiply \eqref{eq:converngence_h_45} by $k(\cdot)$, integrate from $t_0$ to $+\infty$, and apply relation \eqref{eq:integral_estimate_k(t)} to follow
    \begin{align}
    \label{eq:converngence_h_4.5}
        \int_{t_0}^{+\infty} \left( \int_{t}^{+\infty} \frac{ds}{\M(s)} \right) \M(t) k(t) dt < +\infty.
    \end{align}
    By the definition of $\M(\cdot)$, we have $\frac{d}{dt} \M(t) = \M(t) \frac{\alpha}{t^q}$ and then, by \eqref{eq:h''(t)_h'(t)_k(t)},
    \begin{align}
    \label{eq:converngence_h_5}
         \frac{d}{dt} \left( \M(t) h'(t) \right) = \M(t) h''(t) + \M(t) \frac{\alpha}{t^q} h'(t) \le \M(t) k(t) \quad \text{for almost all} \ t \geq t_0.
    \end{align}
We integrate \eqref{eq:converngence_h_5} from $t_0$ to $t \ge t_0$ and observe
    \begin{align*}
        \M(t) h'(t) - \M(t_0) h'(t_0) \le \int_{t_0}^t \M(s) k(s) ds.
    \end{align*}
    The function $k(\cdot)$ takes nonnegative values only and we derive for all $t \ge t_0$
    \begin{align*}
        \left[ h'(t)\right]_+ \le \frac{\lvert\M(t_0) h'(t) \rvert}{\M(t)} + \frac{1}{\M(t)} \int_{t_0}^t \M(s) k(s) ds.
    \end{align*}
    We integrate this inequality from $t_0$ to $+\infty$ and write 
    \begin{align}
    \label{eq:converngence_h_6}
        \int_{t_0}^{+\infty} \left[ h'(t)\right]_+ dt \le \int_{t_0}^t \frac{\lvert\M(t_0) h'(t) \rvert}{\M(t)} dt + \int_{t_0}^{+\infty} \frac{1}{\M(t)} \left( \int_{t_0}^t \M(s) k(s) ds \right) dt.
    \end{align}
    Since $\M(\cdot)$ grows at an exponential rate, we have $\int_{t_0}^{+\infty} \frac{\lvert\M(t_0) h'(t) \rvert}{\M(t)} dt < + \infty$. We apply Fubini's Theorem to the second integral in \eqref{eq:converngence_h_6} and combine it with \eqref{eq:converngence_h_4.5} to conclude
    \begin{align}
    \label{eq:converngence_h_7}
         \int_{t_0}^{+\infty} \frac{1}{\M(t)} \left( \int_{t_0}^t \M(s) k(s) ds \right) dt = \int_{t_0}^{+\infty} \left( \int_{t}^{+\infty} \frac{ds}{\M(s)} \right) \M(t) k(t) dt < +\infty. 
    \end{align}
    Equation \eqref{eq:converngence_h_6} and \eqref{eq:converngence_h_7} imply
    \begin{align*}
        \int_{t_0}^{+\infty} \left[ h'(t)\right]_+ dt < + \infty,
    \end{align*}
    and by the lower boundedness of $h(\cdot)$ we follow that $\lim_{t \to + \infty} h(t)$ exists. 
\end{proof}


\section{The proof of the existence of trajectory solutions of (MTRIGS)}
\label{sec:ap_exist_sol_MTRIGS}

The proof for the existence of solutions of \eqref{eq:MTRIGS} is closely related to the proof given in \cite{Sonntag2024b} (see also \cite{Sonntag2024a}) for the existence of solutions of the system \eqref{eq:MAVD_intro}. 

\subsection{Existence of trajectory solutions of a related differential inclusion (DI)}
\label{subsec:ap_exist_sol_DI}

Consider the set-valued map
\begin{align}
    \hspace{-10mm}G: [t_0, + \infty) \times \H \times \H \rightrightarrows \H\times \H, \quad (t,u,v) \mapsto \lbrace v \rbrace \times \left( - \frac{\alpha}{t^q} v - \argmin_{g \in C(u) + \frac{\beta}{t^p} u} \langle g, -  v \rangle \right),
\label{eq:set_valued_G}
\end{align}
with $C(u) \coloneqq \conv\left(\left\lbrace \nabla f_i(u)\,:\, i =1,\dots, m \right\rbrace\right)$, and the differential inclusion

\leqnomode
\begin{align}
    \label{eq:DI}
    \tag{DI}
    \left\lvert 
    \begin{array}{l}
        (\dot{u}(t), \dot{v}(t)) \in G(t, u(t), v(t)),  \\
        \\
        (u(t_0), v(t_0)) = (u_0, v_0),
    \end{array}
    \right.
\end{align}
\reqnomode

with initial data $t_0 > 0$ and $(u_0, v_0) \in \H \times \H$. In the following proposition, we collect the main properties of $G$ and point out that statement $iii)$, which will play a crucial role in the existence result, requires $\H$ to be finite dimensional. Its proof can be done in the lines of the proof of \cite[Proposition 3.1]{Sonntag2024b}.

\begin{myprop}
\label{prop:set_valued_G}
    The set-valued map $G$ has the following properties:
    \begin{enumerate}[i)]
        \item For all $(t,u,v) \in [t_0, +\infty) \times \H \times \H$, the set $G(t,u,v) \subseteq \H \times \H$ is convex, compact and nonempty.
        \item $G$ is upper semicontinuous.
        \item  If $\H$ is finite dimensional, then the map
        \begin{align*}
            \phi:[t_0, + \infty) \times \H \times \H \to \H\times \H, \quad (t,u,v) \mapsto \proj_{G(t,u,v)}(0)
        \end{align*}
        is locally compact.
        \item If the gradients $\nabla f_i$ are Lipschitz continuous for $i=1, \dots, m$, then there exists $c > 0$ such that for all $(t,u,v) \in [t_0, + \infty) \times \H \times \H \to \H$ it holds
        \begin{align*}
            \sup_{\xi \in G(t,u,v)} \lVert \xi \rVert_{\H \times \H} \le c\left(1 + \lVert (u,v) \rVert_{\H \times \H}\right).
        \end{align*}
    \end{enumerate}
\end{myprop}

The following theorem from \cite{Aubin2012} gives a criterion for the existence of solutions of the differential inclusion \eqref{eq:DI} on compact intervals.

\begin{theorem}
\label{thm:existence_aubin}
    Let $\mathcal{X}$ be a real Hilbert space and let $\Omega \subset \R \times \mathcal{X}$ be an open set containing $(t_0, x_0)$. Let $G : \Omega \rightrightarrows \mathcal{X}$ be an upper semicontinuous set-valued map which takes as values nonempty, closed and convex subsets of $\mathcal{X}$. Assume that the map $(t,x) \mapsto \proj_{G(t,x)}(0)$ is locally compact. Then, there exists $T > t_0$ and an absolutely continuous function $x(\cdot)$ defined on $[t_0, T]$ which is a solution of the differential inclusion
    \begin{align*}
        \dot{x}(t) \in G(t,x(t)) \ \forall t \in [t_0, T], \quad x(t_0) = x_0.
    \end{align*}
\end{theorem}

Building on Theorem \ref{thm:existence_aubin}, we can formulate the following existence result for \eqref{eq:DI}, which can be proven similar to \cite[Theorem 3.4]{Sonntag2024b}.

\begin{theorem}
\label{thm:DI_sol_exist_fin}
    Assume $\H$ is finite dimensional. Then, for all $(u_0, v_0) \in \H \times \H$ there exists $T > t_0$ and an absolutely continuous function $(u, v)$ defined on $[t_0, T]$ which is a solution of the differential inclusion \eqref{eq:DI} on $[t_0, T]$. 
\end{theorem}

In a next step we extend the solutions of \eqref{eq:DI} to $[t_0, +\infty)$ by using a standard argument that relies on Zorn's Lemma. The proof is a refinement of the one given for \cite[Theorem 3.5]{Sonntag2024b}.

\begin{theorem}
\label{thm:DI_sol_global}
    Assume $\H$ is finite dimensional. Then, for all $(u_0, v_0) \in \H \times \H$ there exists a function $(u, v)$ defined on $[t_0, +\infty)$ which is absolutely continuous on $[t_0, T]$ for all $T > t_0$ and is a solution to the differential inclusion \eqref{eq:DI}. 
\end{theorem}

\begin{proof}
We define the following set
\begin{align*}
    \mathfrak{S} \coloneqq \big \lbrace (u, v, T): & \ T \in (t_0, +\infty] \ \mbox{and} \ (u, v) : [t_0, T) \rightarrow \H \times \H \text{ is absolutely continuous on every}\\ 
    & \text{ compact interval contained in } [t_0, T) \text{ and is a solution of \eqref{eq:DI} on } [t_0, T) \big \rbrace.
\end{align*}
Note that the condition $T \in (t_0, +\infty]$ allows for the value $+\infty$ for $T$. By Theorem \ref{thm:DI_sol_exist_fin}, the set $\mathfrak{S}$ is not empty. On $\mathfrak{S}$ we define the partial order $\preccurlyeq$ as follows: for $(u_1, v_1, T_1), (u_2, v_2, T_2) \in \mathfrak{S}$,
\begin{align*}
    (u_1, v_1, T_1) \preccurlyeq (u_2, v_2, T_2) \quad \Longleftrightarrow \quad T_1 \le T_2 \text{ and } (u_1(t), v_1(t)) = (u_2(t), v_2(t)) \,\text{ for all }\, t \in [t_0, T_1).
\end{align*}
The partial order is reflexive, transitive and antisymmetric. We show that any nonempty totally ordered subset of $\mathfrak{S}$ has an upper bound in $\mathfrak{S}$. Let $\mathfrak{C} \subseteq \mathfrak{S}$ be a totally ordered nonempty subset of $\mathfrak{S}$. We define
\begin{align*}
    T_{\mathfrak{C}} \coloneqq \sup\left\lbrace T :  (u, v, T) \in \mathfrak{C} \right\rbrace
\end{align*}
and
\begin{align*}
    &(u_{\mathfrak{C}}, v_{\mathfrak{C}}) : [t_0, T_{\mathfrak{C}}) \to \H \times \H,\, (u_{\mathfrak{C}}, v_{\mathfrak{C}})(t) := (u(t), v(t)) \text{ for } t < T_{\mathfrak{C}} \ \text{and} \ (u, v, t) \in \mathfrak{C}.
\end{align*}
By construction, $(u_{\mathfrak{C}}, v_{\mathfrak{C}}, T_{\mathfrak{C}}) \in \mathfrak{S}$ and $(u, v, T) \preccurlyeq (u_{\mathfrak{C}}, v_{\mathfrak{C}}, T_{\mathfrak{C}})$, hence there exists an upper bound of $\mathfrak{C}$ in $\mathfrak{S}$. 

According to Zorn's Lemma, there exists a maximal element in $\mathfrak{S}$, which we denote by $(u, v ,T)$. If $T = + \infty$, the proof is complete.  

We assume that $T < + \infty$. We show that this contradicts the maximality of $(u, v ,T)$ in $\mathfrak{S}$. We define on $[t_0, T)$ the function
\begin{align*}
    h(t) \coloneqq \left\lVert (u(t), v(t)) - (u(t_0), v(t_0)) \right\rVert_{\H \times \H}.
\end{align*}
Using the Cauchy-Schwarz inequality, we get for almost all $t \in [t_0, T)$
\reqnomode
\begin{align}
\label{eq:derivative_h_sq}
\begin{split}
    \frac{d}{dt} \left(\frac{1}{2} h^2(t) \right) &= \left\langle (\dot{u}(t), \dot{v}(t)), (u(t), v(t)) - (u(t_0), v(t_0)) \right\rangle_{\H \times \H} \le \left\lVert (\dot{u}(t), \dot{v}(t)) \right\rVert_{\H \times \H} h(t).
\end{split}
\end{align}
Proposition \ref{prop:set_valued_G} \emph{(iii)} guarantees the existence of a constant $c > 0$ with
\begin{align}
\label{eq:bound_h_1}
    \lVert (\dot{u}(t), \dot{v}(t)) \rVert_{\H \times \H} \le c(1 + \lVert (u(t), v(t)) \rVert_{\H \times \H}),
\end{align}
for almost all $t \in [t_0, T)$. Define $\tilde{c} \coloneqq c \left( 1 + \lVert (u(t_0), v(t_0)) \rVert_{\H \times \H} \right)$. By applying the triangle inequality, we have for almost all $t \in [t_0, T)$
\begin{align}
\label{eq:bound_h_2}
    \lVert (\dot{u}(t), \dot{v}(t)) \rVert_{\H \times \H} \le \tilde{c}\left(1 + \lVert(u(t), v(t)) - (u(t_0), v(t_0)) \rVert_{\H \times \H}\right),
\end{align}
which gives
\begin{align}
\label{eq:bound_derivative}
    \frac{d}{dt} \left(\frac{1}{2} h^2(t) \right) \le \tilde{c} \left(1 + h(t) \right) h(t).
\end{align}
Using a Gronwall-type argument (see Lemma A.4 and Lemma A.5 in \cite{Brezis1973} and Theorem 3.5 in \cite{Attouch2015}), we conclude from \eqref{eq:bound_derivative} that for all $t \in [t_0, T)$
\begin{align*}
    h(t) \le \tilde{c}T \exp(\tilde{c}T),
\end{align*}
therefore, $h$ is bounded on $[t_0, T)$. Then, $u$ and $v$ are also bounded on $[t_0, T)$ and from \eqref{eq:bound_h_1} we deduce that $\dot{u}$ and $\dot{v}$ are essentially bounded. This and the fact that $\dot{u}$ and $\dot{v}$ are absolutely continuous guarantee that
\begin{align*}
    u_T \coloneqq u_0 + \int_{t_0}^T \dot{u}(s) ds \in \H \, \text{ and } \, v_T \coloneqq v_0 + \int_{t_0}^T \dot{v}(s) ds \in \H
\end{align*}
are well-defined. Further, considering the differential inclusion
\begin{align}
\label{eq:DI_continuation}
    \left\lvert 
    \begin{array}{l}
        (\dot{u}(t), \dot{v}(t)) \in G(t, u(t), v(t))\, \text{ for } t > T, \\
        \\
        (u(T), v(T)) = (u_T, v_T),
    \end{array}
    \right.
\end{align}
and using Theorem \ref{thm:DI_sol_exist_fin}, we obtain that there exist $\delta >0$ and a solution $(\hat{u}, \hat{v}):[T, T + \delta] \to \H \times \H$ of \eqref{eq:DI_continuation} which is absolutely continuous on compact intervals of $[T, T + \delta]$.  Defining
\begin{align*}
    (u^*, v^*):[t_0, \delta) \to \H\times \H,\, t \mapsto \left\lbrace \begin{array}{cl}
        (u(t),v(t)) & \text{ for } t \in [t_0, T), \\
        (\hat{u}(t), \hat{v}(t)) & \text{ for } t \in [T, T+\delta),
    \end{array} \right.
\end{align*}
we obtain an element $(u^*, v^*, T + \delta) \in \mathfrak{S}$ with the property that $(u, v, T) \neq (u^*, v^*, T + \delta)$ and $(u, v, T) \preccurlyeq (u^*, v^*, T + \delta)$. This is a contradiction to the fact that $(u, v, T)$ is a maximal element in $\mathfrak{S}$.
\end{proof}

\subsection{Existence of trajectory solutions of \eqref{eq:MTRIGS}}
\label{subsec:ap_exist_sol_CP}

In this subsection, we construct trajectory solutions of \eqref{eq:MTRIGS} starting from solutions of the differential inclusion \eqref{eq:DI}. For this purpose, we use the following well-known property of the projection, according to which, for $\H$ a real Hilbert space, $C \subseteq \H$ a nonempty, convex, and closed set, and $\eta \in \H$ a given vector, it holds
    \begin{align*}
        \xi \in \eta - \argmin_{\mu \in C} \langle \mu , \eta \rangle \,\, \text{ if and only if } \,\, \eta = \proj_{C + \xi}(0).
    \end{align*}
Using this result, one can easily see that solutions of the differential inclusions \eqref{eq:DI} lead to solutions that satisfy the equation in \eqref{eq:MTRIGS}.

\begin{theorem}
\label{thm:DI_MTRIGS}
    Let $t_0 >0$ and $x_0, v_0 \in \H$. If $(u, v) :[t_0 , \infty) \to \H \times \H$ is a solution of \eqref{eq:DI} with $(u(t_0), v(t_0)) = (x_0, v_0)$, then $x(t) \coloneqq u(t)$ satisfies the differential equation
    \begin{align*}
        \frac{\alpha}{t^q} \dot{x}(t) + \proj_{C(x(t)) + \frac{\beta}{t^p} x(t) + \ddot{x}(t)}(0) = 0,
    \end{align*}
    for almost all $t \in [t_0, +\infty)$, and $x(t_0) = x_0$, and $\dot{x}(t_0) = v_0$. 
\end{theorem}

We are now in a position to prove the existence of a trajectory solution of \eqref{eq:MTRIGS} in the sense of Definition \ref{def:sol_CP}. The following result is obtained by combing Theorem \ref{thm:DI_sol_global} and Theorem \ref{thm:DI_MTRIGS}. The fact that $x \in C^1([t_0, +\infty))$ is a consequence of the fact that $x(t) = u(t) = u(t_0) + \int_{t_0}^t v(s)ds$ for all $t \geq t_0$ and of the continuity of $v$.

\begin{theorem}
    \label{thm:main_existence_theorem_appendix}
    Assume $\H$ is finite dimensional. Then, for all $x_0, v_0 \in \H$, there exists a function $x:[t_0, +\infty) \to \H$ which is a solution of \eqref{eq:MTRIGS} in the sense of Definition \ref{def:sol_CP}.
\end{theorem}


\section{Computational details for Example \ref{ex:counter_example_cont_proj}}\label{appendix_c}

The gradient of $g(\cdot)$ is given by
    \begin{align*}
        \nabla g: \R^2 \to \R^2, \quad  x \mapsto \left\lbrace\begin{array}{lll}
        x, & \text{ if }\, \lvert x_1 \rvert \le 1,& x_2 + 1 \le \sqrt{1 - x_1^2},\\[2pt] 
       \left[ \begin{array}{c}
                \frac{x_1}{\lvert x_1 \rvert}\\
                x_2         
       \end{array}\right], & \text{ if }\, \lvert x_1 \rvert  > 1,& x_2 + 1 \le 0,\\[2pt]
       \left[ \begin{array}{c}
                \frac{x_1}{\sqrt{x_1^2 + (x_2 + 1)^2}}\\
                \frac{x_2 + 1}{\sqrt{x_1^2 + (x_2 + 1)^2}} - 1        
       \end{array}\right],& \text{ else. }
       \end{array}\right. 
    \end{align*}
            \begin{figure}
    \centering
    \begin{tikzpicture}[scale=1.2]
        
        \draw[black, thick, -stealth] (-3,0) -- (3,0);
        \draw[black, thick, -stealth] (0,-3) -- (0,2);

        \draw[black, thick, -] (0,1) -- (-.1,1);
        \draw[black, thick, -] (0,-1) -- (-.1,-1);
        \draw[black, thick, -] (0,-2) -- (-.1,-2);
        \draw[black, thick, -] (0.0125,-3) -- (-.1,-3);

        \draw[black, thick, -] (1,0) -- (1,-.1);
        \draw[black, thick, -] (2,0) -- (2,-.1);
        \draw[black, thick, -] (-1,0) -- (-1,-.1);
        \draw[black, thick, -] (-2,0) -- (-2,-.1);
        \draw[black, thick, -] (-3,0.0125) -- (-3,-.1);

        \begin{scope}
            \clip (-1.5,-1) rectangle (1.5 ,1);
            \draw (0,-1) [black, thick, dashed] circle(1);
        \end{scope}

        \draw[black, thick, dashed] (-1,-3) -- (-1,-1);
        \draw[black, thick, dashed] (1,-3) -- (1,-1);
        \draw[black, thick, dashed] (-3,-1) -- (-.985,-1);
        \draw[black, thick, dashed] (3,-1) -- (.985,-1);

        \node at (2.5, -.33) {$x_1$};
        \node at (-.33, 1.5) {$x_2$};

        \node at (-.5, -2) {$M_1$};
        \node at (-2, -2) {$M_2$};
        \node at (2, -2) {$M_2$};
        \node at (-2, 1) {$M_3$};

        \node at (2, 1) {$\R^2$};
        
    \end{tikzpicture}    
    \caption{The sets $M_i \subseteq \R^2$ for $i = 1,2,3$.}
    \label{fig:Sets_M_i}
\end{figure}
Denoting
    \begin{align*}
        M_1 \coloneqq \left\lbrace x \in \R^2 : \lvert x_1 \rvert \le 1,\, x_2 + 1 \le \sqrt{1 - x_1^2} \right\rbrace,
        M_2 \coloneqq \left\lbrace x \in \R^2 : \lvert x_1 \rvert  > 1, x_2 + 1 \le 0 \right\rbrace, 
        M_3 \coloneqq \R^2 \setminus \left( M_1 \cup M_2 \right),
    \end{align*}
we see that $\nabla g(\cdot)$ is Lipschitz continuous on $\closure(M_i)$ for $i = 1,2,3$. Since $\nabla g\big\lvert_{\closure(M_i)}(\cdot)$ and $\nabla g\big\lvert_{\closure(M_j)}(\cdot)$ coincide on $\closure(M_i) \cap \closure(M_j)$ for $i \neq j  \in \{1,2,3\}$, the Lipschitz continuity of $\nabla g(\cdot)$ follows. In fact, $\nabla g(\cdot) = \proj_{M_1}(\cdot)$, hence the Lipschitz constant of the gradient is $1$. In the following, we show that for $t \ge t_0$
    \begin{align}
    \label{eq:reg_path_z(t)3}
        z(t) = \left[\begin{array}{c}
        - (\omega(t) + 1)\sqrt{\left(\frac{t^p}{ t^p - \beta \omega(t)}\right)^2 - 1}\\
        \omega(t)   
        \end{array}\right] \in \argmin_{z \in \R^2} \max \left( f_1(z) - q_1(t), f_2(z) - q_2(t) \right) + \frac{\beta}{2 t^p}\lVert z \rVert^2.
    \end{align}
For all $t \ge t_0$, the function
    \begin{align*}
        \Phi_t : \R^2 \to \R, \quad z \mapsto \max \left( f_1(z) - q_1(t), f_2(z) - q_2(t) \right) + \frac{\beta}{2 t^p}\lVert z \rVert^2,
    \end{align*}
is strongly convex and therefore has a unique minimizer. We show that 
    \begin{align}
    \label{eq:partial_Phi_z0}
        0 \in \partial_z \Phi_t(z(t)),
    \end{align}
    where $\partial_z \Phi_t(z(t))$ denotes the convex subdifferential of $\Phi_t(\cdot)$ evaluated at $z(t)$. Note that $z_2(t) \in [2.25,2.75]$ for all $t \ge t_0$ and hence
    \begin{align*}
        \Phi_t(z) = \frac{1}{2} z_1^2 + \frac{1}{2} + g(z) + \frac{\beta}{2t^p} \lVert z \rVert^2 + \max \left( -z_1 - q_1(t), z_1 \right),
    \end{align*}
    on an open neighborhood of $z(t)$. We have
    \begin{align*}
        \partial_z \Phi_t(z(t)) = \left[ \begin{array}{c}
            z_1(t) + \frac{z_1(t)}{\sqrt{z_1(t)^2 + (z_2(t) + 1)^2}} + \frac{\beta}{t^p}z_1(t)  \\
            \frac{z_2(t) + 1}{\sqrt{z_1(t)^2 + (z_2(t) + 1)^2}} - 1 + \frac{\beta}{t^p}z_2(t) 
        \end{array} \right] + \partial_z \max\left( -z_1(t) - q_1(t), z_1(t) \right).
    \end{align*}
     Since $z_1(t) = - \frac{1}{2} q_1(t)$ we have $\partial_z \max\left( -z_1(t) - q_1(t), z_1(t) \right) = [-1,1] \times \{ 0 \}$ and hence
    \begin{align}
    \label{eq:partial_Phi_z1}
         \partial_z \Phi_t(z(t)) = \left[ \begin{array}{c}
            z_1(t) + \frac{z_1(t)}{\sqrt{z_1(t)^2 + (z_2(t) + 1)^2}} + \frac{\beta}{t^p}z_1(t)  \\
            \frac{z_2(t) + 1}{\sqrt{z_1(t)^2 + (z_2(t) + 1)^2}} - 1 + \frac{\beta}{t^p}z_2(t) 
        \end{array} \right] + [-1,1] \times \{ 0 \}.
    \end{align}
For all $t \ge t_0 = \left( 192 \beta \right)^{\frac{1}{p}}$, taking into account the definition of $z_1(t)$ and $z_2(t) \in [2.25,2.75]$, it holds
    \begin{align*}
        z_1(t) + \frac{z_1(t)}{\sqrt{z_1(t)^2 + (z_2(t) + 1)^2}} + \frac{\beta}{t^p}z_1(t) \in [-1, 1].
    \end{align*}
On the other hand, since
$$z_1(t)= -(z_2(t)+ 1)\sqrt{\left(\frac{t^p}{ t^p - \beta z_2(t)}\right)^2 - 1},$$
we have
$$\frac{z_2(t)+ 1}{\sqrt{z_1(t)^2 + (z_2(t)+ 1)^2}} = 1 - \frac{\beta}{t^p}z_2(t),$$
which proves that \eqref{eq:partial_Phi_z1}, and therefore \eqref{eq:reg_path_z(t)3} are satisfied.

\bibliographystyle{siamplain}
{\footnotesize
\bibliography{literature}
}
\vspace{30mm}



\end{document}